\documentclass[11pt,a4paper]{article}
\date{}
\makeatletter
\@addtoreset{equation}{section}
\makeatother
\usepackage{float}
\usepackage{titlesec}
\usepackage[colorlinks,
linkcolor=blue,
anchorcolor=blue,
citecolor=blue]{hyperref}
\usepackage{graphicx}
\usepackage{booktabs}
\usepackage{hyperref}
\usepackage{tocloft}
\hypersetup{hypertex=true,
	colorlinks=true,
	linkcolor=blue,
	anchorcolor=blue,
	citecolor=blue}
\usepackage[nottoc,notlot,notlof]{tocbibind}
\usepackage{hyperref,amsmath,amssymb,amscd,booktabs,arydshln}
\usepackage{amssymb,amsmath,enumerate}
\usepackage{indentfirst}
\usepackage{url}
\usepackage[capitalise]{cleveref}
\usepackage{amsthm}
\usepackage{cite}
\newtheorem{theorem}{Theorem}[section]

\newtheorem{lemma}[theorem]{Lemma}
\newtheorem{proposition}[theorem]{Proposition}
\theoremstyle{definition}

\newtheorem{remark}[theorem]{\bf Remark}

\newcommand{\supp}{\operatorname{supp}}
\allowdisplaybreaks[4]

\titleformat{\section}
{\normalfont\large\bfseries}
{\thesection}
{0.3em}
{}
\begin{document}
	\title{
		{Sixth-Order Q-Curvature of Conformally Euclidean Metrics: Improved Dimensional Thresholds, Quantitative Lower Bounds and Applications}}
	\author{\Large Jianlun Liu, Hong-Rui Sun and Ziheng Zhang
	}

	\maketitle
	\baselineskip 17pt
		\vspace{-0.7cm}
	{\bf Abstract}: This paper studies the sign of the sixth-order \(Q\)-curvature of conformally Euclidean metrics on \(\mathbb{R}^n\) for integers \(n\ge2m\) and \(m\ge4\). Under suitable assumptions, we refine the sixfold integral estimates and prove that \(Q_g^{(6)}>0\) whenever \(2m\le n\le4m-6+[\widehat h_m]\), where \(\widehat h_m\) is the unique positive root of \(G_m\) and \(\widehat h_m\sim\frac{2\sqrt7-4}{3}m\) as \(m\to\infty\). Moreover, we strengthen this positivity result by establishing the quantitative lower bound \(Q_g^{(6)}(x)\ge\Theta_{m,n}\bigl(Q_g^{(2)}(x)\bigr)^3>0\) for every \(x\in\mathbb R^n\), where \(\Theta_{m,n}>0\) is an explicit constant depending only on \(m\) and \(n\). These quantitative inequalities also yield explicit lower bounds for intermediate Laplacians in higher-order elliptic equations. On the non-positivity side, we construct metrics satisfying \(Q_g^{(6)}(0)<0\) whenever \(n>\widetilde N_m\), where \(\widetilde N_m\sim\widetilde\Lambda m\) as \(m\to\infty\) and \(\widetilde\Lambda\simeq9.74268\). These metrics extend smoothly to \(\mathbb S^n\) and are conformal to the round metric. In the same dimensional range, we further construct metrics with \(Q_g^{(6)}(0)>0\) and establish uniform positive lower and upper bounds for \(Q_g^{(2m)}\) on \(\mathbb R^n\). Our results extend the positivity range and lower the counterexample threshold obtained by V\'etois and Zeitler (arXiv preprint, 2026). The sixth-order estimate also answers a question raised in our
	previous work on quantitative fourth-order \(Q\)-curvature bounds. These results provide new insights into Conjecture 1 of Li and Xu (J. Funct. Anal. 289 (2025), 111011) in the case \(k=3\), clarifying its validity and limitations across different dimensions.

	
	\vspace{-1cm}
	\tableofcontents
	
	\section{Introduction and Statement of the Main Results}
	Let \(\phi\) be a smooth positive function in \(\mathbb{R}^n\), and consider the conformally Euclidean metric \(g=\phi^2|dx|^2\). Assume that \(m,n\in\mathbb N\) satisfy \(n\ge 2m\). Before turning to our main results for \(m\ge 4\), we recall some basic facts about \(Q\)-curvature. As is well known, \(Q\)-curvature arises naturally in the study of conformal changes of metrics, extends scalar curvature to higher orders, and is closely related to the GJMS operators. In the non-critical case \(n>2m\), the \(Q\)-curvature of order \(2m\) is given by
	\[
	Q_g^{(2m)}=\frac{2}{n-2m}\phi^{-\frac{n+2m}{2}}(-\Delta)^m\left(\phi^{\frac{n-2m}{2}}\right),
	\]
	whereas in the critical case \(n=2m\), it takes the form
	\[
	Q_g^{(2m)}=\phi^{-n}(-\Delta)^m(\ln\phi).
	\]
	Here \(\Delta=\sum_{j=1}^n\partial_{x_j}^2\) is the Euclidean Laplacian. These curvatures were introduced by Branson and {\O}rsted \cite{Branson-Orsted1991} in the case \(n=2m=4\) and generalized to general dimensions by Branson \cite{Branson1995}. They are natural higher-order analogues of scalar curvature.
	Their construction is based on Fefferman and Graham's ambient metric \cite{Fefferman-Graham1985,Fefferman-Graham2012} and the GJMS operators introduced by Graham et al. \cite{GJMS-1992}.
	These operators form a family of conformally covariant differential operators whose leading terms are powers of the Laplacian. Further background on \(Q\)-curvature and conformally covariant operators can be found in \cite{Branson1995,GJMS-1992,Case-Gover2026} and the references therein.
	
	Within this framework, a classical problem is the study of conformal metrics with constant \(Q\)-curvature. In the second-order case, the corresponding equation is closely related to the Yamabe equation. Important symmetry and classification results were obtained by Obata \cite{Obata-1971} and Caffarelli et al. \cite{Caffarelli-Gidas-Spruck1989}. Higher-order analogues have also been studied extensively. Lin \cite{Lin1998} classified positive solutions of a conformally invariant fourth-order equation on \(\mathbb{R}^n\). Wei and Xu \cite{Wei-Xu1999} subsequently extended this classification to higher-order conformally invariant equations. In the critical dimension, Martinazzi \cite{Martinazzi2009} obtained a classification for the higher-order Liouville equation on \(\mathbb{R}^{2m}\) under a finite total \(Q\)-curvature assumption. These results constitute an important part of the classification theory for constant \(Q\)-curvature metrics. 
	
	A different question concerns the signs of the lower-order \(Q\)-curvatures. More precisely, one may ask whether positivity of \(Q_g^{(2m)}\) forces the lower-order quantities \(Q_g^{(2k)}\) to remain positive. This question already plays an important role in fourth-order conformal geometry. For closed manifolds of dimension at least five, Gursky and Malchiodi \cite{Gursky-Malchiodi2015} established a strong maximum principle for the Paneitz-Branson operator under suitable positivity assumptions on the scalar curvature and the fourth-order \(Q\)-curvature. Their result also led to positivity of the Paneitz-Branson operator and its Green function. These properties were then used in the constant fourth-order \(Q\)-curvature problem. Related results were obtained by Hang and Yang \cite{Hang-Yang2015,Hang-Yang2016} through a different approach. Such positivity properties have consequences beyond the pointwise sign of the curvature itself. They have been used in existence, compactness, and uniqueness problems for prescribed \(Q\)-curvature equations. Higher-order existence results were obtained, for example, by Mazumdar and V\'etois \cite{Mazumdar-Vetois2024}. Compactness and noncompactness phenomena have been studied by Gong et al. \cite{Gong-Kim-Wei2025} and by Mazumdar and Premoselli \cite{Mazumdar-Premoselli2025}. Positivity identities also enter rigidity and uniqueness arguments. Related results can be found in Case \cite{Case-2024}, V\'etois \cite{Vetois2024} and Li and Wei \cite{Li-Wei2025}. Thus, lower-order curvature positivity is closely connected with the analytic structure of the conformally invariant equations. 
	
	The same question has a natural interpretation from the viewpoint of nonlinear elliptic equations on \(\mathbb R^n\). In the non-critical case \(n>2m\), setting \(u=\phi^{\frac{n-2m}{2}}\) transforms the constant-curvature condition \(Q_g^{(2m)}\equiv q\) into the higher-order Lane–Emden equation
	\[
	(-\Delta)^m u=\frac{n-2m}{2}q\,u^{\frac{n+2m}{n-2m}}.
	\]
	In the critical case \(n=2m\), the same condition takes the form of the higher-order Liouville equation with \(w=\ln\phi\):
	\[
	(-\Delta)^m w=qe^{nw}.
	\]
	Since a maximum principle is generally unavailable for higher-order equations, it is natural to ask whether positivity of a solution can still force positivity of its intermediate Laplacians. For the fourth-order H\'enon equation, Fazly et al. \cite{Fazly-Wei-Xu2015} established a pointwise estimate for bounded positive solutions under suitable assumptions. Their estimate implies positivity of the scalar curvature of the associated conformal metric. For the following higher-order Hardy-H\'enon equation
	\[
	(-\Delta)^m u=|x|^\sigma u^p,
	\]
	Ng\^o and Ye \cite{Ngo-Ye2022} established integral representations and super-polyharmonic properties under suitable assumptions. These results indicate that potential-theoretic representations provide an effective substitute for a higher-order maximum principle. This approach was developed systematically by Li and Xu \cite{Li-Xu2025}. In the non-critical case \(n>2m\), they assumed that \(Q_g^{(2m)}\) is positive and satisfies a polynomial lower
	bound at infinity. Under these assumptions, they proved \(Q_g^{(2)}>0\) and 
	\(Q_g^{(4)}>0\).
	Their proof is based on an integral representation of the conformal factor. It should be noted that the slow decay condition provides the necessary estimates to establish this representation and justify differentiation under the integral sign. They further conjectured that the same positivity phenomenon should persist throughout the lower-order hierarchy:
	\[
	Q_g^{(2k)}>0,
	\qquad
	1\le k\le m-1.
	\]
	The corresponding critical case \(n=2m\) was subsequently studied by Li et al. \cite{Li-Wei-Xu2025}. Under suitable assumptions, they obtained analogous positivity results for \(Q_g^{(2)}\) and \(Q_g^{(4)}\). These qualitative positivity results were strengthened by Liu et al. \cite{Liu-Lyu-Sun-Zhang2026},
	who established an explicit quantitative comparison between
	the fourth-order and second-order \(Q\)-curvatures.
	More precisely, for \(m\geq3\) and \(n\geq2m\), if \(g\) is a
	smooth conformally Euclidean metric on \(\mathbb R^n\)
	satisfying assumption (A) stated below and
	\(Q_g^{(2m)}\geq0\) in \(\mathbb R^n\), then
	\[
	Q_g^{(4)}(x)\geq
	\widetilde{\Theta}_{m,n}\bigl(Q_g^{(2)}(x)\bigr)^2>0
	\quad\text{for every }x\in\mathbb R^n,
	\]
	where
	\[
	\widetilde{\Theta}_{m,n}:=
	\begin{cases}
		\displaystyle
		\frac{(m-2)^3}{(m-1)^2},
		& n=2m,\\[2mm]
		\displaystyle
		\frac{(m-2)\big[(m-2)n+4(m-1)\big]}{2(m-1)^2},
		& n>2m.
	\end{cases}
	\]
	Thus, fourth-order positivity admits a quantitative refinement
	throughout the range \(n\geq2m\).
	
	The sixth-order curvature is the first case in which this picture changes. V\'etois and Zeitler \cite{Vetois-Zeitler-2026} investigated the positivity of \(Q_g^{(6)}\) for \(m\ge4\). Their analysis is carried out under the following assumption:
	
	(A) If $n>2m$, then we assume that there exist $s\in(-2m,0)$ and $ c,R\in(0,+\infty)$  such that
	$Q_g^{(2m)}\geq c|x|^s$ for all $x\in\mathbb{R}^n\setminus B(0,R)$. If $n=2m$, then we assume that there exists $R\in(0,+\infty)$ such that $Q_g^{(2)}\geq0$ and $Q_g^{(n)}\ge 0$ in
	$\mathbb{R}^n\setminus B(0,R)$. If $n=2m$, then we assume moreover that $Q_g^{(n)}\not\equiv0$ in $\mathbb{R}^n$.
	
	Under this assumption, V\'etois and Zeitler proved that \(Q_g^{(2m)}\ge0\)
	implies \(Q_g^{(6)}>0\) whenever \(2m\le n\le4m-6\). Thus, sixth-order positivity holds throughout this dimension-dependent range. The situation changes when the dimension becomes large. Vétois and Zeitler also constructed conformally Euclidean metrics for which \(Q_g^{(2m)}\) is positive everywhere while \(Q_g^{(6)}\) is negative at some point. Their construction applies for \(n>N_m\), where \(N_m\) is the largest real root of an explicit polynomial. They further proved that \(N_m\sim\Lambda m\) as \(m\to\infty\), where \(\Lambda\simeq10.55\).
	Consequently, there remains a substantial interval between the established positivity range and the dimensions covered by their counterexamples. V\'etois and Zeitler explicitly noted that the endpoint \(4m-6\) in their positivity theorem is not expected to be optimal. To study dimensions beyond \(4m-6\), it is enough to consider the non-critical case \(n>2m\). In this setting, let \(u=\phi^{\frac{n-2m}{2}}\).
	Under assumption (A) and the global nonnegativity of \(Q_g^{(2m)}\), the integral representation of Li and Xu \cite{Li-Xu2025} remains valid in the present setting. Thus,
	\[
	u(x)=\int_{\mathbb R^n}|x-y|^{2m-n}\,d\mu_g^{(2m)}(y),
	\]
	where \(d\mu_g^{(2m)}\) is the non-negative measure defined explicitly at the beginning of Section \ref{sec:2}. Differentiation under the integral sign then allows the lower-order \(Q\)-curvatures to be expressed in terms of multiple integrals. In particular, to derive an integral expression for \(Q_g^{(6)}\), V\'etois and Zeitler introduced
	\[
	Y_{i,j}:=\frac{x-y_i}{|x-y_i|^2}-\frac{x-y_j}{|x-y_j|^2}
	\]
	and
	\[
	Z_{i,j}:=|Y_{i,j}|^2=\frac{|y_i-y_j|^2}{|x-y_i|^2|x-y_j|^2}.
	\]
	They then obtained the sixfold integral formula
	\begin{equation}\label{eqn:VZformula}
		\begin{aligned}
			\widetilde Q_g^{(6)}(x)&:=u(x)^{6+\frac{12}{n-2m}}Q_g^{(6)}(x)
			=2(m-3)	\int_{\mathbb{R}^{6n}}\prod_{i=1}^{6}|x-y_i|^{2m-n}P_{n,m}(Z(x,y))\,d\mu_g^{(2m)}(y),
		\end{aligned}
	\end{equation}
	where \(P_{n,m}\) and \(Z(x,y)\) are given in \cite[Proposition 2.1]{Vetois-Zeitler-2026}. Their positivity problem is therefore reduced to controlling the negative contributions of this polynomial. The degree of this polynomial is one of the main differences between the sixth-order problem and the lower-order cases. The corresponding structures for \(Q_g^{(2)}\) and \(Q_g^{(4)}\) have degrees one and two, respectively. For \(Q_g^{(6)}\), the relevant polynomial is cubic. Its coefficients depend on \(n\) and \(m\). V\'etois and Zeitler divided the proof into several dimensional ranges and controlled the unfavorable terms through suitable integral identities and inequalities. Their argument yields positivity up to \(n=4m-6\). When \(n>4m-6\), three coefficients of the cubic polynomial are negative, and their method does not yield positivity in this remaining range. Table \ref{tab:previous-results} summarizes the existing results and unresolved dimensional range for lower-order \(Q\)-curvature positivity.
	\begin{table}[H]
		\centering
		\small
		\setlength{\tabcolsep}{7pt}
		\renewcommand{\arraystretch}{1.4}
		\begin{tabular}{clll}
			\toprule
			Reference
			& Assumptions
			& Dimensional conditions
			& Results \\
			\midrule
			
			\cite{Li-Xu2025}
			& \(Q_g^{(2m)}>0\) and
			& \(m\ge 2,\ n>2m\)
			& \(Q_g^{(2k)}>0\), \\
			& slow decay barrier
			&
			& \(1\le k\le\min\{2,m-1\}\) \\[4pt]
			\cmidrule(lr){1-4}		
			\cite{Li-Wei-Xu2025}
			& \(Q_g^{(n)}\ge 0\) in \(\mathbb R^n\) and
			& \(m\ge 2,\ n=2m\)
			& \(Q_g^{(2k)}>0\) or \(Q_g^{(2k)}\equiv 0\), \\
			& \(Q_g^{(2)}\ge 0\) near infinity
			&
			&\(1\le k\le\min\{2,m-1\}\) \\[4pt]		
			\cmidrule(lr){1-4}
			\cite{Vetois-Zeitler-2026}
			& \(Q_g^{(2m)}\ge 0\) and (A)
			& \(m\ge 4,\ 2m\le n\le 4m-6\)
			& \(Q_g^{(6)}>0\) \\[4pt]		
			\cmidrule(lr){1-4}
			\cite{Vetois-Zeitler-2026}
			& \(Q_g^{(2m)}\ge 0\) and (A)
			& \(m\ge 4,\ n>N_m\)
			& Existence of metrics with \\
			&
			&
			& \(Q_g^{(6)}(0)<0\)  \\[4pt]		
			\cmidrule(lr){1-4}
			Remaining range
			& \(Q_g^{(2m)}\ge 0\) and (A)
			& \(m\ge 4,\ 4m-6<n\le N_m\)
			& \textbf{?} \\
			\bottomrule
		\end{tabular}
		\caption{Previous results and unresolved cases for lower-order
			\(Q\)-curvature positivity.}
		\label{tab:previous-results}
	\end{table}
	The present work begins with this remaining range. On the positivity side, we retain the integral representation and the cubic polynomial \(P_{n,m}\) introduced by \cite{Vetois-Zeitler-2026}. The main improvement comes from a more refined comparison of the sixfold integrals associated with its cubic monomials. New inequalities reveal additional relations among these terms and provide stronger control of the negative contributions. This allows the positivity range of \(Q_g^{(6)}\) to be extended beyond \(4m-6\). To state the resulting theorem, define
	\begin{equation}\label{eqn:G-m-h-define}
		G_m(h):=8(m-2)(m-1)^2-(12m^2-40m+36)h-(14m-24)h^2-3h^3,
	\end{equation}
	and let \(\widehat h_m\) denote its unique positive root. We use \([\widehat h_m]\) for the greatest integer not exceeding \(\widehat h_m\). The main theorem is as follows.
	\begin{theorem}\label{thm:main}
		Let \(m,n\in\mathbb{N}\) with \(m\geq4\) and \(2m\leq n\leq4m-6+[\widehat h_m]\). If \(g=\phi^2|dx|^2\) is a smooth conformally Euclidean metric on \(\mathbb{R}^n\) satisfying assumption (A) and \(Q_g^{(2m)}\geq0\) in \(\mathbb{R}^n\), then \(Q_g^{(6)}(x)>0\) for every \(x\in\mathbb{R}^n\). Moreover, for every \(m\geq4\), we have \(\widehat h_m>\frac{3m}{10}\)
		and \(4m-6+[\widehat h_m]\sim\frac{8+2\sqrt7}{3} m\) as \(m\to\infty\).
	\end{theorem}
	\begin{remark}
		The numerical values in Table \ref{tab:comparison-1} are obtained directly from Theorem \ref{thm:main}. For each fixed \(m\ge4\), we first compute numerically the unique positive root \(\widehat h_m\) of
		\[
		G_m(h)=8(m-2)(m-1)^2-(12m^2-40m+36)h-(14m-24)h^2-3h^3=0.
		\]
		Since \(n\) is an integer, Theorem \ref{thm:main} extends the positivity range from the upper endpoint \(4m-6\) in \cite[Theorem 1.1]{Vetois-Zeitler-2026} to
		\(4m-6+[\widehat h_m]\), where \([\widehat h_m]\) denotes the greatest integer not exceeding \(\widehat h_m\).
		For example, when \(m=4\), \(G_4(h)=144-68h-32h^2-3h^3\),
		whose unique positive root is \(\widehat h_4\simeq1.26928\).
		Hence \([\widehat h_4]=1\) and the new positivity upper endpoint is \(4m-6+[\widehat h_m]=10+1=11\).
		Thus, while \cite[Theorem 1.1]{Vetois-Zeitler-2026} guarantees positivity only up to \(n=10\), our Theorem \ref{thm:main} also covers \(n=11\). The remaining values in Table \ref{tab:comparison-1} are obtained in exactly the same way.
		\begin{table}[H]
			\centering
			\begin{tabular}{ccccc}
				\toprule
				\(m\) 
				& V-Z positivity 
				& \(\widehat h_m\) 
				& \([\widehat h_m]\) 
				& New positivity \\
				& upper endpoint 
				& approximately
				& 
				& upper endpoint\\
				& (\(4m-6\))
				& 
				& 
				& (\(4m-6+[\widehat h_m]\))\\
				\midrule
				4  & 10 & 1.26928 & 1 & 11 \\
				5  & 14 & 1.71603 & 1 & 15 \\
				6  & 18 & 2.15499 & 2 & 20 \\
				7  & 22 & 2.59066 & 2 & 24 \\
				8  & 26 & 3.02464 & 3 & 29 \\
				9  & 30 & 3.45764 & 3 & 33 \\
				10 & 34 & 3.89002 & 3 & 37 \\
				\bottomrule
			\end{tabular}
			\caption{Comparison of the positivity ranges given by \cite[Theorem 1.1]{Vetois-Zeitler-2026} and our Theorem \ref{thm:main}.}
			\label{tab:comparison-1}
		\end{table}
	\end{remark}
	
	The positivity result in Theorem \ref{thm:main} can be further
	strengthened to a quantitative comparison between the sixth-order
	and second-order \(Q\)-curvatures. In \cite[Remark 1.2]{Liu-Lyu-Sun-Zhang2026}, Liu et al.  asked whether such a comparison holds in the range
	\(2m\leq n\leq4m-6\).
	The following theorem answers this question affirmatively
	and extends the estimate to the larger range
	\(2m\leq n\leq4m-6+[\widehat h_m]\).

	\begin{theorem}\label{them:main-3}
		Let \(m,n\in\mathbb N\) with \(m\ge4\) and \(2m\le n\le4m-6+[\widehat h_m]\). Suppose that \(g\) is a smooth conformally Euclidean metric on	\(\mathbb R^n\) satisfying assumption \((A)\) and \(Q_g^{(2m)}\geq0\) in \(\mathbb R^n\). Then,
		\[
		Q_g^{(6)}(x)\ge\Theta_{m,n}\bigl(Q_g^{(2)}(x)\bigr)^3>0 \quad\mbox{for every}\quad x\in\mathbb R^n, 
		\]
		where
		\[
		\Theta_{m,n}:=
		\begin{cases}
			\displaystyle
			\frac{(m-3)^3}{4(m-1)^3},&n=2m;\\[3mm]
			\displaystyle\frac{(m-3)A^2B^2}{4(m-1)^3},&2m<n\le3m-3;\\[3mm]
			\displaystyle\frac{(m-3)\left[A^2B^2+hC(6mA+EF)\right]}{4(m-1)^3},&3m-3<n\le4m-6;\\[3mm]
			\displaystyle\frac{(m-3)\left[AG_m(h)+hCFE\right]}{4(m-1)^3},&4m-6<n\le4m-6+[\widehat h_m].
		\end{cases}
		\]
		Here \(h,A,B,C,E\), and \(F\) are defined in Section \ref{sec-3}, while \(G_m\) is given in (\ref{eqn:G-m-h-define}).
	\end{theorem}
	
	\begin{remark}
		The fourth-order estimate of Liu et al. \cite{Liu-Lyu-Sun-Zhang2026} and
		Theorem \ref{them:main-3} establish quantitative lower bounds
		for \(Q_g^{(4)}\) and \(Q_g^{(6)}\) in terms of \(\bigl(Q_g^{(2)}\bigr)^2\) and
		\(\bigl(Q_g^{(2)}\bigr)^3\), respectively. The fourth-order estimate holds for all \(n\geq2m\), whereas the sixth-order estimate is established here for
		\(2m\leq n\leq4m-6+[\widehat h_m]\). These results motivate the following question at higher orders:
		
		\textbf{Question.} Let \(m,n,k\in\mathbb N\) with \(m\geq 5\), \(n\ge 2m\) and \(4\leq k\leq m-1\). Suppose that \(g\) is a smooth conformally Euclidean metric on	\(\mathbb R^n\) satisfying assumption (A) and \(Q_g^{(2m)}\geq0\) in \(\mathbb R^n\). 
		Under suitable restrictions on \(m,n,k\), does there exist
		a constant \(C_{m,n,k}>0\) such that
		\[
		Q_g^{(2k)}(x)\geq C_{m,n,k}\bigl(Q_g^{(2)}(x)\bigr)^k
		\]
		for every \(x\in\mathbb R^n\)?
		
		A further direction is to investigate whether analogous quantitative comparisons between lower-order curvatures can be established in CR geometry. Such an extension would require suitable analogues of the Euclidean integral identities and inequalities used here, taking into account the anisotropic and noncommutative
		structure of the Heisenberg group. In the CR fractional setting, Liu et al.
		\cite{Liu-Radulescu-Sun-Zhang2026} developed a direct spectral approach to the non-degeneracy, local rigidity and Morse index of bubble solutions to the critical CR fractional Yamabe equation.	It would be interesting to explore whether this spectral
		approach, together with suitable integral representations, can contribute to quantitative curvature estimates in that setting.

	\end{remark}
	
	The quantitative curvature estimate in Theorem \ref{them:main-3} also yields differential inequalities for higher-order elliptic equations. The corresponding fourth-order inequalities were obtained in \cite[Theorem 1.3]{Liu-Lyu-Sun-Zhang2026}.
	Here we derive sixth-order analogues involving powers or exponentials of the solutions.
	We formulate these results for both the non-critical and critical cases.
	
	\begin{theorem}\label{them:main-5}
		Let \(m,n\in\mathbb N\) satisfy \(m\geq4\) and \(n\geq2m\). Let \(\Theta_{m,n}\) be the constant defined in Theorems \(\ref{them:main-3}\). Then the following statements hold:
		
		\((i)\) Non-critical case \(n>2m\). Suppose that \(u\in C^\infty(\mathbb R^n)\) is a positive solution of
		\[
		(-\Delta)^m u=f(x)u^{\frac{n+2m}{n-2m}}\quad\text{in }\mathbb R^n,
		\]
		where \(f\geq0\) in \(\mathbb R^n\). Assume further that there exist \(s\in(-2m,0)\) and \(c,R>0\) such that \(f(x)\geq c|x|^s\) for \(|x|>R\). If \(2m<n\leq 4m-6+[\widehat h_m]\), then
		\[
		\begin{aligned}
		(-\Delta)^3\left(u^{\frac{n-6}{n-2m}}\right)\geq
		\frac{4(n-6)\Theta_{m,n}}{(n-2)^3}u^{-\frac{2n}{n-2m}}
		\left((-\Delta)\left(u^{\frac{n-2}{n-2m}}\right)\right)^3>0.
		\end{aligned}
		\]
		
		\((ii)\) Critical case \(n=2m\). Suppose that \(u\in C^\infty(\mathbb R^n)\) satisfies
		\[
		(-\Delta)^m u=f(x)e^{nu}\quad\text{in }\mathbb R^n,
		\]
		where \(f\geq0\) and \(f\not\equiv0\) in \(\mathbb R^n\). Assume further that there exists \(R>0\) such that \(\Delta u+\frac{n-2}{2}|\nabla u|^2\leq0\) for \(|x|>R\). Then,
		\[
		\begin{aligned}
		(-\Delta)^3e^{\frac{n-6}{2}u}\geq\frac{4(n-6)\Theta_{m,n}}{(n-2)^3}e^{-nu}
		\left((-\Delta)e^{\frac{n-2}{2}u}\right)^3>0.
		\end{aligned}
		\]
	\end{theorem}

	We now turn to the complementary question of non-positivity. While our Theorem \ref{thm:main} extends the range in which \(Q_g^{(6)}\) remains positive, it is also natural to ask whether the dimensional threshold for the known counterexamples can be reduced. The construction in \cite[Theorem 1.2]{Vetois-Zeitler-2026} uses two concentration points with different weights. 
	By replacing this configuration with identical peaks of equal mass and height centered at mutually orthogonal unit vectors, we obtain an improved non-positivity threshold. The result is stated as follows.
	\begin{theorem}\label{thm:main-2}
		For each integer \(m\ge4\), define \(\widetilde N_m\) as the largest real root of
		\[
		\begin{aligned}
			H_m(n):={}&-7n^4+(82m-174)n^3
			+(-164m^2+660m-252)n^2\\
			&+(296m^3-1512m^2-840m+2856)n-80m^4+1776m^3-272m^2-6864m+5440.
		\end{aligned}
		\]
		For every integer \(n>\widetilde N_m\), there exists a conformally Euclidean metric \(g\) on \(\mathbb R^n\) satisfying
		\[
		Q_g^{(6)}(0)<0
		\]
		and
		\[
		C^{-1}\le Q_g^{(2m)}\le C
		\quad\text{in }\mathbb R^n
		\]
		for some constant \(C>0\). Under stereographic projection, \(g\) extends to a smooth metric on \(\mathbb S^n\) conformal to the standard metric. Moreover,
		\(\widetilde N_m>4m-6\) and \(\widetilde N_m\sim\widetilde\Lambda m\) as \(m\to\infty\),
		where \(\widetilde\Lambda\simeq 9.74268\) is the largest real root of \(7X^4-82X^3+164X^2-296X+80\).
	\end{theorem}
	
	\begin{remark}
		The improvement is already visible at \(m=4\). In this case, \cite[Theorem 1.2]{Vetois-Zeitler-2026} gives a counterexample threshold \(N_4\simeq 29.1\), whereas our Theorem \ref{thm:main-2} gives \(\widetilde N_4\simeq 19.38\). Thus, the previous construction yields counterexamples for integer dimensions \(n\ge30\), while our construction already applies for \(n\ge20\). The difference becomes larger as \(m\) increases. Table \ref{tab:comparison-2} provides a direct comparison between the thresholds in \cite[Theorem 1.2]{Vetois-Zeitler-2026} and those obtained in Theorem \ref{thm:main-2} for \(m\in\{4,5,6,7,8,9,10\}\).
		\begin{table}[H]
			\centering
			\begin{tabular}{ccccc}
				\toprule
				\(m\) & \(N_m\) & V--Z counterexample & \(\widetilde N_m\) & New counterexample\\
				& approximately & lower endpoint & approximately & lower endpoint \\
				\midrule
				4  & 29.1 & 30 & 19.38 & 20 \\
				5  & 39.7 & 40 & 28.96 & 29 \\
				6  & 50.3 & 51 & 38.58 & 39 \\
				7  & 60.9 & 61 & 48.24 & 49 \\
				8  & 71.4 & 72 & 57.92 & 58 \\
				9  & 82.0 & 82 & 67.61 & 68 \\
				10 & 92.6 & 93 & 77.31 & 78 \\
				\bottomrule
			\end{tabular}
			\caption{Comparison of the counterexample ranges given by \cite[Theorem 1.2]{Vetois-Zeitler-2026} and our Theorem \ref{thm:main-2}.}
			\label{tab:comparison-2}
		\end{table}
	\end{remark}

	\begin{remark}
		Vétois and Zeitler \cite{Vetois-Zeitler-2026} explicitly noted that they do not expect the dimension \(4m-6\) in their positivity theorem to be optimal. Our Theorem \ref{thm:main} provides a partial improvement in this direction by extending the positivity range to \(2m\le n\le 4m-6+ [\widehat h_m]\),
		where \({4m-6+[\widehat h_m]}\sim \frac{8+2\sqrt7}{3}m \simeq 4.4305m\) as 
		\(m\to\infty\).
		On the non-positivity side, \cite[Theorem 1.2]{Vetois-Zeitler-2026} gives counterexamples for \(n>N_m\), with \(N_m\sim\Lambda m\simeq 10.55m\) as \(m\to\infty\).
		Our Theorem \ref{thm:main-2} improves this threshold by giving counterexamples for \(n>\widetilde N_m\), and \(\widetilde N_m\sim\widetilde\Lambda m\simeq 9.74268m\) as \(m\to\infty\).
		Thus, our Theorems \ref{thm:main} and \ref{thm:main-2} improve the known ranges from both sides. Our Theorem \ref{thm:main} extends the positivity range to higher dimensions, while our Theorem \ref{thm:main-2} lowers the threshold for the non-positivity examples. Nevertheless, a substantial gap remains between the two regimes. Determining the precise transition between positivity and non-positivity remains an open problem. The resulting positivity and counterexample ranges, together with the remaining open region, are summarized schematically in Figure \ref{fig:range-comparison}.
		\begin{figure}[H]
			\centering
			\includegraphics[width=0.95\textwidth]{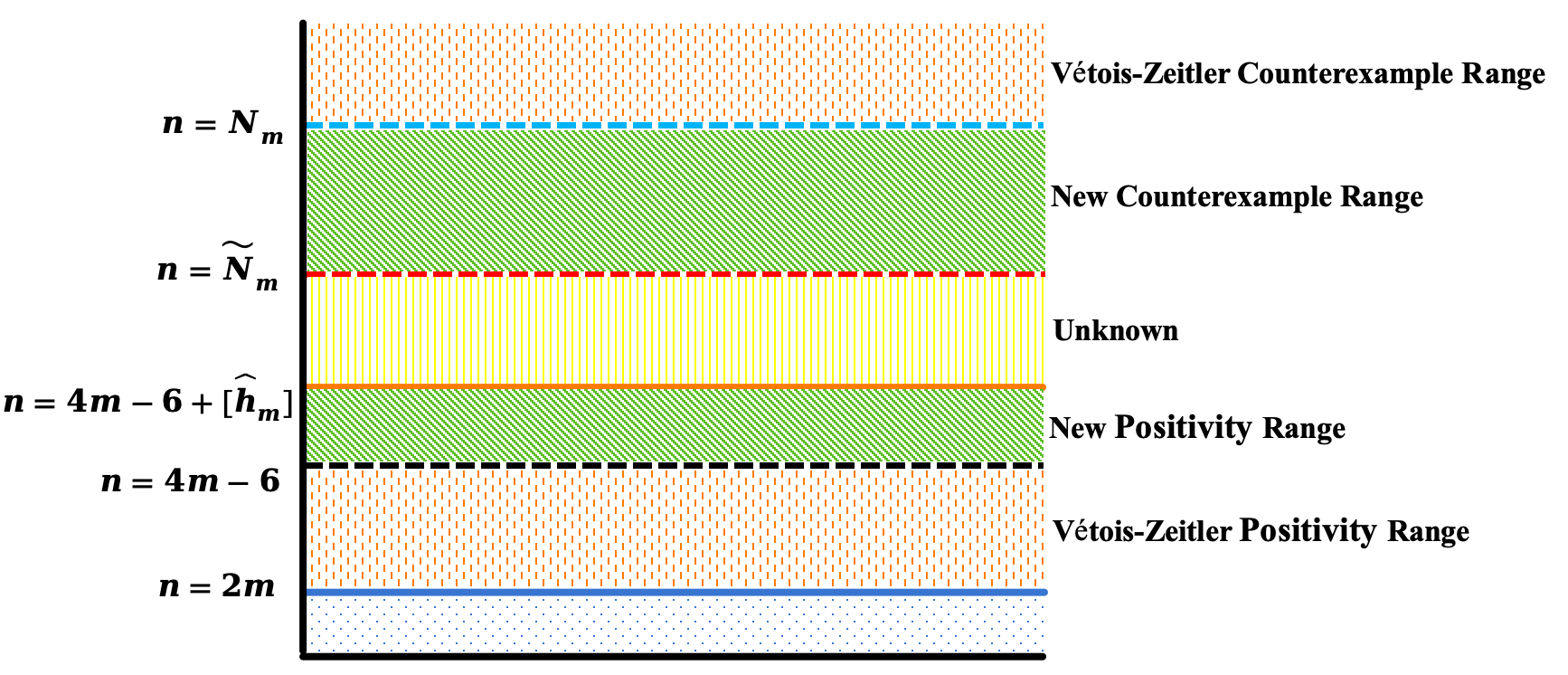}
			\caption{Schematic comparison of the positivity and counterexample ranges for \(Q_g^{(6)}\).}
			\label{fig:range-comparison}
		\end{figure}
	\end{remark}
	
	The multi-peak construction also yields positive sixth-order $Q$-curvature at the origin. Indeed, taking $L=2$, a direct calculation gives $\Phi_{m,2}(n)>0$ for $m\ge 4$ and $n>\widetilde N_m$. Therefore, we obtain a positive example in the dimensional range of Theorem \ref{thm:main-2}. The calculation follows the same argument as in the proof of Theorem \ref{thm:main-2}, with \(L=11\) replaced by \(L=2\). To avoid repetition, we omit the details.
	\begin{theorem}
		Let \(m,n\in\mathbb N\) satisfy \(m\ge4\) and \(n>\widetilde N_m\). Then there exists a smooth conformally Euclidean metric \(\tilde{g}\) on \(\mathbb R^n\) satisfying assumption (A) such that
		\[
		Q_{\tilde{g}}^{(6)}(0)>0.
		\]
		Moreover, there also exists a constant $\widetilde{C}>0$ such that
		\[
		\widetilde{C}^{-1}\le Q_{\tilde{g}}^{(2m)}\le \widetilde{C}
		\quad\text{in }\mathbb R^n.
		\]
	\end{theorem}

	\section{Refined Integral Inequalities for Q-Curvature Estimates}\label{sec:2}
	The positivity result for \(2m\le n\le4m-6\) was established in \cite[Theorem 1.1]{Vetois-Zeitler-2026}. It therefore remains to consider the range \(4m-6<n\le4m-6+[\widehat h_m]\). Since \(m\ge4\), this range lies entirely in the non-critical case \(n>2m\ge8\). Accordingly, throughout this section, we consider a smooth conformally Euclidean metric \(g=\phi^2|dx|^2\) on \(\mathbb R^n\) with \(\phi>0\). As in Theorem \ref{thm:main}, we assume that \(g\) satisfies (A) and \(Q_g^{(2m)}\ge0\) in \(\mathbb R^n\). Under these assumptions, the integral representation of Li and Xu \cite[Theorem 2.6]{Li-Xu2025} remains valid. Thus, setting \(u=\phi^{\frac{n-2m}{2}}\), we obtain that
	\[
	u(x)=\int_{\mathbb R^n}|x-y|^{2m-n}\,d\mu_g^{(2m)}(y)
	\]
	for every \(x\in\mathbb R^n\), where \(d\mu_g^{(2m)}\) is the measure defined by
	\[
	d\mu_g^{(2m)}(y)
	:=
	\frac{\Gamma\left(\frac n2-m+1\right)}
	{2^{2m-1}(m-1)!\Gamma\left(\frac n2\right)\omega_{n-1}}
	Q_g^{(2m)}(y)
	u(y)^{\frac{n+2m}{n-2m}}\,dy,
	\]
	and \(\omega_{n-1}:=|\mathbb S^{n-1}|=\frac{2\pi^\frac{n}{2}}{\Gamma(\frac{n}{2})}\) denotes the surface area of the unit sphere in \(\mathbb{R}^n\). As stated in \cite[Page 3]{Vetois-Zeitler-2026}, the integral representation of \(u\) can be differentiated under the integral sign up to order \(2m-1\). Since \(2m\ge8\), this justifies the differentiations required for the sixth-order \(Q\)-curvature formula used below. Since \(u>0\) and \(Q_g^{(2m)}\ge0\), the measure \(d\mu_g^{(2m)}\) is non-negative. For \(l\geq1\), we use \(d\mu_g^{(2m)}(y)\) as shorthand for the product measure \((\mu_g^{(2m)})^{\otimes l}\) on \((\mathbb{R}^n)^l\), with the integration variables specified in each formula. When \(l=1\), this is the original measure. For conciseness, we use the notation
	\[
	|x - y|^{[\theta_1, \dots, \theta_l]} := \prod_{i=1}^l |x - y_i|^{\theta_i}
	\]
	for all $l \in \mathbb{N}$, $\theta_1, \dots, \theta_l \in \mathbb{R}$, $x \in \mathbb{R}^n$, and $y = (y_1, \dots, y_l) \in \mathbb{R}^{ln}$. With these notations in place, we recall the following integral representation of the sixth-order \(Q\)-curvature established in \cite[Proposition 2.1]{Vetois-Zeitler-2026}.
	\begin{proposition}\label{pro:vetois-Zeitler-result-pro2.1}
		Let \(m,n\in\mathbb N\) satisfy \(n>2m\geq8\), and let \(g\) be a conformally Euclidean metric on \(\mathbb R^n\) satisfying assumption (A). Then
		\[
		\begin{aligned}
			\widetilde Q_g^{(6)}(x)=2(m-3)\Big[&
			(2m-n-2)^2(2m-n-4)^2 I_1 \\
			&-6(4m-n-6)(2m-n-2)^2(2m-n-4) I_2 \\
			&+3(4m-n-6)(3m-n-3)(2m-n-2)^2 I_3 \\
			&+4(4m-n-6)(3m-n-3)(2m-n-2)(2m-n-4) I_4 \\
			&-2(4m-n-6)(2m-n-2)^3 I_5 \\
			&+12(4m-n-6)(3m-n-3)(2m-n-2)^2 I_6 \\
			&-6(4m-n-6)(3m-n-3)(8m-3n-6)(2m-n-2) I_7 \\
			&+(4m-n-6)(3m-n-3)(5m-2n-3)(8m-3n-6) I_8
			\Big],
		\end{aligned}
		\]
		where
		\[I_1=I_1(x):=\int_{\mathbb R^{6n}}Z_{1,2}^{3}\prod_{i=1}^{6}|x-y_i|^{2m-n}d\mu_g^{(2m)}(y),
		\]
		\[
		I_2=I_2(x):=\int_{\mathbb R^{6n}}Z_{1,2}^{2}Z_{1,3}\prod_{i=1}^{6}|x-y_i|^{2m-n}d\mu_g^{(2m)}(y),
		\]
		\[
		I_3=I_3(x):=\int_{\mathbb R^{6n}}Z_{1,2}^{2}Z_{3,4}\prod_{i=1}^{6}|x-y_i|^{2m-n}d\mu_g^{(2m)}(y),
		\]
		\[
		I_4=I_4(x):=\int_{\mathbb R^{6n}}Z_{1,2}Z_{1,3}Z_{1,4}\prod_{i=1}^{6}|x-y_i|^{2m-n}d\mu_g^{(2m)}(y),
		\]
		\[
		I_5=I_5(x):=\int_{\mathbb R^{6n}}Z_{1,2}Z_{1,3}Z_{2,3}\prod_{i=1}^{6}|x-y_i|^{2m-n}d\mu_g^{(2m)}(y),
		\]
		\[
		I_6=I_6(x):=\int_{\mathbb R^{6n}}Z_{1,2}Z_{1,3}Z_{3,4}\prod_{i=1}^{6}|x-y_i|^{2m-n}d\mu_g^{(2m)}(y),
		\]
		\[
		I_7=I_7(x):=\int_{\mathbb R^{6n}}Z_{1,2}Z_{1,3}Z_{5,6}\prod_{i=1}^{6}|x-y_i|^{2m-n}d\mu_g^{(2m)}(y),
		\]
		and
		\[
		I_8=I_8(x):=\int_{\mathbb R^{6n}}Z_{1,2}Z_{3,4}Z_{5,6}\prod_{i=1}^{6}|x-y_i|^{2m-n}d\mu_g^{(2m)}(y).
		\]
	\end{proposition}
	We first verify that the integrals \(I_1,\ldots,I_8\) are finite. Using the inequality \(|a-b|^2\le 2|a|^2+2|b|^2\) and the definition of \(Z_{i,j}\), we obtain that
	\[
	Z_{i,j}\le 2\left(|x-y_i|^{-2}+|x-y_j|^{-2}\right).
	\]
	Since each integrand in \(I_1,\ldots,I_8\) contains three \(Z\)-factors, it follows from \cite[(9.1)]{Steele-2004} that
	
	\[
	\begin{aligned}
		Z_{i_1,j_1}Z_{i_2,j_2}Z_{i_3,j_3}&\le8\prod_{l=1}^{3}\left(|x-y_{i_l}|^{-2}+|x-y_{j_l}|^{-2}\right)\\
		&\le 8\prod_{l=1}^{3}\left(\sum_{\eta=1}^{6}|x-y_\eta|^{-2}\right)
		=8\left(\sum_{\eta=1}^{6}|x-y_\eta|^{-2}\cdot1\right)^3
		&\le 288\sum_{\eta=1}^{6}|x-y_\eta|^{-6},
	\end{aligned}
	\]
	for any \(i_l,j_l\in\{1,\ldots,6\}\) with \(i_l\ne j_l\) and \(l=1,2,3\). By Tonelli's theorem and the integral representation of \(u\), we have
	\[
	\begin{aligned}
		I_r(x)&\le 288\int_{\mathbb R^{6n}}\left(\sum_{\eta=1}^{6}|x-y_\eta|^{-6}\right)
		\prod_{i=1}^{6}|x-y_i|^{2m-n}\,d\mu_g^{(2m)}(y)\\
		&=288\sum_{\eta=1}^{6}\int_{\mathbb R^{6n}}|x-y_\eta|^{2m-n-6}\prod_{\substack{i=1\\i\ne \eta}}^{6}|x-y_i|^{2m-n}\,d\mu_g^{(2m)}(y)\\
		&=288\sum_{\eta=1}^{6}\left(\int_{\mathbb R^n}|x-y_\eta|^{2m-n-6}\,d\mu_g^{(2m)}(y_\eta)\right)\times u(x)^5\\
		&=1728u(x)^5\int_{\mathbb R^n}|x-y|^{2m-n-6}\,d\mu_g^{(2m)}(y)\\
	\end{aligned}
	\]
	for \(r=1,\dots,8\). Since \(m\ge4\), we have \(6<2m\). The weighted integrability argument of Li and Xu \cite[Lemma 2.2]{Li-Xu2025} also applies under assumption (A), yielding
	\[
	\int_{\mathbb R^n}|x-y|^{2m-n-6}\,d\mu_g^{(2m)}(y)<\infty.
	\]
	Moreover, \(u=\varphi^{\frac{n-2m}{2}}\), \(\varphi\in C^\infty(\mathbb R^n)\) and \(\varphi>0\) in \(\mathbb R^n\) imply that \(0<u(x)<\infty\) for every \(x\in\mathbb R^n\). Hence, \(I_r(x)<\infty\) for \(r=1,\ldots,8\).
	
	The positivity of \(u(x)\) also implies that \(Q_g^{(6)}(x)\) and \(\widetilde Q_g^{(6)}(x)\) have the same sign. Thus, it suffices to consider \(\widetilde Q_g^{(6)}\) in the subsequent proofs. To this end, we first exploit the symmetry of the sixfold integrals. Relabeling the integration variables yields
	\[
	\int_{\mathbb R^{6n}}Z_{i,j}^3 \prod_{k=1}^6 |x-y_k|^{2m-n}d\mu_g^{(2m)}(y)=I_1
	\]
	for any distinct \(i,j\in\{1,\ldots,6\}\). Furthermore, H\"older's inequality gives
	\[
	\begin{aligned}
		&\int_{\mathbb R^{6n}}Z_{i_1,j_1}Z_{i_2,j_2}Z_{i_3,j_3}\prod_{k=1}^6|x-y_k|^{2m-n}\,d\mu_g^{(2m)}(y)\\
		&\le	\prod_{l=1}^3\left(\int_{\mathbb R^{6n}}Z_{i_l,j_l}^3\prod_{k=1}^6|x-y_k|^{2m-n}\,d\mu_g^{(2m)}(y)\right)^\frac{1}{3}
		=I_1.
	\end{aligned}
	\]
	Applying the above estimate to the triples \(((i_1,j_1),(i_2,j_2),(i_3,j_3))=((1,2),(1,3),(1,4))\)
	and \(	((i_1,j_1),(i_2,j_2),(i_3,j_3))=((1,2),(1,3),(2,3))\),
	we obtain that
	\begin{equation}\label{eqn:I4-leq-I1}
		I_4=\int_{\mathbb R^{6n}}Z_{1,2}Z_{1,3}Z_{1,4}\prod_{k=1}^6 |x-y_k|^{2m-n}\,d\mu_g^{(2m)}(y)
		\le I_1
	\end{equation}
	and
	\begin{equation}\label{eqn:I5-leq-I1}
		I_5=\int_{\mathbb R^{6n}}Z_{1,2}Z_{1,3}Z_{2,3}\prod_{k=1}^6 |x-y_k|^{2m-n}\,d\mu_g^{(2m)}(y)
		\le I_1.
	\end{equation}
	In the following, we establish three inequalities that will be used in the proof of Theorem \ref{thm:main}.
	\begin{lemma}\label{lem:I7-leq-I3-I6-I1-I4}
		For any \(x\in\mathbb R^n\), we have
		\begin{itemize}
			\item[(i)] $I_7\leq I_3$;
			\item[(ii)] $I_2^2\leq I_1 I_4$;
			\item[(iii)] $I_7\leq I_6$.
		\end{itemize}
	\end{lemma}
	\begin{proof}
		$(i)$ We note that the variables \(y_5\) and \(y_6\) do not appear in \(Z_{1,2}^2Z_{3,4}\). By Fubini's theorem, we have
		\[
		\begin{aligned}
			I_3	
			&=\left(
			\int_{\mathbb R^n}|x-y|^{2m-n}\,d\mu_g^{(2m)}(y)\right)^2
			\times	\left(\int_{\mathbb R^{2n}}Z_{1,2}^2\prod_{i=1}^2|x-y_i|^{2m-n}\,d\mu_g^{(2m)}(y)\right)\\
			&\ \ \times
			\left(\int_{\mathbb R^{2n}}Z_{3,4}\prod_{i=3}^4|x-y_i|^{2m-n}\,d\mu_g^{(2m)}(y)\right)\\
			&=\left(\int_{\mathbb R^n}|x-y|^{2m-n}\,d\mu_g^{(2m)}(y)\right)^2
			\times\left(\int_{\mathbb R^{2n}}Z_{1,2}^2\prod_{i=1}^2|x-y_i|^{2m-n}\,d\mu_g^{(2m)}(y)\right)\\
			&\ \ \times\left(\int_{\mathbb R^{2n}}Z_{1,2}\prod_{i=1}^2|x-y_i|^{2m-n}\,d\mu_g^{(2m)}(y)\right).
		\end{aligned}
		\]
		Similarly, \(y_4\) does not appear in \(Z_{1,2}Z_{1,3}Z_{5,6}\), while \(Z_{5,6}\) depends only on \(y_5\) and \(y_6\). Thus,
		\[
		\begin{aligned}
			I_7&=\left(\int_{\mathbb R^n}|x-y|^{2m-n}\,d\mu_g^{(2m)}(y)\right)\times\left(\int_{\mathbb R^{3n}}Z_{1,2}Z_{1,3}\prod_{i=1}^3|x-y_i|^{2m-n}\,d\mu_g^{(2m)}(y)\right)\\
			&\ \ \times\left(\int_{\mathbb R^{2n}}Z_{5,6}\prod_{i=5}^6|x-y_i|^{2m-n}\,d\mu_g^{(2m)}(y)\right)\\
			&=\left(\int_{\mathbb R^n}|x-y|^{2m-n}\,d\mu_g^{(2m)}(y)\right)
			\times\left(\int_{\mathbb R^{3n}}Z_{1,2}Z_{1,3}\prod_{i=1}^3|x-y_i|^{2m-n}\,d\mu_g^{(2m)}(y)\right)\\
			&\ \ \times\left(\int_{\mathbb R^{2n}}Z_{1,2}\prod_{i=1}^2|x-y_i|^{2m-n}\,d\mu_g^{(2m)}(y)\right).
		\end{aligned}
		\]
		Since
		\[
		\begin{aligned}
			0&\leq\int_{\mathbb R^{3n}}(Z_{1,2}-Z_{1,3})^2\prod_{i=1}^3|x-y_i|^{2m-n}\,d\mu_g^{(2m)}(y)\\
			&=2\int_{\mathbb R^{3n}}Z_{1,2}^2\prod_{i=1}^3|x-y_i|^{2m-n}\,d\mu_g^{(2m)}(y)-2\int_{\mathbb R^{3n}}Z_{1,2}Z_{1,3}\prod_{i=1}^3|x-y_i|^{2m-n}\,d\mu_g^{(2m)}(y)\\
			&=2\left(\int_{\mathbb R^n}|x-y|^{2m-n}\,d\mu_g^{(2m)}(y)\right)\left(\int_{\mathbb R^{2n}}
			Z_{1,2}^2\prod_{i=1}^2|x-y_i|^{2m-n}\,d\mu_g^{(2m)}(y)\right)\\
			&\ \ -2\int_{\mathbb R^{3n}}Z_{1,2}Z_{1,3}\prod_{i=1}^3|x-y_i|^{2m-n}\,d\mu_g^{(2m)}(y),\\
		\end{aligned}
		\]
		we have
		\[
		\begin{aligned}
			&\int_{\mathbb R^{3n}}Z_{1,2}Z_{1,3}\prod_{i=1}^3|x-y_i|^{2m-n}\,d\mu_g^{(2m)}(y)\\
			&\leq\left(\int_{\mathbb R^n}|x-y|^{2m-n}\,d\mu_g^{(2m)}(y)\right)
			\left(\int_{\mathbb R^{2n}}Z_{1,2}^2\prod_{i=1}^2|x-y_i|^{2m-n}\,d\mu_g^{(2m)}(y)\right).
		\end{aligned}
		\]
		Combining this with the facts that
		\[
		\int_{\mathbb R^{2n}}Z_{1,2}\prod_{i=1}^2|x-y_i|^{2m-n}\,d\mu_g^{(2m)}(y)\geq0
		\]
		and
		\[
		u(x)=\int_{\mathbb R^n}|x-y|^{2m-n}\,d\mu_g^{(2m)}(y)>0,
		\]
		we deduce that \(I_7\leq I_3\).
		
		$(ii)$ For fixed \(y_1\), applying the Cauchy-Schwarz inequality with respect to \(y_2\) gives
		\[
		\begin{aligned}
			&\left(\int_{\mathbb R^n}Z_{1,2}^2|x-y_2|^{2m-n}\,d\mu_g^{(2m)}(y_2)\right)^2\\
			&=\left( \int_{\mathbb R^n} Z_{1,2}^\frac{1}{2}|x-y_2|^\frac{2m-n}{2}\cdot Z_{1,2}^\frac{3}{2} |x-y_2|^\frac{2m-n}{2}\,d\mu_g^{(2m)}(y_2) \right)^2\\
			&\le\left(\int_{\mathbb R^n}Z_{1,2}|x-y_2|^{2m-n}\,d\mu_g^{(2m)}(y_2)\right)
			\left(\int_{\mathbb R^n}Z_{1,2}^3|x-y_2|^{2m-n}\,d\mu_g^{(2m)}(y_2)\right).
		\end{aligned}
		\]
		Hence,
		\begin{equation}\label{eqn:Z122-Z12-Z12-Z123}
			\begin{aligned}
				&\left(\int_{\mathbb R^n}Z_{1,2}^2|x-y_2|^{2m-n}\,d\mu_g^{(2m)}(y_2)\right)
				\left(\int_{\mathbb R^n}Z_{1,2}|x-y_2|^{2m-n}\,d\mu_g^{(2m)}(y_2)\right)\\
				&\le\left(\int_{\mathbb R^n}Z_{1,2}|x-y_2|^{2m-n}\,d\mu_g^{(2m)}(y_2)\right)^\frac{3}{2}
				\left(\int_{\mathbb R^n}Z_{1,2}^3|x-y_2|^{2m-n}\,d\mu_g^{(2m)}(y_2)\right)^\frac{1}{2}.
			\end{aligned}
		\end{equation}
		The variables \(y_4,y_5\) and \(y_6\) do not appear in \(Z_{1,2}^2Z_{1,3}\). By Fubini's theorem, we derive that
		\[
		\begin{aligned}
			I_2&=\int_{\mathbb R^{6n}}Z_{1,2}^{2}Z_{1,3}\prod_{i=1}^{6}|x-y_i|^{2m-n}d\mu_g^{(2m)}(y)\\
			&=\left(\int_{\mathbb R^n}|x-y|^{2m-n}\,d\mu_g^{(2m)}(y)\right)^3\int_{\mathbb R^{3n}}Z_{1,2}^2Z_{1,3}\prod_{i=1}^3|x-y_i|^{2m-n}\,d\mu_g^{(2m)}(y).
		\end{aligned}
		\]
		For fixed \(y_1\), \(Z_{1,2}\) and \(Z_{1,3}\) depend separately on \(y_2\) and \(y_3\). Applying Fubini's theorem again gives
		\[
		\begin{aligned}
			I_2
			&=\left(\int_{\mathbb R^n}|x-y|^{2m-n}\,d\mu_g^{(2m)}(y)\right)^3\\
			&\quad\times\int_{\mathbb R^n}|x-y_1|^{2m-n}\left(\int_{\mathbb R^n}Z_{1,2}^2|x-y_2|^{2m-n}\,d\mu_g^{(2m)}(y_2)\right)\\
			&\quad\times\left(\int_{\mathbb R^n}Z_{1,2}|x-y_2|^{2m-n}\,d\mu_g^{(2m)}(y_2)\right)d\mu_g^{(2m)}(y_1).
		\end{aligned}
		\]
		Combining this with (\ref{eqn:Z122-Z12-Z12-Z123}) and the Cauchy-Schwarz inequality, we have
		\begin{equation}\label{eqn:I2-1/2-1/2}
			\begin{aligned}
				I_2
				&\leq\left(\int_{\mathbb R^n}|x-y|^{2m-n}\,d\mu_g^{(2m)}(y)\right)^3\\
				&\quad\times\left[\int_{\mathbb R^n}|x-y_1|^{2m-n}\left(\int_{\mathbb R^n}Z_{1,2}|x-y_2|^{2m-n}\,d\mu_g^{(2m)}(y_2)\right)^3d\mu_g^{(2m)}(y_1)\right]^\frac{1}{2}\\
				&\quad\times\left[\int_{\mathbb R^n}|x-y_1|^{2m-n}\left(\int_{\mathbb R^n}Z_{1,2}^3|x-y_2|^{2m-n}\,d\mu_g^{(2m)}(y_2)\right)d\mu_g^{(2m)}(y_1)\right]^\frac{1}{2}\\
			\end{aligned}
		\end{equation}
		By Fubini's theorem, the definitions of \(I_4\) and \(I_1\) give
		\[
		\begin{aligned}
			I_4&=\int_{\mathbb R^{6n}}Z_{1,2}Z_{1,3}Z_{1,4}\prod_{i=1}^{6}|x-y_i|^{2m-n}\,d\mu_g^{(2m)}(y)\\
			&=\left(\int_{\mathbb R^n}|x-y|^{2m-n}\,d\mu_g^{(2m)}(y)\right)^2\int_{\mathbb R^n}|x-y_1|^{2m-n}\left(\int_{\mathbb R^n}Z_{1,2}|x-y_2|^{2m-n}\,d\mu_g^{(2m)}(y_2)\right)\\
			&\quad\times\left(\int_{\mathbb R^n}Z_{1,3}|x-y_3|^{2m-n}\,d\mu_g^{(2m)}(y_3)\right)\times\left(\int_{\mathbb R^n}Z_{1,4}|x-y_4|^{2m-n}\,d\mu_g^{(2m)}(y_4)\right)d\mu_g^{(2m)}(y_1)\\
			&=\left( \int_{\mathbb R^n}\!\!\!\!\! |x-y|^{2m-n} \,d\mu_g^{(2m)}(y) \right)^2 
			\times \int_{\mathbb R^n} \!\!\!|x-y_1|^{2m-n} \left( \int_{\mathbb R^n}\!\!\! Z_{1,2}|x-y_2|^{2m-n} d\mu_g^{(2m)}(y_2) \right)^3\!\!\! d\mu_g^{(2m)}(y_1).
		\end{aligned}
		\]
		and
		\[
		\begin{aligned}
			I_1
			&=\left(\int_{\mathbb R^n}|x-y|^{2m-n}\,d\mu_g^{(2m)}(y)\right)^4\\
			&\quad\times\int_{\mathbb R^{2n}}Z_{1,2}^3|x-y_1|^{2m-n}|x-y_2|^{2m-n}\,d\mu_g^{(2m)}(y_1)d\mu_g^{(2m)}(y_2)\\
			&=\left(\int_{\mathbb R^n}|x-y|^{2m-n}\,d\mu_g^{(2m)}(y)\right)^4\\
			&\quad\times\int_{\mathbb R^n}|x-y_1|^{2m-n}
			\left(\int_{\mathbb R^n}Z_{1,2}^3|x-y_2|^{2m-n}\,d\mu_g^{(2m)}(y_2)\right)
			d\mu_g^{(2m)}(y_1).
		\end{aligned}
		\]
		Substituting these identities into (\ref{eqn:I2-1/2-1/2}) yields
		\[
		\begin{aligned}
			I_2&\leq
			\frac{I_4^\frac{1}{2}I_1^\frac{1}{2}\left(\int_{\mathbb R^n}|x-y|^{2m-n}d\mu_g^{(2m)}(y)\right)^3}{\left(\int_{\mathbb R^n}|x-y|^{2m-n}\,d\mu_g^{(2m)}(y)\right)\left(\int_{\mathbb R^n}|x-y|^{2m-n}\,d\mu_g^{(2m)}(y)\right)^2},
		\end{aligned}
		\]
		which implies that \(I_2^2\leq I_1I_4\).
	
		$(iii)$ For fixed \(x\in\mathbb R^n\), define
		\[
		V(y):=\frac{x-y}{|x-y|^2}
		\]
		and let
		\[
		\overline V(x):=\frac{\displaystyle\int_{\mathbb R^n}V(y)|x-y|^{2m-n}\,d\mu_g^{(2m)}(y)}
		{\displaystyle\int_{\mathbb R^n}|x-y|^{2m-n}\,d\mu_g^{(2m)}(y)}.
		\]
		We further set
		\[
		X(y):=V(y)-\overline V(x)
		\]
		and
		\[
		Y_{i,j}(y_i,y_j):=V(y_i)-V(y_j).
		\]
		The definition of \(\overline V(x)\) immediately yields
		\begin{equation}\label{eqn:Xy-int-equ-0}
			\begin{aligned}
				&\int_{\mathbb R^n}X(y)|x-y|^{2m-n}\,d\mu_g^{(2m)}(y)\\
				&=\int_{\mathbb R^n}V(y)|x-y|^{2m-n}\,d\mu_g^{(2m)}(y)
				-\frac{\displaystyle\int_{\mathbb R^n}V(y)|x-y|^{2m-n}\,d\mu_g^{(2m)}(y)}{\displaystyle\int_{\mathbb R^n}|x-y|^{2m-n}\,d\mu_g^{(2m)}(y)}\int_{\mathbb R^n}|x-y|^{2m-n}\,d\mu_g^{(2m)}(y)\\
				&=0.
			\end{aligned}
		\end{equation}
		Moreover, there also hold
		\[
		Y_{i,j}(y_i,y_j)=V(y_i)-V(y_j)=X(y_i)+\overline V(x)-X(y_j)-\overline V(x)=X(y_i)-X(y_j)
		\]
		and 
		\[
		Z_{i,j}=|Y_{i,j}|^2=|X(y_i)-X(y_j)|^2.
		\]
		
		We first evaluate \(I_7\). The variable \(y_4\) does not appear in \(Z_{1,2}Z_{1,3}Z_{5,6}\), while \(Z_{5,6}\) depends only on \(y_5\) and \(y_6\). Fubini's theorem gives
		\begin{equation}\label{eqn:I7-one-term}
			\begin{aligned}
				I_7&=\left(\int_{\mathbb R^n}|x-y|^{2m-n}\,d\mu_g^{(2m)}(y)\right)\\
				&\quad\times\left(\int_{\mathbb R^{3n}}Z_{1,2}Z_{1,3}\prod_{i=1}^3|x-y_i|^{2m-n}\,d\mu_g^{(2m)}(y)\right)\\
				&\quad\times\left(\int_{\mathbb R^{2n}}Z_{5,6}\prod_{i=5}^6|x-y_i|^{2m-n}\,d\mu_g^{(2m)}(y)\right).
			\end{aligned}
		\end{equation}
		Using \(Z_{5,6}=|X(y_5)-X(y_6)|^2\), the symmetry between \(y_5\) and \(y_6\) and (\ref{eqn:Xy-int-equ-0}), one has
		
		\begin{equation}\label{eqn:I7-two-term}
			\begin{aligned}
				&\int_{\mathbb R^{2n}}Z_{5,6}\prod_{i=5}^6|x-y_i|^{2m-n}\,d\mu_g^{(2m)}(y)\\
				&=2\int_{\mathbb R^{2n}} |X(y_5)|^2 |x-y_5|^{2m-n} |x-y_6|^{2m-n} \,d\mu_g^{(2m)}(y)\\
				&=2\left(\int_{\mathbb R^n}|x-y|^{2m-n}\,d\mu_g^{(2m)}(y)\right)
				\left(\int_{\mathbb R^n}|X(y)|^2|x-y|^{2m-n}\,d\mu_g^{(2m)}(y)\right).
			\end{aligned}
		\end{equation}
		Similarly, fixing \(y_1\) and using (\ref{eqn:Xy-int-equ-0}), we arrive at
		\[
		\begin{aligned}
			&\int_{\mathbb R^n}Z_{1,2}|x-y_2|^{2m-n}\,d\mu_g^{(2m)}(y_2)\\
			&=\int_{\mathbb R^n}|X(y_1)-X(y_2)|^2|x-y_2|^{2m-n}\,d\mu_g^{(2m)}(y_2)\\
			&=\left(\int_{\mathbb R^n}|x-y|^{2m-n}\,d\mu_g^{(2m)}(y)\right)|X(y_1)|^2+\int_{\mathbb R^n}|X(y)|^2|x-y|^{2m-n}\,d\mu_g^{(2m)}(y).
		\end{aligned}
		\]
		Applying the same formula to \(Z_{1,3}\) and integrating with respect to \(y_1\), we obtain
		\begin{equation}\label{eqn:I7-three-term}
			\begin{aligned}
				&\int_{\mathbb R^{3n}}Z_{1,2}Z_{1,3}\prod_{i=1}^3|x-y_i|^{2m-n}d\mu_g^{(2m)}(y)\\
				&=\left(\int_{\mathbb R^n}|x-y|^{2m-n}\,d\mu_g^{(2m)}(y)\right)^2
				\left(\int_{\mathbb R^n}|X(y)|^4|x-y|^{2m-n}\,d\mu_g^{(2m)}(y)\right)\\
				&\quad+3\left(\int_{\mathbb R^n}|x-y|^{2m-n}\,d\mu_g^{(2m)}(y)\right)
				\left(\int_{\mathbb R^n}|X(y)|^2|x-y|^{2m-n}\,d\mu_g^{(2m)}(y)\right)^2.
			\end{aligned}
		\end{equation}
		Combining (\ref{eqn:I7-one-term}), (\ref{eqn:I7-two-term}) and (\ref{eqn:I7-three-term}), one has
		\begin{equation}\label{eqn:I7-new-expression}
			\begin{aligned}
				I_7&=2\left(\int_{\mathbb R^n}|x-y|^{2m-n}\,d\mu_g^{(2m)}(y)\right)^4
				\left(\int_{\mathbb R^n}|X(y)|^2|x-y|^{2m-n}\,d\mu_g^{(2m)}(y)\right)\\
				&\quad\times\left(\int_{\mathbb R^n}|X(y)|^4|x-y|^{2m-n}\,d\mu_g^{(2m)}(y)\right)\\
				&\quad+6\left(\int_{\mathbb R^n}|x-y|^{2m-n}\,d\mu_g^{(2m)}(y)\right)^3
				\left(\int_{\mathbb R^n}|X(y)|^2|x-y|^{2m-n}\,d\mu_g^{(2m)}(y)\right)^3.
			\end{aligned}
		\end{equation}
		
		We now turn to \(I_6\). It follows from Fubini's theorem and the fact \(y_5\) and \(y_6\) do not appear in \(Z_{1,2}Z_{1,3}Z_{3,4}\) that
		\begin{equation}\label{eqn:I6-A1-A2-A3-expression}
			\begin{aligned}
				I_6&=\left(	\int_{\mathbb R^n}|x-y|^{2m-n}\,d\mu_g^{(2m)}(y)\right)^2
				\int_{\mathbb R^{4n}}Z_{1,2}Z_{1,3}Z_{3,4}\prod_{i=1}^4|x-y_i|^{2m-n}\,d\mu_g^{(2m)}(y)\\
				&=\left(\int_{\mathbb R^n}|x-y|^{2m-n}\,d\mu_g^{(2m)}(y)\right)^2\\
				&\quad\times\int_{\mathbb R^{2n}}\left[\left(\int_{\mathbb R^n}|x-y|^{2m-n}\,d\mu_g^{(2m)}(y)
				\right)|X(y_1)|^2+\int_{\mathbb R^n}|X(y)|^2|x-y|^{2m-n}\,d\mu_g^{(2m)}(y)\right]\\
				&\quad\times\left[\left(\int_{\mathbb R^n}|x-y|^{2m-n}\,d\mu_g^{(2m)}(y)\right)|X(y_3)|^2
				+\int_{\mathbb R^n}|X(y)|^2|x-y|^{2m-n}\,d\mu_g^{(2m)}(y)\right]\\
				&\quad\times|X(y_1)-X(y_3)|^2|x-y_1|^{2m-n}|x-y_3|^{2m-n}d\mu_g^{(2m)}(y_1)d\mu_g^{(2m)}(y_3)\\
				&=\left(\int_{\mathbb R^n}|x-y|^{2m-n}\,d\mu_g^{(2m)}(y)\right)^2(\mathcal A_1-\mathcal A_2+\mathcal A_3),
			\end{aligned}
		\end{equation}
		where 
		\[
		\begin{aligned}
			\mathcal A_1&:=\int_{\mathbb R^{2n}}\left[\left(\int_{\mathbb R^n}|x-y|^{2m-n}\,d\mu_g^{(2m)}(y)
			\right)|X(y_1)|^2+\int_{\mathbb R^n}|X(y)|^2|x-y|^{2m-n}\,d\mu_g^{(2m)}(y)\right]\\
			&\quad\times\left[\left(\int_{\mathbb R^n}|x-y|^{2m-n}\,d\mu_g^{(2m)}(y)\right)|X(y_3)|^2
			+\int_{\mathbb R^n}|X(y)|^2|x-y|^{2m-n}\,d\mu_g^{(2m)}(y)\right]\\
			&\quad\times|X(y_1)|^2|x-y_1|^{2m-n}|x-y_3|^{2m-n}d\mu_g^{(2m)}(y_1)d\mu_g^{(2m)}(y_3),
		\end{aligned}
		\]
		\[
		\begin{aligned}
			\mathcal A_2&:=2\int_{\mathbb R^{2n}}\left[\left(\int_{\mathbb R^n}|x-y|^{2m-n}\,d\mu_g^{(2m)}(y)
			\right)|X(y_1)|^2+\int_{\mathbb R^n}|X(y)|^2|x-y|^{2m-n}\,d\mu_g^{(2m)}(y)\right]\\
			&\quad\times\left[\left(\int_{\mathbb R^n}|x-y|^{2m-n}\,d\mu_g^{(2m)}(y)\right)|X(y_3)|^2
			+\int_{\mathbb R^n}|X(y)|^2|x-y|^{2m-n}\,d\mu_g^{(2m)}(y)\right]\\
			&\quad\times (X(y_1)\cdot X(y_3))
			|x-y_1|^{2m-n}|x-y_3|^{2m-n}d\mu_g^{(2m)}(y_1)d\mu_g^{(2m)}(y_3),
		\end{aligned}
		\]
		and
		\[
		\begin{aligned}
			\mathcal A_3&:=\int_{\mathbb R^{2n}}\left[\left(\int_{\mathbb R^n}|x-y|^{2m-n}\,d\mu_g^{(2m)}(y)\right)|X(y_1)|^2+\int_{\mathbb R^n}|X(y)|^2|x-y|^{2m-n}\,d\mu_g^{(2m)}(y)\right]\\
			&\quad\times\left[\left(\int_{\mathbb R^n}|x-y|^{2m-n}\,d\mu_g^{(2m)}(y)\right)|X(y_3)|^2
			+\int_{\mathbb R^n}|X(y)|^2|x-y|^{2m-n}\,d\mu_g^{(2m)}(y)\right]\\
			&\quad\times|X(y_3)|^2|x-y_1|^{2m-n}|x-y_3|^{2m-n}d\mu_g^{(2m)}(y_1)d\mu_g^{(2m)}(y_3).
		\end{aligned}
		\]
		We begin with \(\mathcal A_1\). The product structure of its integrand allows us to separate the \(y_1\)- and \(y_3\)-integrations. Fubini's theorem then yields
		\[
		\begin{aligned}
			\mathcal A_1
			&=\Bigg[\left(\int_{\mathbb R^n}|x-y|^{2m-n}\,d\mu_g^{(2m)}(y)\right)
			\left(\int_{\mathbb R^n}|X(y)|^4|x-y|^{2m-n}\,d\mu_g^{(2m)}(y)\right)\\
			&\quad+\left(\int_{\mathbb R^n}|X(y)|^2|x-y|^{2m-n}\,d\mu_g^{(2m)}(y)\right)^2\Bigg]\\
			&\quad\times2\left(\int_{\mathbb R^n}|x-y|^{2m-n}\,d\mu_g^{(2m)}(y)\right)
			\left(\int_{\mathbb R^n}|X(y)|^2|x-y|^{2m-n}\,d\mu_g^{(2m)}(y)\right)\\
			&=2\left(\int_{\mathbb R^n}|x-y|^{2m-n}\,d\mu_g^{(2m)}(y)\right)^2\left(\int_{\mathbb R^n}|X(y)|^2|x-y|^{2m-n}\,d\mu_g^{(2m)}(y)\right)\\
			&\quad\times
			\left(\int_{\mathbb R^n}|X(y)|^4|x-y|^{2m-n}\,d\mu_g^{(2m)}(y)\right)\\
			&\quad+2\left(\int_{\mathbb R^n}|x-y|^{2m-n}\,d\mu_g^{(2m)}(y)\right)
			\left(\int_{\mathbb R^n}|X(y)|^2|x-y|^{2m-n}\,d\mu_g^{(2m)}(y)\right)^3.
		\end{aligned}
		\]
		By symmetry in \(y_1\) and \(y_3\),
		\[
		\mathcal A_3=\mathcal A_1.
		\]
		Consequently,
		\begin{equation}\label{eqn:A1+A3-expression}
			\begin{aligned}
				\mathcal A_1+\mathcal A_3
				&=4\left(\int_{\mathbb R^n}|x-y|^{2m-n}\,d\mu_g^{(2m)}(y)\right)^2
				\left(\int_{\mathbb R^n}|X(y)|^2|x-y|^{2m-n}\,d\mu_g^{(2m)}(y)\right)\\
				&\quad\times
				\left(\int_{\mathbb R^n}|X(y)|^4|x-y|^{2m-n}\,d\mu_g^{(2m)}(y)\right)\\
				&\quad+4\left(\int_{\mathbb R^n}|x-y|^{2m-n}\,d\mu_g^{(2m)}(y)\right)
				\left(\int_{\mathbb R^n}|X(y)|^2|x-y|^{2m-n}\,d\mu_g^{(2m)}(y)\right)^3.
			\end{aligned}
		\end{equation}
		We next consider \(\mathcal A_2\). For \(\mathcal A_2\), the expressions involving \(y_1\) and \(y_3\) have the same form. Hence, Fubini's theorem and (\ref{eqn:Xy-int-equ-0}) give
		\begin{equation}\label{eqn:A2-refine-expression}
			\begin{aligned}
				\mathcal A_2
				&=2\Bigg| \left( \int_{\mathbb R^n} |x-y|^{2m-n}\,d\mu_g^{(2m)}(y) \right) \int_{\mathbb R^n} |X(y_1)|^2X(y_1)|x-y_1|^{2m-n}\, d\mu_g^{(2m)}(y_1)\\
				&\quad+ \left( \int_{\mathbb R^n} |X(y)|^2|x-y|^{2m-n}\,d\mu_g^{(2m)}(y) \right) \int_{\mathbb R^n} X(y_1)|x-y_1|^{2m-n}\, d\mu_g^{(2m)}(y_1) \Bigg|^2\\
				&=2\left(\int_{\mathbb R^n}|x-y|^{2m-n}\,d\mu_g^{(2m)}(y)\right)^2
				\times\left|\int_{\mathbb R^n}|X(y)|^2X(y)|x-y|^{2m-n}\,d\mu_g^{(2m)}(y)\right|^2.
			\end{aligned}
		\end{equation}
		Substituting \((\ref{eqn:A1+A3-expression})\) and \((\ref{eqn:A2-refine-expression})\) into (\ref{eqn:I6-A1-A2-A3-expression}), we arrive at
		\begin{equation}\label{eqn:I6-new-expression}
			\begin{aligned}
				I_6&=
				4\left(	\int_{\mathbb R^n}|x-y|^{2m-n}\,d\mu_g^{(2m)}(y)\right)^4\left(\int_{\mathbb R^n}|X(y)|^2|x-y|^{2m-n}\,d\mu_g^{(2m)}(y)\right)\\
				&\quad\times
				\left(\int_{\mathbb R^n}|X(y)|^4|x-y|^{2m-n}\,d\mu_g^{(2m)}(y)\right)\\
				&\quad+4\left(\int_{\mathbb R^n}|x-y|^{2m-n}\,d\mu_g^{(2m)}(y)\right)^3
				\left(\int_{\mathbb R^n}|X(y)|^2|x-y|^{2m-n}\,d\mu_g^{(2m)}(y)\right)^3\\
				&\quad-2\left(\int_{\mathbb R^n}|x-y|^{2m-n}\,d\mu_g^{(2m)}(y)\right)^4\left|\int_{\mathbb R^n}|X(y)|^2X(y)|x-y|^{2m-n}\,d\mu_g^{(2m)}(y)\right|^2.
			\end{aligned}
		\end{equation}
		Subtracting the expression for \(I_7\) gives
		\begin{equation}\label{eqn:I6-I6-minus-I7-expression}
			\begin{aligned}
				I_6-I_7&=
				2\left(\int_{\mathbb R^n}|x-y|^{2m-n}\,d\mu_g^{(2m)}(y)\right)^3\\
				&\quad\times\Bigg[\left(\int_{\mathbb R^n}|x-y|^{2m-n}\,d\mu_g^{(2m)}(y)\right)
				\left(\int_{\mathbb R^n}|X(y)|^2|x-y|^{2m-n}\,d\mu_g^{(2m)}(y)\right)\\
				&\quad\times\left(\int_{\mathbb R^n}|X(y)|^4|x-y|^{2m-n}\,d\mu_g^{(2m)}(y)\right)-
				\left(\int_{\mathbb R^n}|X(y)|^2|x-y|^{2m-n}\,d\mu_g^{(2m)}(y)\right)^3\\
				&\quad-\left(\int_{\mathbb R^n}|x-y|^{2m-n}\,d\mu_g^{(2m)}(y)\right)
				\left|\int_{\mathbb R^n}|X(y)|^2X(y)|x-y|^{2m-n}\,d\mu_g^{(2m)}(y)\right|^2\Bigg].
			\end{aligned}
		\end{equation}
		Using (\ref{eqn:Xy-int-equ-0}) and the Cauchy-Schwarz inequality, we obtain that
		\[
		\begin{aligned}
			&\left(\int_{\mathbb R^n}|x-y|^{2m-n}\,d\mu_g^{(2m)}(y)\right)
			\left|\int_{\mathbb R^n}|X(y)|^2X(y)|x-y|^{2m-n}\,d\mu_g^{(2m)}(y)\right|^2\\
			&=\left(\int_{\mathbb R^n}	|x-y|^{2m-n}\,d\mu_g^{(2m)}(y)\right)
			\Bigg|\int_{\mathbb R^n}|X(y)|^2X(y)|x-y|^{2m-n}\,d\mu_g^{(2m)}(y)\\
			&\quad	-\frac{\displaystyle\int_{\mathbb R^n}|X(z)|^2|x-z|^{2m-n}\,d\mu_g^{(2m)}(z)}
			{\displaystyle\int_{\mathbb R^n}|x-z|^{2m-n}\,d\mu_g^{(2m)}(z)}
			\int_{\mathbb R^n}X(y)|x-y|^{2m-n}\,d\mu_g^{(2m)}(y)\Bigg|^2\\
			&\leq\left(\int_{\mathbb R^n}|x-y|^{2m-n}\,d\mu_g^{(2m)}(y)\right)\left(\int_{\mathbb R^n}|X(y)|^2|x-y|^{2m-n}\,d\mu_g^{(2m)}(y)\right)\\
			&\quad\times\Bigg(\int_{\mathbb R^n}\Bigg(|X(y)|^2-\frac{\displaystyle\int_{\mathbb R^n}|X(z)|^2|x-z|^{2m-n}\,d\mu_g^{(2m)}(z)}{\displaystyle\int_{\mathbb R^n}|x-z|^{2m-n}\,
				d\mu_g^{(2m)}(z)}\Bigg)^2|x-y|^{2m-n}\,d\mu_g^{(2m)}(y)\Bigg)\\
		\end{aligned}
		\]
		Expanding the square in the last integral and relabeling the variables, we further obtain that
		\[
		\begin{aligned}
			&\left(\int_{\mathbb R^n}|x-y|^{2m-n}\,d\mu_g^{(2m)}(y)\right)
			\left|\int_{\mathbb R^n}|X(y)|^2X(y)|x-y|^{2m-n}\,d\mu_g^{(2m)}(y)\right|^2\\
			&\leq\left(\int_{\mathbb R^n}|x-y|^{2m-n}\,d\mu_g^{(2m)}(y)\right)\left(\int_{\mathbb R^n}|X(y)|^2|x-y|^{2m-n}\,d\mu_g^{(2m)}(y)\right)\\
			&\quad\times\Bigg[\int_{\mathbb R^n}|X(y)|^4|x-y|^{2m-n}\,d\mu_g^{(2m)}(y)\\
			&\quad-2\frac{\displaystyle\int_{\mathbb R^n}|X(z)|^2|x-z|^{2m-n}\,d\mu_g^{(2m)}(z)}
			{\displaystyle\int_{\mathbb R^n}|x-z|^{2m-n}\,d\mu_g^{(2m)}(z)}
			\left(\int_{\mathbb R^n}|X(y)|^2|x-y|^{2m-n}\,d\mu_g^{(2m)}(y)\right)\\
			&\quad+\left(\frac{\displaystyle\int_{\mathbb R^n}|X(z)|^2|x-z|^{2m-n}\,d\mu_g^{(2m)}(z)}
			{\displaystyle\int_{\mathbb R^n}|x-z|^{2m-n}\,d\mu_g^{(2m)}(z)}\right)^2
			\left(\int_{\mathbb R^n}|x-y|^{2m-n}\,d\mu_g^{(2m)}(y)\right)\Bigg]\\
			&=\left(\int_{\mathbb R^n}|x-y|^{2m-n}\,d\mu_g^{(2m)}(y)\right)
			\left(\int_{\mathbb R^n}|X(y)|^2|x-y|^{2m-n}\,d\mu_g^{(2m)}(y)\right)\\
			&\quad\times\left(\int_{\mathbb R^n}|X(y)|^4|x-y|^{2m-n}\,d\mu_g^{(2m)}(y)	\right)\\
			&\quad-\left(\int_{\mathbb R^n}|X(y)|^2|x-y|^{2m-n}\,d\mu_g^{(2m)}(y)\right)^3.
		\end{aligned}
		\]
		Combining this with (\ref{eqn:I6-I6-minus-I7-expression}), we have \(I_6-I_7\geq0\),
		that is \(I_7\leq I_6\).
	\end{proof}

	\section{Extended Positivity Range for Sixth-Order Q-Curvature}\label{sec-3}
	
	The range \(2m\le n\le 4m-6\) is covered by Theorem 1.1 of \cite{Vetois-Zeitler-2026}. Therefore, it remains to treat the case \(4m-6<n\le 4m-6+[\widehat h_m]\).
	We temporarily assume \(4m-6<n\le6m-6\) and will verify below that this restriction holds throughout the above range. For convenience, set
	\[
	h:=n-4m+6,\qquad
	A:=n-2m+2,\qquad
	B:=n-2m+4,\qquad
	C:=n-3m+3,
	\]
	and
	\[
	E:=3n-8m+6,\qquad
	F:=2n-5m+3.
	\]
	Using \(n=4m-6+h\), these parameters can be rewritten as
	\[
	A=2m-4+h,\qquad
	B=2m-2+h,\qquad
	C=m-3+h,
	\]
	\[
	E=4m-12+3h,\qquad
	F=3m-9+2h.
	\]
	These expressions also show that \(A,B,C,E\) and \(F\) are positive, because \(m\ge4\) and \(h>0\). Substituting these identities into the curvature formula of Proposition \ref{pro:vetois-Zeitler-result-pro2.1} yields
	\begin{equation}\label{eqn:sumI-Q6g}
		\begin{aligned}
			\frac{\widetilde{Q}_g^{(6)}(x)}{2(m-3)}
			&=(2m-n-2)^2(2m-n-4)^2 I_1 \\
			&\quad-6(4m-n-6)(2m-n-2)^2(2m-n-4) I_2 \\
			&\quad+3(4m-n-6)(3m-n-3)(2m-n-2)^2 I_3 \\
			&\quad+4(4m-n-6)(3m-n-3)(2m-n-2)(2m-n-4) I_4 \\
			&\quad-2(4m-n-6)(2m-n-2)^3 I_5 \\
			&\quad+12(4m-n-6)(3m-n-3)(2m-n-2)^2 I_6 \\
			&\quad-6(4m-n-6)(3m-n-3)(8m-3n-6)(2m-n-2) I_7 \\
			&\quad+(4m-n-6)(3m-n-3)(5m-2n-3)(8m-3n-6) I_8\\
			&=
			a_1I_1+a_2I_2+a_3I_3+a_4I_4+a_5I_5+a_6I_6+a_7I_7+a_8I_8,
		\end{aligned}
	\end{equation}
	where
	\[
	\begin{aligned}
		a_1&:=A^2B^2,
		& a_2&:=-6hA^2B,
		& a_3&:=3hCA^2,
		& a_4&:=4hCAB,\\
		a_5&:=-2hA^3,
		& a_6&:=12hCA^2,
		& a_7&:=-6hCEA,
		& a_8&:=hCFE.
	\end{aligned}
	\]
	Among the coefficients in (\ref{eqn:sumI-Q6g}), only \(a_2\), \(a_5\), and \(a_7\) are negative when \(h>0\). We first handle the contribution involving \(I_7\). Since \(a_3,a_6>0\), Lemma \ref{lem:I7-leq-I3-I6-I1-I4} \((i)\) and \((iii)\) imply that \(I_7\le I_3\) and \(I_7\le I_6\). Hence,
	\[
	\begin{aligned}
		a_3I_3+a_6I_6+a_7I_7&\ge(a_3+a_6+a_7)I_7\\
		&=(3hCA^2+12hCA^2-6hCEA)I_7\\
		&=3hCA(6m-n-2)I_7.
	\end{aligned}
	\]
	Therefore, whenever \(n\le 6m-2\),
	\begin{equation}\label{eqn:a3-a6-a7-ge-0}
		a_3I_3+a_6I_6+a_7I_7\ge0.
	\end{equation}
	We now estimate the remaining terms in \((\ref{eqn:sumI-Q6g})\). Using \((\ref{eqn:I5-leq-I1})\), \((\ref{eqn:a3-a6-a7-ge-0})\), and Lemma \ref{lem:I7-leq-I3-I6-I1-I4} \((ii)\), together with \(a_2<0\) and \(a_5<0\), we obtain that
	\begin{equation}\label{eqn:Q6-want-divide-I1}
		\begin{aligned}
			\frac{\widetilde{Q}_g^{(6)}(x)}{2(m-3)}
			&\ge a_1I_1+a_2I_2+a_4I_4+a_5I_5+a_8I_8\\
			&\ge (a_1+a_5)I_1+a_2I_2+a_4I_4+a_8I_8\\
			&\ge(a_1+a_5)I_1+a_2\sqrt{I_1I_4}+a_4I_4+a_8I_8
		\end{aligned}
	\end{equation}
	for \(n\le6m-2\). Meanwhile, assumption (A) ensures that \(Q_g^{(2m)}>0\) in \(\mathbb R^n\setminus B(0,R)\), see also \cite[Page 15]{Vetois-Zeitler-2026}. Because \(u>0\), the measure \(d\mu_g^{(2m)}\) has a strictly positive density on this region. Moreover, \(V(y)=\frac{x-y}{|x-y|^2}\) is not constant on any open subset of \(\mathbb R^n\setminus\{x\}\). Hence \(Z_{1,2}=|Y_{1,2}|^2=|V(y_1)-V(y_2)|^2\) is positive on a set of positive product measure, and the definition of \(I_1\) gives \(I_1>0\). Define 
	\[
	t:=\sqrt{\frac{I_4}{I_1}}.
	\]
	 By \((\ref{eqn:I4-leq-I1})\), we have \(I_4\le I_1\) and hence \(0\le t\le1\). Dividing (\ref{eqn:Q6-want-divide-I1}) by \(I_1\), we arrive at
	\begin{equation}\label{eqn:Qg6-over-I1-geq-qt}
		\frac{\widetilde{Q}_g^{(6)}(x)}{2(m-3)I_1}\ge
		q(t)+a_8\frac{I_8}{I_1},
	\end{equation}
	where \(q(t):=a_4t^2+a_2t+a_1+a_5\).
	Since \(a_4>0\) and \(t\in[0,1]\), we have
	\[
	\begin{aligned}
		q'(t)&=2a_4t+a_2
		\le 2a_4+a_2
		=8hCAB-6hA^2B
		=2hAB(n-6m+6).
	\end{aligned}
	\]
	If \(n\le6m-6\), the derivative \(q'(t)\) is non-positive in \([0,1]\). Thus \(q\) is decreasing on this interval and \(q(t)\ge q(1)\).
	Combining this with (\ref{eqn:Qg6-over-I1-geq-qt}), we derive that
	\begin{equation}\label{eqn:Qg6-geq-q1I1-I8}
		\begin{aligned}
			\frac{\widetilde{Q}_g^{(6)}(x)}{2(m-3)}&\ge q(t)I_1+a_8I_8\\
			&\ge q(1)I_1+a_8I_8\\
			&=(a_1+a_2+a_4+a_5)I_1+a_8I_8\\
			&=A(AB^2-6hAB+4hCB-2hA^2)I_1+a_8I_8\\
			&=A\Big[ 8(m-2)(m-1)^2 -(12m^2-40m+36)h-(14m-24)h^2 -3h^3 \Big]I_1+a_8I_8\\ &=A\,G_m(h)I_1+a_8I_8,
		\end{aligned}
	\end{equation}
	where \(G_m(h)\) is given in (\ref{eqn:G-m-h-define}). For the quadratic polynomial \(12m^2-40m+36\), the discriminant is \(\Delta=(-40)^2-4\cdot12\cdot36=-128<0\),
	and its leading coefficient is positive. Thus, for \(m\in\mathbb R\), \(12m^2-40m+36>0\).
	Moreover, it follows from the fact \(14m-24>0\) for \(m\ge4\) that
	\[
	\begin{aligned}
		G_m'(h)
		&=-(12m^2-40m+36)-2(14m-24)h-9h^2<0.
	\end{aligned}
	\]
	Hence, \(G_m\) is strictly decreasing on \((0,\infty)\). Combining this with
	\[
	G_m(0)=8(m-2)(m-1)^2>0
	\]
	and
	\[
	G_m(2m)
	=
	-16(m-1)(6m^2-3m-1)<0,
	\]
	we know that \(G_m\) has a unique positive zero \(\widehat h_m\) satisfying \(0<\widehat h_m<2m\).
	This bound also verifies the dimensional restrictions used above. Indeed, if \(0<h\le\widehat h_m\), then \(h<2m\), and hence \(n=4m-6+h<6m-6<6m-2\).
	Therefore, both the condition required in \((\ref{eqn:a3-a6-a7-ge-0})\) and the monotonicity condition for \(q\) hold throughout this range. Moreover, the strict decrease of \(G_m\) gives \(G_m(h)\ge0\) for \(0<h\le\widehat h_m\).
	Suppose that \(G_m(h)>0\). From \((\ref{eqn:Qg6-geq-q1I1-I8})\), 
	\[
	\frac{\widetilde{Q}_g^{(6)}(x)}{2(m-3)}\ge A\,G_m(h)I_1+a_8I_8. 
	\]
	Here \(A>0\), \(I_1>0\), and \(a_8I_8\ge0\). Consequently, \(\widetilde{Q}_g^{(6)}(x)>0\).
	It remains only to consider the possible endpoint at which \(G_m(h)=0\). By the uniqueness of the positive zero, this can occur only when \(h=\widehat h_m\). In this case, \((\ref{eqn:Qg6-geq-q1I1-I8})\) reduces to \(\frac{\widetilde{Q}_g^{(6)}(x)}{2(m-3)} \ge a_8 I_8\).
	By the definitions of \(C\), \(F\), \(E\) and \(h>0\), we get \(a_8=hCFE>0\).
	As in the argument for \(I_1\), we also have \(I_8>0\). Therefore, \(\widetilde{Q}_g^{(6)}(x)>0\)
	at \(h=\widehat h_m\). Recalling that \(h=n-4m+6\), the condition \(0<h\le \widehat h_m\) is equivalent to \(4m-6<n\le 4m-6+\widehat h_m\).
	Because \(n\) is an integer, this becomes \(4m-6<n\le 4m-6+[\widehat h_m]\).
	Since \(Q_g^{(6)}\) and \(\widetilde Q_g^{(6)}\) have the same sign in the non-critical case, the preceding argument establishes the desired positivity in the range \(4m-6<n\le 4m-6+[\widehat h_m]\). Combining this with the previously known positivity range \(2m\le n\le 4m-6\) in \cite[Theorem 1.1]{Vetois-Zeitler-2026}, we conclude that \({Q}_g^{(6)}(x)>0\) whenever \(2m\le n\le 4m-6+[\widehat h_m]\).
	
	Next, we derive a numerical estimate for \(\widehat h_m\).
	Since \(m\ge4\), a direct computation gives
	\[
	\begin{aligned}
		G_m\left(\frac{3m}{10}\right)
		&=\frac{3059}{1000}(m-4)^3+\frac{4717}{250}(m-4)^2+\frac{4164}{125}(m-4)+\frac{1392}{125}
		>0.
	\end{aligned}
	\]Since \(G_m\) is strictly decreasing and \(G_m(\widehat h_m)=0\), we obtain \(\widehat h_m>\frac{3m}{10}\).
	To determine the asymptotic behavior of \(\widehat h_m\), set 
	\[
	\rho_m:=\frac{\widehat h_m}{m}.
	\]
	We first show that \(\{\rho_m\}_{m\ge4}\) is strictly increasing. Define 
	\[
	F_m(\rho):=\frac{G_m(m\rho)}{m^3}.
	\]
	A direct computation gives
	\[
	\begin{aligned}
		F_m(\rho)&=\frac{8(m-2)(m-1)^2}{m^3}-\left(12-\frac{40}{m}+\frac{36}{m^2}\right)\rho-\left(14-\frac{24}{m}\right)\rho^2-3\rho^3\\
		&=8-12\rho-14\rho^2-3\rho^3+\frac{8(3\rho^2+5\rho-4)}{m}+\frac{4(10-9\rho)}{m^2}-\frac{16}{m^3}.
	\end{aligned}
	\]
	Differentiating with respect to \(\rho\), we obtain that
	\[
	\begin{aligned}
		F'_m(\rho)&=-\left(12-\frac{40}{m}+\frac{36}{m^2}\right)-2\left(14-\frac{24}{m}\right)\rho-9\rho^2\\
		&=-\frac{12m^2-40m+36}{m^2}-2\left(14-\frac{24}{m}\right)\rho-9\rho^2<0
	\end{aligned}
	\]
	for \(m\ge4\) and \(\rho\ge0\). By definition, \(F_m(\rho_m)=\frac{G_m(m\rho_m)}{m^3}=\frac{G_m(\widehat h_m)}{m^3}=0\).
	To identify the limiting behavior of \(\rho_m\), we first consider \(F_m(\rho)\) for fixed \(\rho\). As \(m\to\infty\),
	\[
	F_m(\rho)\rightarrow8-12\rho-14\rho^2-3\rho^3=-(\rho+2)(3\rho^2+8\rho-4).
	\]
	Hence, the above limiting polynomial has a unique positive root
	\(\rho_0:=\frac{2\sqrt7-4}{3}\).
	Substituting \(\rho=\rho_0\) into \(F_m\) gives
	\[
	\begin{aligned}
		F_m(\rho_0)
		&=\frac{8}{m^3}\left[(4-2\sqrt7)m^2+(11-3\sqrt7)m-2\right]
		=\frac{8}{m^3}\Theta(m),
	\end{aligned}
	\]
	where \(\Theta(m):=(4-2\sqrt7)m^2+(11-3\sqrt7)m-2\).
	For \(m\ge4\), we have
	\[
	\begin{aligned}
		\Theta'(m)&=2(4-2\sqrt7)m+11-3\sqrt7
		\leq 8(4-2\sqrt7)+11-3\sqrt7
		=43-19\sqrt7<0,
	\end{aligned}
	\]
	which implies that \(\Theta(m)\) is strictly decreasing on \([4,\infty)\). Together with \(\Theta(4)=106-44\sqrt7<0\), this yields \(\Theta(m)<0\) for \(m\ge4\).
	Consequently, \(F_m(\rho_0)=\frac{8}{m^3}\Theta(m)<0\).
	Recalling that \(F_m\) is strictly decreasing with respect to \(\rho\), while \(F_m(\rho_m)=0\) and \(F_m(\rho_0)<0\), the zero \(\rho_m\) must lie to the left of \(\rho_0\). Hence, \(0<\rho_m<\rho_0\).
	
	We next examine how the zero \(\rho_m\) varies with \(m\). Define
	\[
	F(m,\rho):=F_m(\rho)=8-12\rho-14\rho^2-3\rho^3+\frac{8(3\rho^2+5\rho-4)}{m}+\frac{4(10-9\rho)}{m^2}-\frac{16}{m^3}.
	\]
	Differentiation with respect to \(m\) gives
	\begin{equation}\label{eqn:F-partial-m}
		\begin{aligned}
			\frac{\partial F (m,\rho)}{\partial m}&=-\frac{8(3\rho^2+5\rho-4)}{m^2}-\frac{8(10-9\rho)}{m^3}+\frac{48}{m^4}\\
			&=\frac{8}{m^4}\left[(4-5\rho-3\rho^2)m^2-(10-9\rho)m+6\right].
		\end{aligned}
	\end{equation}
	Thus, the sign of \(\frac{\partial F}{\partial m}\) is determined by \(H(m,\rho):=(4-5\rho-3\rho^2)m^2-(10-9\rho)m+6\).
	For \(m\ge4\) and \(0\le\rho\le\rho_0\),
	\[
	\begin{aligned}
		\frac{\partial H(m,\rho)}{\partial\rho}&=-(5+6\rho)m^2+9m
		=m\bigl[9-(5+6\rho)m\bigr]\leq m(9-5m)<0,
	\end{aligned}
	\]
	which means that \(H(m,\rho)\) is strictly decreasing in \(\rho\) on \([0,\rho_0]\). Therefore,
	\[
	H(m,\rho)\ge H(m,\rho_0)=(2\sqrt7-4)m^2+(6\sqrt7-22)m+6.
	\]
	The positivity of \(2\sqrt7-4\) gives, for \(m\ge4\),
	\[
	\begin{aligned}
		\frac{\partial H(m,\rho_0)}{\partial m}&=2(2\sqrt7-4)m+6\sqrt7-22\\
		&\ge 8(2\sqrt7-4)+6\sqrt7-22
		=22\sqrt7-54>0,
	\end{aligned}
	\]
	which means that \(H(m,\rho_0)\) is strictly increasing with respect to \(m\) on \([4,\infty)\). Together with the monotonicity of \(H(m,\rho)\) in \(\rho\), this yields
	\[
	\begin{aligned}
		H(m,\rho)
		&\ge H(m,\rho_0)
		\ge H(4,\rho_0)
		=16(2\sqrt7-4)+4(6\sqrt7-22)+6
		=56\sqrt7-146>0
	\end{aligned}
	\]
	for \(m\ge4\) and \(0\le\rho\le\rho_0\). In view of \((\ref{eqn:F-partial-m})\), it follows that \(\frac{\partial F(m,\rho)}{\partial m}>0\) throughout the same range. Since \(0<\rho_m<\rho_0\), the preceding monotonicity in \(m\) applies at \(\rho=\rho_m\). Hence,
	\[
	F_{m+1}(\rho_m)=F(m+1,\rho_m)>F(m,\rho_m)=F_m(\rho_m)=0.
	\]
	On the other hand, \(F_{m+1}(\rho_{m+1})=0\), while \(F_{m+1}\) is strictly decreasing with respect to \(\rho\). Hence, \(\rho_{m+1}>\rho_m\) and the sequence \(\{\rho_m\}_{m\ge4}\) is strictly increasing. Moreover, it follows from 
	\(0<\rho_m<\rho_0\) that there exists \(\rho_*\in(0,\rho_0]\) such that \(\rho_m\rightarrow\rho_*\).
	Recalling that \(F_m(\rho_m)=0\), we pass to the limit as \(m\to\infty\) in \(F_m(\rho_m)\) and this gives
	\[
	\begin{aligned}
		0&=\lim_{m\to\infty}F_m(\rho_m)\\
		&=\lim_{m\to\infty}\Big[8-12\rho_m-14\rho^2_m-3\rho^3_m+\frac{8(3\rho^2_m+5\rho_m-4)}{m}+\frac{4(10-9\rho_m)}{m^2}-\frac{16}{m^3}\Big]\\
		&=8-12\rho_*-14\rho_*^2-3\rho_*^3\\
		&=-(\rho_*+2)(3\rho_*^2+8\rho_*-4).
	\end{aligned}
	\]
	The positivity of \(\rho_*\) rules out the roots \(-2\) and \(\frac{-4-2\sqrt7}{3}\). Therefore, \(\rho_*\) must be the unique positive root, that is \(\rho_*=\rho_0=\frac{2\sqrt7-4}{3}\).
	Consequently, 
	\[
	\frac{\widehat h_m}{m}=\rho_m\nearrow\frac{2\sqrt7-4}{3}
	\]
	as \(m\to\infty\). We finally return to the integer part \([\widehat h_m]\). By definition, \(0\le \widehat h_m-[\widehat h_m]<1\).
	After dividing by \(m\), we obtain that 
	\[
	0\le\frac{\widehat h_m}{m}-\frac{[\widehat h_m]}{m}<\frac{1}{m}.
	\]
	Together with the limit above, this yields \({[\widehat h_m]}\sim\frac{2\sqrt7-4}{3}m\)
	as \(m\to\infty\).
	It follows that
	\[
	4m-6+[\widehat h_m]\sim\Big(4+\frac{2\sqrt7-4}{3}\Big)m
	=\frac{8+2\sqrt7}{3} m\quad\mbox{as}\quad m\to\infty.
	\]
	
\begin{figure}[H]
		\centering
		\begin{minipage}{0.48\textwidth}
			\centering
			\includegraphics[width=\textwidth]{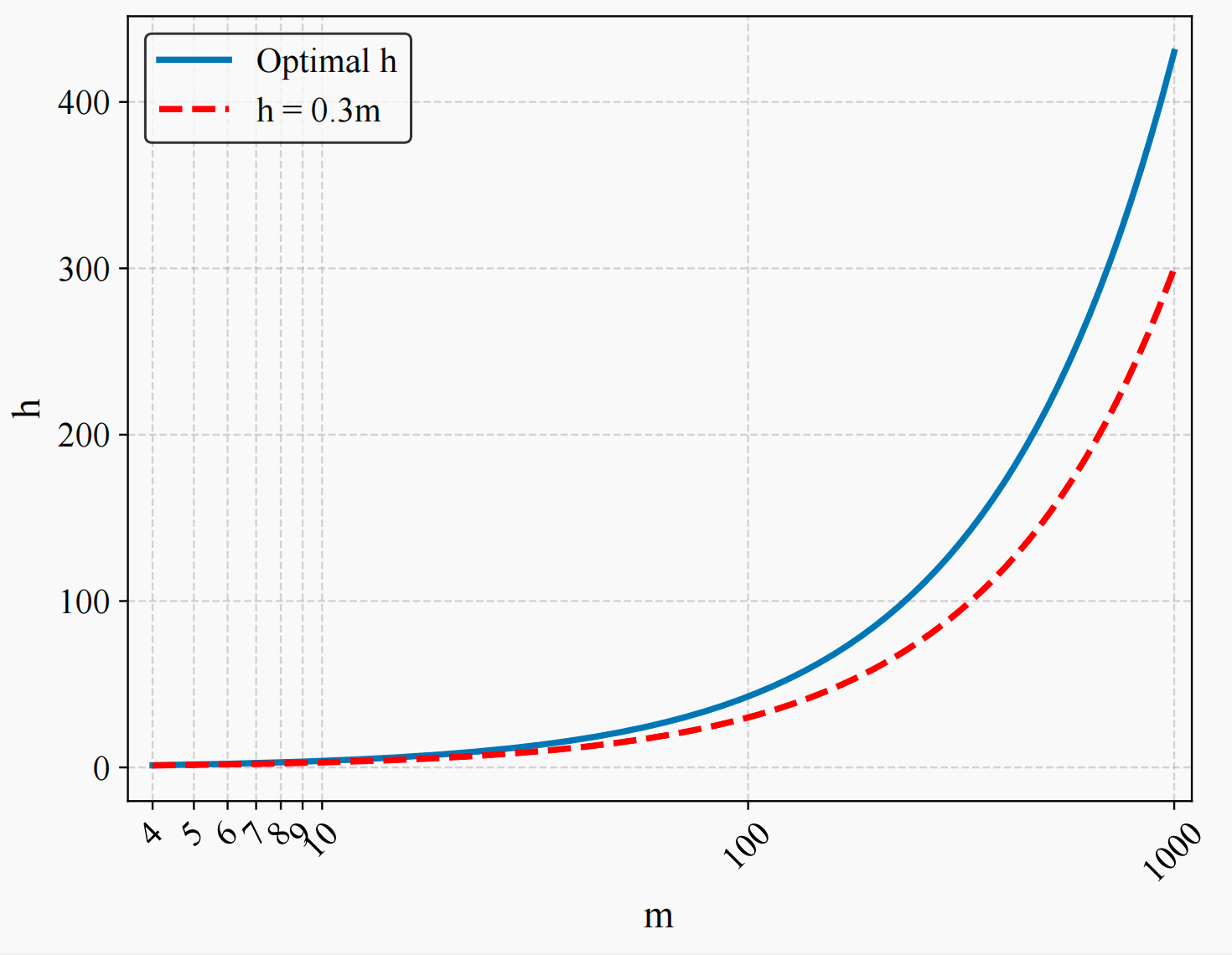}
			\caption{Comparison between the threshold parameter $\widehat h_m$ and the choice $\frac{3m}{10}$ on a logarithmic horizontal scale.}
			\label{fig:h_comparison-1}
		\end{minipage}
		\hfill
		\begin{minipage}{0.48\textwidth}
			\centering
			\includegraphics[width=\textwidth]{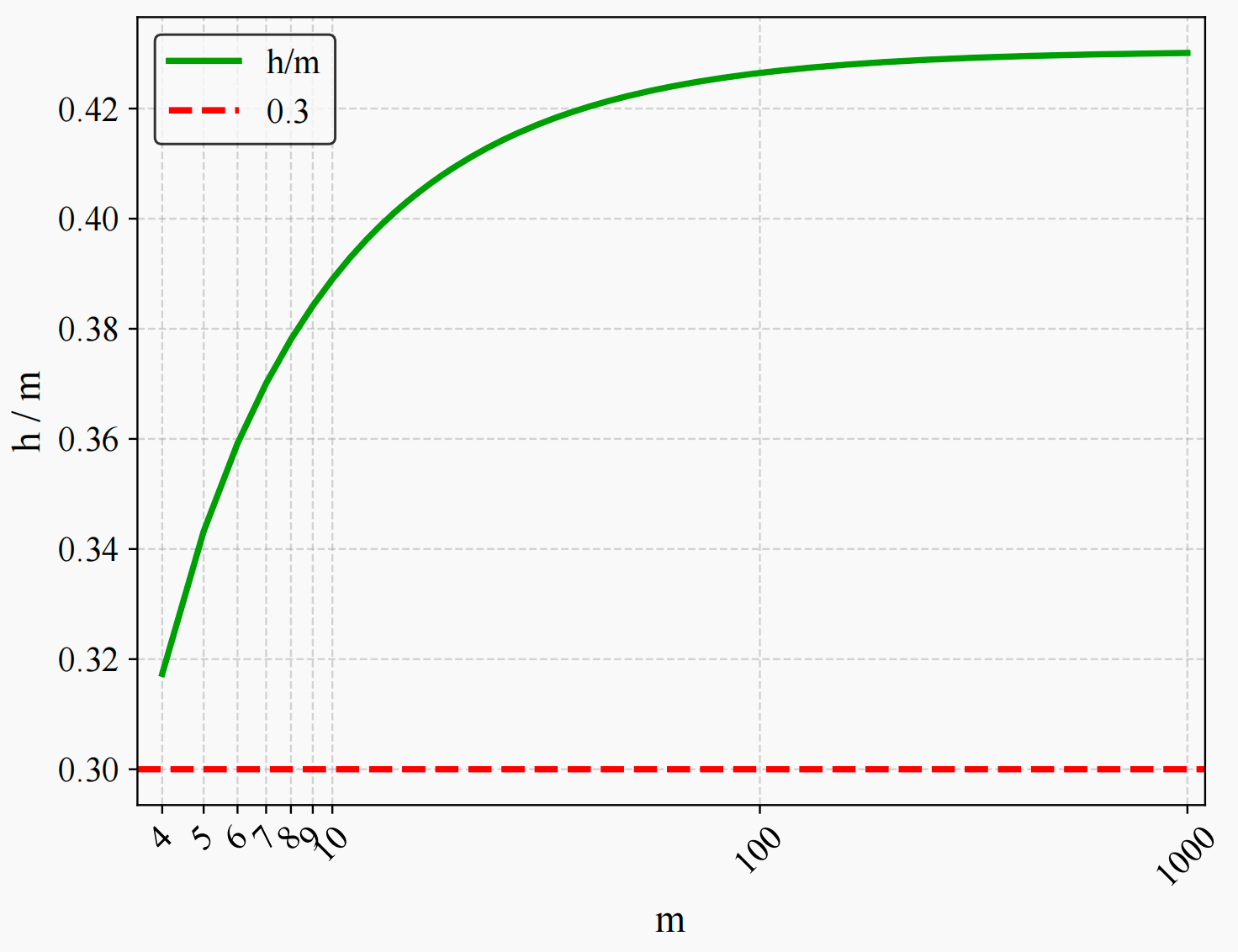}
			\caption{Comparison between the normalized threshold parameter $\frac{\widehat h_m}{m}$ and the choice $\frac{3}{10}$.}
			\label{fig:h_comparison-2}
		\end{minipage}
	\end{figure}

	To further illustrate the behavior of the threshold parameter, we provide two numerical plots generated by mathematical software. Figure \ref{fig:h_comparison-1} presents the relationship between \(m\), the selected value \(h=\frac{3m}{10}\), and the threshold value \(\widehat h_m\), where \(\widehat h_m\) denotes the unique positive root of \(G_m(h)=0\). For visualization over a wide range of \(m\), the horizontal axis is displayed on a logarithmic scale. The value \(\frac{3m}{10}\) provides a real-valued lower comparison for \(\widehat h_m\). In integer dimensions, it
	yields the guaranteed endpoint \(4m-6+[\frac{3m}{10}]\), while the theorem uses the sharper endpoint \(4m-6+[\widehat h_m]\). Figure \ref{fig:h_comparison-2} further illustrates the behavior of the normalized quantity \(\frac{\widehat h_m}{m}\). The numerical results show that \(\frac{\widehat h_m}{m}\) increases rapidly for small values of \(m\) and gradually stabilizes as \(m\) becomes large. This behavior agrees with the theoretical limiting value 
	\[
	\frac{\widehat h_m}{m}\nearrow\frac{2\sqrt7-4}{3}\simeq 0.4305
	\]
	as \(m\to\infty\). Meanwhile, we use mathematical software to perform numerical computations and further illustrate the relationships among \(m\), the threshold value \(\widehat h_m\), and the normalized quantity \(\frac{\widehat h_m}{m}\), see Table \ref{fig:h_comparison-3}.
	\\
	\qed
	
	\begin{table}[H]
		\centering
		\includegraphics[width=0.5\textwidth]{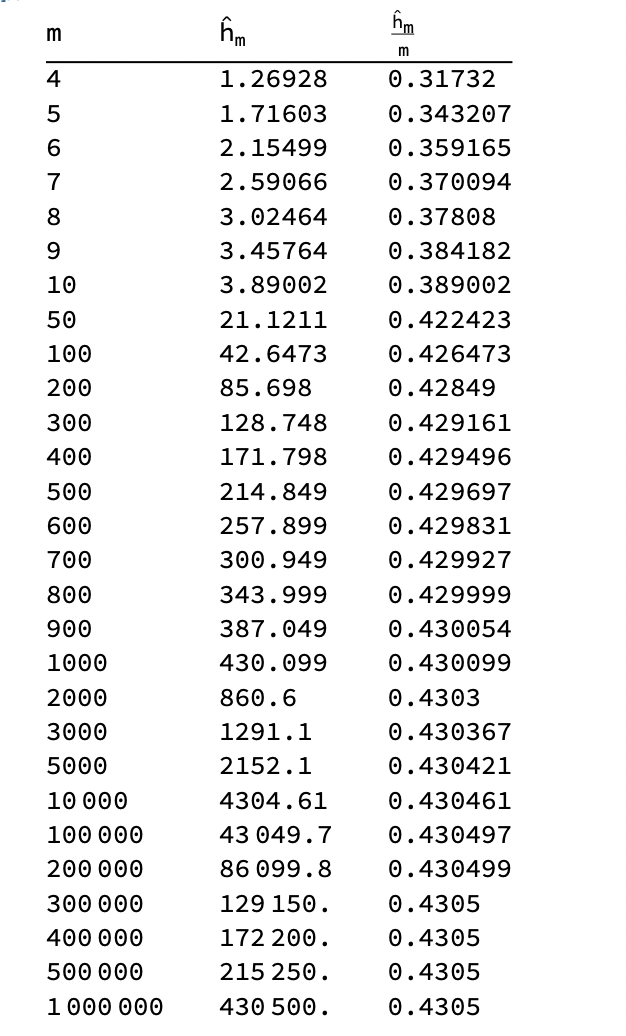}
		\caption{Numerical verification of the threshold parameter $\widehat h_m$ and the normalized quantity $\frac{\widehat h_m}{m}$.}
		\label{fig:h_comparison-3}
	\end{table}

	\section{Quantitative Lower Bounds for Sixth-Order Q-Curvatures}
	Let
	\begin{equation}\label{eqn:section-4-u-express}
		u(x)=
		\begin{cases}
			\displaystyle\int_{\mathbb R^n}|x-y|^{2m-n}\,d\mu_g^{(2m)}(y),&\quad \mbox{if}\quad n>2m;\\[2mm]
			\displaystyle\int_{\mathbb R^n}\ln \Big(\frac{|y|}{|x-y|}\Big)d\mu_g^{(n)}(y)+C,&\quad \mbox{if}\quad n=2m.
		\end{cases}
	\end{equation}
	The assumption (A), together with the total curvature estimate established by Li et al. \cite[Theorem 3.8]{Li-Wei-Xu2025}, ensures that
	\[
	0<\lambda:=\mu_g^{(n)}(\mathbb R^n)\le2
	\]
	in the case \(n=2m\). For fixed \(x\in\mathbb R^n\), we recall the following notation:
	\[
	V(y)=\frac{x-y}{|x-y|^2},\quad Z_i=\frac{1}{|x-y_i|^2},\quad X(y)=V(y)-\overline V(x),
	\]
	\[
	Z_{i,j}=\left|\frac{x-y_i}{|x-y_i|^2}-\frac{x-y_j}{|x-y_j|^2}\right|^2,\quad\overline V(x)=\frac{\displaystyle\int_{\mathbb R^n}V(y)|x-y|^{2m-n}\,d\mu_g^{(2m)}(y)}{\displaystyle\int_{\mathbb R^n}|x-y|^{2m-n}\,d\mu_g^{(2m)}(y)}.
	\]
	Direct differentiation gives
	\[
	\begin{aligned}
		\nabla\ln|x-y|=\frac{x-y}{|x-y|^2},\ \ \Delta\ln|x-y|=\frac{n-2}{|x-y|^2},
	\end{aligned}
	\]
	\[
	\begin{aligned}
		\nabla|x-y|^{2m-n}
		=-(n-2m)|x-y|^{2m-n}V(y)
	\end{aligned}
	\]
	and
	\[
	\begin{aligned}
		\Delta|x-y|^{2m-n}
		=-2(n-2m)(m-1)|x-y|^{2m-n}|V(y)|^2.
	\end{aligned}
	\]
	As stated in \cite[Page 3]{Vetois-Zeitler-2026}, the integral representations in (\ref{eqn:section-4-u-express}) can be differentiated under the integral sign up to order \(2m-1\). Since \(2m\ge8\), this justifies the following calculations of \(\nabla u(x)\) and \(\Delta u(x)\) in both the critical and non-critical cases. For the non-critical case \(n>2m\), the integral representation of \(u\) gives
	\begin{equation}\label{eqn:nabla-u-x}
		\begin{aligned}
			\nabla u(x)
			&=-(n-2m)\int_{\mathbb R^n}V(y)|x-y|^{2m-n}\,d\mu_g^{(2m)}(y)
			=-(n-2m)u(x)\overline V(x),
		\end{aligned}
	\end{equation}
	and
	\begin{equation}\label{eqn:Delta-u-x}
		\Delta u(x)=-2(n-2m)(m-1)\int_{\mathbb R^n}|V(y)|^2|x-y|^{2m-n}\,d\mu_g^{(2m)}(y).
	\end{equation}
	In the critical case \(n=2m\), the logarithmic integral representation yields
	\[
	\begin{aligned}
		\nabla u(x)
		=-\int_{\mathbb R^n}\frac{x-y}{|x-y|^2}\,d\mu_g^{(n)}(y)
	\end{aligned}
	\]
	and
	\[
	\begin{aligned}
		\Delta u(x)
		=-(n-2)\int_{\mathbb R^n}\frac{1}{|x-y|^2}\,d\mu_g^{(n)}(y).
	\end{aligned}
	\]
	
	\textbf{Proof of Theorem \ref{them:main-3}}. We first establish the relationship between \(Q_g^{(6)}\) and \(Q_g^{(2)}\) in the critical case \(n=2m\). Expanding \(Z_{1,2}\) and applying Fubini's theorem, we find that
	\[
	\begin{aligned}
		&\int_{\mathbb R^{2n}}Z_{1,2}\,d\mu_g^{(n)}(y_1)d\mu_g^{(n)}(y_2)\\
		&=\int_{\mathbb R^{2n}}\left|\frac{x-y_1}{|x-y_1|^2}-\frac{x-y_2}{|x-y_2|^2}\right|^2
		d\mu_g^{(n)}(y_1)d\mu_g^{(n)}(y_2)\\
		&=2\lambda\int_{\mathbb R^n}\frac{1}{|x-y|^2}\,d\mu_g^{(n)}(y)
		-2\left|\int_{\mathbb R^n}\frac{x-y}{|x-y|^2}\,d\mu_g^{(n)}(y)\right|^2\\
		&=2\lambda	\int_{\mathbb R^n}\frac1{|x-y|^2}\,d\mu_g^{(n)}(y)-2|\nabla u(x)|^2.
	\end{aligned}
	\]
	Rearranging this identity gives
	\[
	|\nabla u(x)|^2=\lambda\int_{\mathbb R^n}\frac{1}{|x-y|^2}\,d\mu_g^{(n)}(y)\\
	-\frac12\int_{\mathbb R^{2n}}Z_{1,2}\,d\mu_g^{(n)}(y_1)d\mu_g^{(n)}(y_2).
	\]
	We now substitute these identities into the conformal transformation formula
	\[
	Q_g^{(2)}(x)=-e^{-2u(x)}\left(\Delta u(x)+\frac{n-2}{2}|\nabla u(x)|^2\right),
	\]
	and this gives
	\begin{equation}\label{eqn:n-2n-critical-Q2-expression}
		\begin{aligned}
			Q_g^{(2)}(x)&=-e^{-2u(x)}
			\Bigg[-(n-2)\int_{\mathbb R^n}\frac{1}{|x-y|^2}\,d\mu_g^{(n)}(y)\\
			&\quad +\frac{n-2}{2}\Big(\lambda\int_{\mathbb R^n}\frac{1}{|x-y|^2}\,d\mu_g^{(n)}(y)-\frac12\int_{\mathbb R^{2n}}Z_{1,2}\,d\mu_g^{(n)}(y_1)d\mu_g^{(n)}(y_2)\Big)\Bigg]\\
			&=(n-2)e^{-2u(x)}\Bigg[\left(1-\frac{\lambda}{2}\right)\int_{\mathbb R^n}
			\frac{1}{|x-y|^2}\,d\mu_g^{(n)}(y)+\frac{1}{4}\int_{\mathbb R^{2n}}Z_{1,2}\,
			d\mu_g^{(n)}(y_1)d\mu_g^{(n)}(y_2)\Bigg],
		\end{aligned}
	\end{equation}
	which implies that
	\[
	\begin{aligned}
		\left(1-\frac{\lambda}{2}\right)\int_{\mathbb R^n}\frac{1}{|x-y|^2}\,d\mu_g^{(n)}(y)
		+\frac{1}{4}\int_{\mathbb R^{2n}}Z_{1,2}\,d\mu_g^{(n)}(y_1)d\mu_g^{(n)}(y_2)
		=\frac{e^{2u(x)}}{n-2}Q_g^{(2)}(x).
	\end{aligned}
	\]
	On the other hand, Proposition 2.1 of \cite{Vetois-Zeitler-2026} gives
	\[
	Q_g^{(6)}(x)
	=
	e^{-6u(x)}
	\int_{\mathbb R^{6n}}
	P_{n,m}(Z(x,y))\,d\mu_g^{(n)}(y).
	\]
	As shown in the proof of Theorem 1.1 in \cite{Vetois-Zeitler-2026}, all coefficients of \(P_{n,m}\) are nonnegative when \(n=2m\ge8\) and \(0<\lambda\le2\). Therefore, we have
	\begin{equation}\label{eqn:n-2n-critical-Q6-expression}
		\begin{aligned}
			e^{6u(x)}Q_g^{(6)}(x)&\ge\frac{(n-6)^3(2-\lambda)}{8\lambda^3}
			\left(n-2-\frac{n-6}{2}\lambda\right)\left(n-4-\frac{n-6}{2}\lambda\right)
			\int_{\mathbb R^{6n}}Z_1Z_2Z_3\,d\mu_g^{(n)}(y)\\
			&\quad +\frac{(n-6)^5}{256}\int_{\mathbb R^{6n}}Z_{1,2}Z_{3,4}Z_{5,6}\,d\mu_g^{(n)}(y).
		\end{aligned}
	\end{equation}
	The two integrals on the right-hand side can be evaluated separately. Since \(Z_1Z_2Z_3\) depends only on \(y_1,y_2,y_3\), Fubini's theorem gives
	\[
	\begin{aligned}
		\int_{\mathbb R^{6n}}Z_1Z_2Z_3\,d\mu_g^{(n)}(y)&=
		\left(\int_{\mathbb R^n}\frac{1}{|x-y|^2}\,d\mu_g^{(n)}(y)\right)^3
		\left(\int_{\mathbb R^n}d\mu_g^{(n)}(y)\right)^3
		=\lambda^3\left(\int_{\mathbb R^n}\frac{1}{|x-y|^2}\,d\mu_g^{(n)}(y)\right)^3.
	\end{aligned}
	\]
	Similarly, the three factors \(Z_{1,2}\), \(Z_{3,4}\), and \(Z_{5,6}\) involve disjoint pairs of integration variables. Thus,
	\[
	\begin{aligned}
		\int_{\mathbb R^{6n}}Z_{1,2}Z_{3,4}Z_{5,6}\,d\mu_g^{(n)}(y)
		=\left(\int_{\mathbb R^{2n}}Z_{1,2}\,d\mu_g^{(n)}(y_1)d\mu_g^{(n)}(y_2)\right)^3.
	\end{aligned}
	\]
	Substituting these identities into (\ref{eqn:n-2n-critical-Q6-expression}) leads to
	\begin{equation}\label{eqn:n-2n-critical-Q6-expression-2}
		\begin{aligned}
			e^{6u(x)}Q_g^{(6)}(x)&\ge\frac{(n-6)^3(2-\lambda)}8\left(n-2-\frac{n-6}{2}\lambda\right)\left(n-4-\frac{n-6}{2}\lambda\right)\left(\int_{\mathbb R^n}\frac{1}{|x-y|^2}\,d\mu_g^{(n)}(y)\right)^3\\
			&\quad +\frac{(n-6)^5}{256}\left(\int_{\mathbb R^{2n}}Z_{1,2}\,d\mu_g^{(n)}(y_1)d\mu_g^{(n)}(y_2)\right)^3.
		\end{aligned}
	\end{equation}
	To estimate the coefficients in this expression, we use \(0<\lambda\le2\) and \(n\ge8\).
	First,
	\[
	n-2-\frac{n-6}{2}\lambda\ge n-2-(n-6)=4,
	\]
	while
	\[
	n-4-\frac{n-6}{2}\lambda\ge n-4-(n-6)=2.
	\]
	Consequently,
	\[
	\frac{(n-6)^3(2-\lambda)}8\left(n-2-\frac{n-6}{2}\lambda\right)\left(n-4-\frac{n-6}{2}\lambda\right)\ge(n-6)^3(2-\lambda).
	\]
	Furthermore, it follows from \(0\le1-\frac{\lambda}{2}\le1\) that \(2-\lambda=2\left(1-\frac{\lambda}{2}\right)\ge\left(1-\frac{\lambda}{2}\right)^3\).
	Combining these inequalities gives
	\begin{equation}\label{eqn:term-1}
		\begin{aligned}
			\frac{(n-6)^3(2-\lambda)}8\left(n-2-\frac{n-6}{2}\lambda\right)\left(n-4-\frac{n-6}{2}\lambda\right)\ge(n-6)^3\left(1-\frac{\lambda}{2}\right)^3.
		\end{aligned}
	\end{equation}
	The coefficient of the second term in (\ref{eqn:n-2n-critical-Q6-expression-2}) admits a similar estimate. Since \(n\ge8\) and \((n-6)^2\ge4\), we have
	\begin{equation}\label{eqn:term-2}
		\begin{aligned}
			\frac{(n-6)^5}{256}=\frac{(n-6)^3(n-6)^2}{256}\ge\frac{(n-6)^3}{64}=(n-6)^3\left(\frac14\right)^3.
		\end{aligned}
	\end{equation}
	Applying (\ref{eqn:term-1}) and (\ref{eqn:term-2}) to (\ref{eqn:n-2n-critical-Q6-expression-2}), we arrive at
	\begin{equation}\label{eqn:n-2n-critical-Q6-expression-3}
		\begin{aligned}
			e^{6u(x)}Q_g^{(6)}(x)\ge(n-6)^3\left[\left(\left(1-\frac{\lambda}{2}\right)\!\!\!
			\int_{\mathbb R^n}\!\!\frac{1}{|x-y|^2}\,d\mu_g^{(n)}(y)\right)^3
			+\left(\frac14\int_{\mathbb R^{2n}}\!\!\!\!\!Z_{1,2}\,d\mu_g^{(n)}(y_1)d\mu_g^{(n)}(y_2)\right)^3\right].
		\end{aligned}
	\end{equation}
	To combine the two terms on the right-hand side of (\ref{eqn:n-2n-critical-Q6-expression-3}), we use the elementary inequality
	\[
	a^3+b^3\ge\frac14(a+b)^3,\qquad a,b\ge0,
	\]
	which follows from \(4(a^3+b^3)-(a+b)^3=3(a+b)(a-b)^2\ge0\).
	Applying this inequality to (\ref{eqn:n-2n-critical-Q6-expression-3}) yields
	\[
	\begin{aligned}
		e^{6u(x)}Q_g^{(6)}(x)\ge\frac{(n-6)^3}{4}\Bigg[\left(1-\frac{\lambda}{2}\right)
		\int_{\mathbb R^n}\frac{1}{|x-y|^2}\,d\mu_g^{(n)}(y)+\frac{1}{4}
		\int_{\mathbb R^{2n}}Z_{1,2}\,d\mu_g^{(n)}(y_1)d\mu_g^{(n)}(y_2)\Bigg]^3.
	\end{aligned}
	\]
	Combining this with (\ref{eqn:n-2n-critical-Q2-expression}), we deduce that
	\[
	\begin{aligned}
		e^{6u(x)}Q_g^{(6)}(x)\ge
		\frac{(n-6)^3}{4}\left(\frac{e^{2u(x)}}{n-2}Q_g^{(2)}(x)\right)^3
		=\frac{(m-3)^3}{4(m-1)^3}e^{6u(x)}\bigl(Q_g^{(2)}(x)\bigr)^3.
	\end{aligned}
	\]
	Since \(e^{6u(x)}>0\), division by this factor gives
	\[
	Q_g^{(6)}(x)\ge\frac{(m-3)^3}{4(m-1)^3}\bigl(Q_g^{(2)}(x)\bigr)^3.
	\]
	It remains to verify that the right-hand side is strictly positive.
	By assumption (A), the measure \(d\mu_g^{(n)}\) is non-negative and non-zero. It also has a smooth density. Since \(Z_{1,2}=\left|\frac{x-y_1}{|x-y_1|^2}-\frac{x-y_2}{|x-y_2|^2}\right|^2\)
	is strictly positive whenever \(y_1\ne y_2\), the non-triviality of the measure implies
	\[
	\int_{\mathbb R^{2n}}Z_{1,2}\,d\mu_g^{(n)}(y_1)d\mu_g^{(n)}(y_2)>0.
	\]
	Returning to the integral representation of \(Q_g^{(2)}\) in (\ref{eqn:n-2n-critical-Q2-expression}), we find that \(Q_g^{(2)}(x)>0\) for every \(x\in\mathbb R^n\). Thus,
	\[
	Q_g^{(6)}(x)\ge\frac{(m-3)^3}{4(m-1)^3}\bigl(Q_g^{(2)}(x)\bigr)^3>0
	\quad\text{in }\mathbb R^{2m}.
	\]
	This completes the proof in the critical case.
	
	We now consider the non-critical case \(n>2m\) and begin by expressing the second-order \(Q\)-curvature in terms of \(X(y)\). By using (\ref{eqn:nabla-u-x}) and (\ref{eqn:Delta-u-x}), together with the definitions of \(\overline V(x)\) and \(X(y)\), we deduce that
	\[
	\begin{aligned}
		Q_g^{(2)}(x)&=\frac{2}{n-2}u(x)^{-\frac{n+2}{n-2m}}(-\Delta)\left(u^{\frac{n-2}{n-2m}}\right)(x)\\
		&=-\frac{2}{n-2m}u(x)^{-1-\frac4{n-2m}}\Delta u(x)
		-\frac{4(m-1)}{(n-2m)^2}u(x)^{-2-\frac4{n-2m}}|\nabla u(x)|^2\\
		&=4(m-1)u(x)^{-1-\frac4{n-2m}}\int_{\mathbb R^n}|V(y)-\overline V(x)|^2|x-y|^{2m-n}\,d\mu_g^{(2m)}(y)\\
		&=4(m-1)u(x)^{-1-\frac4{n-2m}}\int_{\mathbb R^n}|X(y)|^2|x-y|^{2m-n}\,d\mu_g^{(2m)}(y).
	\end{aligned}
	\]
	Hence,
	\begin{equation}\label{eqn:Xy-Qg2}
		\int_{\mathbb R^n}|X(y)|^2|x-y|^{2m-n}\,d\mu_g^{(2m)}(y)=
		\frac{u(x)^{1+\frac{4}{n-2m}}}{4(m-1)}Q_g^{(2)}(x).
	\end{equation}
	
	We now estimate \(I_1\). Since \(Z_{1,2}\) depends only on \(y_1\) and \(y_2\), Fubini's theorem gives
	\begin{equation}\label{eqn:I1-new-expression}
		\begin{aligned}
			I_1&=\int_{\mathbb R^{6n}}Z_{1,2}^3\prod_{i=1}^{6}|x-y_i|^{2m-n}\,d\mu_g^{(2m)}(y)\\
			&=\left(\int_{\mathbb R^n}|x-y|^{2m-n}\,d\mu_g^{(2m)}(y)\right)^4\times
			\int_{\mathbb R^{2n}}Z_{1,2}^3\prod_{i=1}^{2}|x-y_i|^{2m-n}
			\,d\mu_g^{(2m)}(y_1)d\mu_g^{(2m)}(y_2)\\
			&=u(x)^4\int_{\mathbb R^{2n}}Z_{1,2}^3	\prod_{i=1}^{2}|x-y_i|^{2m-n}
			\,d\mu_g^{(2m)}(y_1)d\mu_g^{(2m)}(y_2).
		\end{aligned}
	\end{equation}
	Applying H\"older's inequality with exponents \(3\) and \(\frac{3}{2}\), we have
	\[
	\begin{aligned}
		&\int_{\mathbb R^{2n}}Z_{1,2}\prod_{i=1}^{2}|x-y_i|^{2m-n}\,d\mu_g^{(2m)}(y_1)d\mu_g^{(2m)}(y_2)\\
		&=\int_{\mathbb R^{2n}}	\left(Z_{1,2}^3\prod_{i=1}^{2}|x-y_i|^{2m-n}\right)^\frac{1}{3}
		\left(\prod_{i=1}^{2}|x-y_i|^{2m-n}\right)^\frac{2}{3}\,d\mu_g^{(2m)}(y_1)d\mu_g^{(2m)}(y_2)\\
		&\le\left(\int_{\mathbb R^{2n}}Z_{1,2}^3\prod_{i=1}^{2}|x-y_i|^{2m-n}
		\,d\mu_g^{(2m)}(y_1)d\mu_g^{(2m)}(y_2)\right)^\frac{1}{3}\\
		&\quad\times\left(\int_{\mathbb R^{2n}}\prod_{i=1}^{2}|x-y_i|^{2m-n}
		\,d\mu_g^{(2m)}(y_1)d\mu_g^{(2m)}(y_2)	\right)^\frac{2}{3}.
	\end{aligned}
	\]
	Cubing both sides and using (\ref{eqn:I1-new-expression}) yields
	\[
	\begin{aligned}
		I_1\ge&u(x)^4\frac{\left(\displaystyle\int_{\mathbb R^{2n}}Z_{1,2}
			\prod_{i=1}^{2}|x-y_i|^{2m-n}\,d\mu_g^{(2m)}(y_1)d\mu_g^{(2m)}(y_2)	\right)^3}
		{\left(\displaystyle\int_{\mathbb R^{2n}}\prod_{i=1}^{2}|x-y_i|^{2m-n}
			\,d\mu_g^{(2m)}(y_1)d\mu_g^{(2m)}(y_2)\right)^2}.
	\end{aligned}
	\]
	By Fubini's theorem, the denominator satisfies
	\[
	\begin{aligned}
		\int_{\mathbb R^{2n}}\prod_{i=1}^{2}|x-y_i|^{2m-n}\,d\mu_g^{(2m)}(y_1)d\mu_g^{(2m)}(y_2)
		=\left(\int_{\mathbb R^n}|x-y|^{2m-n}\,d\mu_g^{(2m)}(y)\right)^2=u(x)^2.
	\end{aligned}
	\]
	Therefore,
	\begin{equation}\label{eqn:i1-inequaity-ux}
		I_1\ge\left(\int_{\mathbb R^{2n}}Z_{1,2}\prod_{i=1}^{2}|x-y_i|^{2m-n}\,d\mu_g^{(2m)}(y_1)d\mu_g^{(2m)}(y_2)
		\right)^3.
	\end{equation}
	Moreover, it follows from \(Z_{1,2}=|X(y_1)-X(y_2)|^2\), Fubini's theorem and (\ref{eqn:Xy-int-equ-0}) that
	\begin{equation}\label{eqn:Z12-2-ux-relation}
		\begin{aligned}
			&\int_{\mathbb R^{2n}}Z_{1,2}\prod_{i=1}^{2}|x-y_i|^{2m-n}
			\,d\mu_g^{(2m)}(y_1)d\mu_g^{(2m)}(y_2)\\
			&=\int_{\mathbb R^{2n}}|X(y_1)-X(y_2)|^2\prod_{i=1}^{2}|x-y_i|^{2m-n}
			\,d\mu_g^{(2m)}(y_1)d\mu_g^{(2m)}(y_2)\\
			&=
			2u(x)\int_{\mathbb R^n}
			|X(y)|^2|x-y|^{2m-n}
			\,d\mu_g^{(2m)}(y).
		\end{aligned}
	\end{equation}
	Combining this with (\ref{eqn:i1-inequaity-ux}) and (\ref{eqn:Xy-Qg2}), we arrive at
	\begin{equation}\label{eqn:I1-bounded-lower}
		\begin{aligned}
			I_1&\ge\left(2u(x)\int_{\mathbb R^n}	|X(y)|^2|x-y|^{2m-n}\,d\mu_g^{(2m)}(y)\right)^3\\
			&=8u(x)^3\left(\int_{\mathbb R^n}|X(y)|^2|x-y|^{2m-n}\,d\mu_g^{(2m)}(y)\right)^3
			=\frac{u(x)^{6+\frac{12}{n-2m}}}{8(m-1)^3}\Big(Q_g^{(2)}(x)\Big)^3.
		\end{aligned}
	\end{equation}
	
	We next estimate \(I_8\). Unlike \(I_1\), its three \(Z\)-factors involve disjoint pairs of integration variables. Hence, Fubini's theorem and (\ref{eqn:Z12-2-ux-relation}) give
	\begin{equation}\label{eqn:I8-bounded-lower}
		\begin{aligned}
			I_8&=\int_{\mathbb R^{6n}}Z_{1,2}Z_{3,4}Z_{5,6}\prod_{i=1}^{6}|x-y_i|^{2m-n}
			\,d\mu_g^{(2m)}(y)\\
			&=\left(\int_{\mathbb R^{2n}}Z_{1,2}\prod_{i=1}^{2}|x-y_i|^{2m-n}\,d\mu_g^{(2m)}(y_1)d\mu_g^{(2m)}(y_2)\right)^3\\
			&=8u(x)^3\left(\frac{u(x)^{1+\frac4{n-2m}}}{4(m-1)}Q_g^{(2)}(x)\right)^3
			=\frac{u(x)^{6+\frac{12}{n-2m}}}{8(m-1)^3}\bigl(Q_g^{(2)}(x)\bigr)^3,
		\end{aligned}
	\end{equation}
	which implies that
	\begin{equation}\label{eqn:I1-geq-I8}
		I_1\ge I_8.
	\end{equation}
	A further comparison will be useful below. By Fubini's theorem and the Cauchy-Schwarz inequality, we get
	\begin{equation}\label{eqn:I3-bounded-lower}
		\begin{aligned}
			I_3
			&\geq u(x)^2\times\frac{\displaystyle\left(\int_{\mathbb R^{2n}}Z_{1,2}
				\prod_{i=1}^{2}|x-y_i|^{2m-n}\,d\mu_g^{(2m)}(y_1)d\mu_g^{(2m)}(y_2)\right)^2}{\displaystyle\left(\int_{\mathbb R^{2n}}\prod_{i=1}^{2}|x-y_i|^{2m-n}
				\,d\mu_g^{(2m)}(y_1)d\mu_g^{(2m)}(y_2)\right)}\\
			&\quad \times\left(\int_{\mathbb R^{2n}}Z_{3,4}\prod_{i=3}^{4}|x-y_i|^{2m-n}
			\,d\mu_g^{(2m)}(y_3)d\mu_g^{(2m)}(y_4)\right)\\
			&=\left(\int_{\mathbb R^{2n}}Z_{1,2}
			\prod_{i=1}^{2}|x-y_i|^{2m-n}\,d\mu_g^{(2m)}(y_1)d\mu_g^{(2m)}(y_2)\right)^3\\
			&=I_8.
		\end{aligned}
	\end{equation}
	Meanwhile, it follows from (\ref{eqn:sumI-Q6g}) that the integral representation of the sixth-order \(Q\)-curvature takes the form
	\begin{equation}\label{eqn:Q6-section4-define}
		\begin{aligned}
			\frac{\widetilde Q_g^{(6)}(x)}{2(m-3)}
			=a_1I_1+a_2I_2+a_3I_3+a_4I_4+a_5I_5+a_6I_6+a_7I_7+a_8I_8,
		\end{aligned}
	\end{equation}
	where
	\[
	\begin{gathered}
		h=n-4m+6,\qquad A=n-2m+2,\qquad B=n-2m+4,\\
		C=n-3m+3,\qquad E=3n-8m+6,\qquad F=2n-5m+3,
	\end{gathered}
	\]
	and
	\[
	\begin{gathered}
		a_1=A^2B^2,\qquad
		a_2=-6hA^2B,\qquad
		a_3=3hCA^2,\qquad
		a_4=4hCAB,\\
		a_5=-2hA^3,\qquad
		a_6=12hCA^2,\qquad
		a_7=-6hCEA,\qquad
		a_8=hCFE.
	\end{gathered}
	\]
	The signs of these coefficients depend on the dimension. We therefore divide the remainder of the proof into three cases.
	
	\textbf{Case \((i)\). \(2m<n\le3m-3\)}.
	
	In this range, \(h=n-4m+6<0\) and \(C=n-3m+3\le0\). Since \(A,B>0\), these inequalities imply \(a_1>0\) and \(a_2,a_3,a_4,a_5,a_6\ge0\). The remaining coefficients \(a_7\) and \(a_8\) require additional consideration. However, they cannot be negative simultaneously. Indeed, the relation \(F=\frac23E+\frac{m-3}{3}\) shows that \(E>0\) implies \(F>0\). Since \(hC\ge0\), a negative \(a_7\) forces \(a_8\ge0\).
	Suppose first that \(a_7<0\). The inequalities \(I_7\le I_3\) and \(I_7\le I_6\), established in Lemma \ref{lem:I7-leq-I3-I6-I1-I4}, imply that
	\[
	\begin{aligned}
		a_3I_3+a_6I_6+a_7I_7
		&\ge(a_3+a_6+a_7)I_7\\
		&=(3hCA^2+12hCA^2-6hCEA)I_7\\
		&=3hCA(6m-n-2)I_7
		\ge0.
	\end{aligned}
	\]
	If \(a_8<0\), the inequality \(I_3\ge I_8\) gives
	\[
	\begin{aligned}
		a_3I_3+a_8I_8
		&\ge(a_3+a_8)I_8\\
		&=hC(3A^2+EF)I_8\\
		&=hC\Big[3(n-2m+2)^2+(3n-8m+6)(2n-5m+3)\Big]I_8\\
		&=hC\Big[9n^2+(33-43m)n+52m^2-78m+30\Big]I_8.
	\end{aligned}
	\]
	Together with the fact that \(\Delta=(33-43m)^2-36(52m^2-78m+30)=-23m^2-30m+9<0\)
	for \(m\ge4\), we have \(3A^2+EF>0\),
	which means \(a_3I_3+a_8I_8\ge0\).
	Therefore, whether \(a_7\) or \(a_8\) is negative, its effect can be controlled by the remaining non-negative terms. Combining this with (\ref{eqn:Q6-section4-define}) and (\ref{eqn:I1-geq-I8}), we deduce that \(\frac{\widetilde Q_g^{(6)}(x)}{2(m-3)}\ge a_1I_1=A^2B^2I_1\ge A^2B^2I_8\).
	Substituting (\ref{eqn:I8-bounded-lower}) into this estimate gives
	\[
	\begin{aligned}
		u(x)^{6+\frac{12}{n-2m}}Q_g^{(6)}(x)&=\widetilde Q_g^{(6)}(x)
		\ge2(m-3)A^2B^2I_8(x)
		=\frac{(m-3)A^2B^2}{4(m-1)^3}
		u(x)^{6+\frac{12}{n-2m}}
		\bigl(Q_g^{(2)}(x)\bigr)^3.
	\end{aligned}
	\]
	Since \(u(x)>0\), we have
	\[
	Q_g^{(6)}(x)\ge\frac{(m-3)A^2B^2}{4(m-1)^3}\bigl(Q_g^{(2)}(x)\bigr)^3.
	\]

	\textbf{Case \((ii)\). \(3m-3<n\le4m-6\)}.
	
	We first consider \(3m-3<n<4m-6\). The dimensional conditions now imply
	\[
	h<0,\qquad C>0,\qquad E>0,\qquad F>0.
	\]
	Following the rearrangement of \(P_{n,m}\) by \cite[Section 3]{Vetois-Zeitler-2026}, the integral representation of \(\widetilde Q_g^{(6)}(x)\) can be expressed as follows:
	\begin{equation}\label{eqn:two-case-tilde-Q-6-expression}
		\begin{aligned}
			\frac{\widetilde Q_g^{(6)}(x)}{2(m-3)}
			={}&A^2B^2I_1
			-(4m-n-6)C(6mA+EF)I_3\\
			&+(4m-n-6)A
			\left[6(m+1)AI_2-4CBI_4\right]\\
			&+2(4m-n-6)A^3I_5\\
			&+6(4m-n-6)CA^2(I_3+I_2-2I_6)\\
			&+3(4m-n-6)CEA(2I_7-I_3)\\
			&+(4m-n-6)CEF(I_3-I_8).
		\end{aligned}
	\end{equation}
	We use the non-negative integral combinations established in \cite[Section 3]{Vetois-Zeitler-2026}. Specifically,
	\[
	6(m+1)AI_2-4CBI_4\ge0,\quad I_3+I_2-2I_6\ge0,\quad 2I_7-I_3\ge0.
	\]
	Together this with (\ref{eqn:I3-bounded-lower}), we arrive at
	\[
	\begin{aligned}
		\frac{\widetilde Q_g^{(6)}(x)}{2(m-3)}\ge A^2B^2I_1-(4m-n-6)C(6mA+EF)I_3.
	\end{aligned}
	\]
	The coefficient of \(I_3\) is negative. To control it, we use the estimate \(I_3\le I_1\), which follows from H\"older's inequality
	\[
	\begin{aligned}
		I_3&=
		\int_{\mathbb R^{6n}}\left(Z_{1,2}^3\right)^\frac{2}{3}\left(Z_{3,4}^3\right)^\frac{1}{3}
		\prod_{i=1}^{6}|x-y_i|^{2m-n}\,d\mu_g^{(2m)}(y)\\
		&\le\left(\int_{\mathbb R^{6n}}Z_{1,2}^3\prod_{i=1}^{6}|x-y_i|^{2m-n}
		\,d\mu_g^{(2m)}(y)\right)^\frac{2}{3}\times
		\left(\int_{\mathbb R^{6n}}Z_{3,4}^3\prod_{i=1}^{6}|x-y_i|^{2m-n}
		\,d\mu_g^{(2m)}(y)\right)^\frac{1}{3}\\
		&=I_1^\frac{2}{3}I_1^\frac{1}{3}=I_1.
	\end{aligned}
	\]
	Consequently,
	\begin{equation}\label{eqn:two-case-Q6-I1}
		\begin{aligned}
			\frac{\widetilde Q_g^{(6)}(x)}{2(m-3)}
			&\ge A^2B^2I_1-(4m-n-6)C(6mA+EF)I_3\\
			&\ge\left[A^2B^2-(4m-n-6)C(6mA+EF)\right]I_1.
		\end{aligned}
	\end{equation}
	Substituting the definitions of \(A\), \(E\), and \(F\) into \(6mA+EF\) gives
	\[
	\begin{aligned}
		6mA+EF&=6m(n-2m+2)+(3n-8m+6)(2n-5m+3)\\
		&=6n^2-(25m-21)n+28m^2-42m+18.
	\end{aligned}
	\]
	On the other hand,
	\[
	\begin{aligned}
		7B^2-(6mA+EF)&=7(n-2m+4)^2-6n^2+(25m-21)n-28m^2+42m-18\\
		&=n^2-(3m-35)n-70m+94\\
		&=\int_{3m-3}^{n}(2t-3m+35)\,dt+26m-2>0,
	\end{aligned}
	\]
	since \(n>3m-3\) and \(m\ge4\). Thus,
	\begin{equation}\label{eqn:6mA-EF-7B2}
		6mA+EF<7B^2.
	\end{equation}
	Another calculation gives
	\[
	\begin{aligned}
		A^2-7(4m-n-6)C&=(n-2m+2)^2-7(4m-n-6)(n-3m+3)\\
		&=8n^2-(53m-67)n+88m^2-218m+130.
	\end{aligned}
	\]
	Viewed as a quadratic polynomial in \(n\), its discriminant is \(-7(m^2+18m-47)<0\)
	for every \(m\ge4\). Hence,
	\begin{equation}\label{eqn:A-2-7C}
		A^2-7(4m-n-6)C>0.
	\end{equation}
	Combining (\ref{eqn:6mA-EF-7B2}) and (\ref{eqn:A-2-7C}), we obtain
	\[
	\begin{aligned}
		A^2B^2-(4m-n-6)C(6mA+EF)&\ge A^2B^2-7(4m-n-6)CB^2\\
		&=B^2\left[A^2-7(4m-n-6)C\right]>0.
	\end{aligned}
	\]
	The positivity of this coefficient allows us to apply \(I_1\ge I_8\) to (\ref{eqn:two-case-Q6-I1}). Therefore,
	\[
	\begin{aligned}
		\frac{\widetilde Q_g^{(6)}(x)}{2(m-3)}
		&\ge\left[A^2B^2-(4m-n-6)C(6mA+EF)\right]I_1\\
		&\ge\left[A^2B^2-(4m-n-6)C(6mA+EF)\right]I_8\\
		&=\left[A^2B^2+hC(6mA+EF)\right]I_8\\
		&=\left[A^2B^2+hC(6mA+EF)\right]\frac{u(x)^{6+\frac{12}{n-2m}}}{8(m-1)^3}\bigl(Q_g^{(2)}(x)\bigr)^3,
	\end{aligned}
	\]
	which implies that
	\[
	\begin{aligned}
		u(x)^{6+\frac{12}{n-2m}}Q_g^{(6)}(x)&=\widetilde Q_g^{(6)}(x)
		\ge\frac{m-3}{4(m-1)^3}\left[A^2B^2+hC(6mA+EF)\right]
		u(x)^{6+\frac{12}{n-2m}}\bigl(Q_g^{(2)}(x)\bigr)^3.
	\end{aligned}
	\]
	Dividing both sides by the positive factor \(u(x)^{6+\frac{12}{n-2m}}\), we arrive at
	\[
	\begin{aligned}
		Q_g^{(6)}(x)\ge\frac{m-3}{4(m-1)^3}\left[A^2B^2+hC(6mA+EF)\right]\bigl(Q_g^{(2)}(x)\bigr)^3.
	\end{aligned}
	\]
	
	The endpoint \(n=4m-6\) can be handled directly. In this case \(h=0\), so the coefficients \(a_2,\ldots,a_8\) vanish. Formula (\ref{eqn:two-case-tilde-Q-6-expression}) consequently reduces to
	\[
	\frac{\widetilde Q_g^{(6)}(x)}{2(m-3)}
	=A^2B^2I_1\geq A^2B^2I_8=\frac{A^2B^2 u(x)^{6+\frac{12}{n-2m}}}{8(m-1)^3}\bigl(Q_g^{(2)}(x)\bigr)^3,
	\]
	This gives the same estimate as in the preceding case when \(h=0\).
	
	\textbf{Case \((iii)\). \(4m-6<n\le4m-6+[\widehat h_m]\)}.
	
	We finally turn to the dimensional range newly covered by Theorem \ref{thm:main}. The hypothesis gives \(0<h=n-4m+6\le\widehat h_m\).
	Recalling that \(G_m\) is strictly decreasing on \((0,\infty)\) and satisfies \(G_m(\widehat h_m)=0\), we know that \(G_m(h)\ge0\).
	In particular, (\ref{eqn:Qg6-geq-q1I1-I8}), (\ref{eqn:I1-geq-I8}) and (\ref{eqn:I8-bounded-lower}) give that
	\[
	\begin{aligned}
		u(x)^{6+\frac{12}{n-2m}}Q_g^{(6)}(x)&=\widetilde Q_g^{(6)}(x)\\
		&\ge
		2(m-3)\left[AG_m(h)I_1+hCFEI_8\right]\\
		&\ge 2(m-3)\left[AG_m(h)I_8+hCFEI_8\right]\\
		&=2(m-3)\left[AG_m(h)+hCFE\right]\frac{u(x)^{6+\frac{12}{n-2m}}}{8(m-1)^3}\bigl(Q_g^{(2)}(x)\bigr)^3\\
		&=\frac{(m-3)\left[AG_m(h)+hCFE\right]}{4(m-1)^3}u(x)^{6+\frac{12}{n-2m}}\bigl(Q_g^{(2)}(x)\bigr)^3.
	\end{aligned}
	\]
	Hence, 
	\[
	Q_g^{(6)}(x)\ge
	\frac{(m-3)\left[AG_m(h)+hCFE\right]}{4(m-1)^3}\bigl(Q_g^{(2)}(x)\bigr)^3.
	\]
	
	It remains to justify the strict positivity of \(Q_g^{(2)}\) when \(2m< n\le4m-6+[\widehat h_m]\). By assumption (A), the measure \(d\mu_g^{(2m)}\) has a strictly positive density on an exterior region of \(\mathbb R^n\). Since \(V(y)=\frac{x-y}{|x-y|^2}\) is not constant on any open subset of \(\mathbb R^n\setminus\{x\}\), the function \(X(y)=V(y)-\overline V(x)\) cannot vanish identically on that exterior region. Therefore,
	\[
	\int_{\mathbb R^n}|X(y)|^2|x-y|^{2m-n}\,d\mu_g^{(2m)}(y)>0.
	\]
	The representation of \(Q_g^{(2)}\) established at (\ref{eqn:Xy-Qg2}) implies \(Q_g^{(2)}(x)>0\) for every \(x\in\mathbb R^n\).
	
	Combining the estimates for the critical case \(n=2m\) and Cases \((i)\)-\((iii)\) in the non-critical range \(2m<n\le 4m-6+[\widehat h_m]\), we conclude that \(Q_g^{(6)}(x)\ge\Theta_{m,n}\bigl(Q_g^{(2)}(x)\bigr)^3\)
	for every \(x\in\mathbb R^n\) and \(2m\le n\le4m-6+[\widehat h_m]\), where
	\[
	\Theta_{m,n}:=
	\begin{cases}
		\displaystyle
		\frac{(m-3)^3}{4(m-1)^3},&n=2m;\\[3mm]
		\displaystyle\frac{(m-3)A^2B^2}{4(m-1)^3},&2m<n\le3m-3;\\[3mm]
		\displaystyle\frac{(m-3)\left[A^2B^2+hC(6mA+EF)\right]}{4(m-1)^3},&3m-3<n\le4m-6;\\[3mm]
		\displaystyle\frac{(m-3)\left[AG_m(h)+hCFE\right]}{4(m-1)^3},&4m-6<n\le4m-6+[\widehat h_m].
	\end{cases}
	\]
	The preceding estimates ensure that \(\Theta_{m,n}>0\) throughout \(2m\le n\le4m-6+[\widehat h_m]\), which, together with \(Q_g^{(2)}(x)>0\), yields
	\[
	Q_g^{(6)}(x)\ge\Theta_{m,n}\bigl(Q_g^{(2)}(x)\bigr)^3>0
	\quad\text{in }\mathbb R^n.
	\]
	This completes the proof.
	\\
	\qed

	\section{Applications of Q-Curvature Positivity to Higher-Order Elliptic Equations}
	Let \(u\in C^\infty(\mathbb R^n)\) and consider the conformally Euclidean metric
	\[
	g=
	\begin{cases}
		\displaystyle
		u^{\frac{4}{n-2m}}|dx|^2,& n>2m;\\
		\displaystyle
		e^{2u}|dx|^2,& n=2m,
	\end{cases}
	\]
	where \(u>0\) is assumed in the non-critical case \(n>2m\).
	Under the normalization adopted here, the conformal transformation formula for the highest-order \(Q\)-curvature takes the form
	\[
	Q_g^{(2m)}
	=
	\begin{cases}
		\displaystyle
		\frac{2}{n-2m}u^{-\frac{n+2m}{n-2m}}(-\Delta)^m u,& n>2m;\\
		\displaystyle
		e^{-nu}(-\Delta)^m u,& n=2m.
	\end{cases}
	\]
	For the lower-order \(Q\)-curvatures, the corresponding formulas are
	\[
	Q_g^{(2k)}
	=
	\begin{cases}
		\displaystyle
		\frac{2}{n-2k}u^{-\frac{n+2k}{n-2m}}(-\Delta)^k\left(u^{\frac{n-2k}{n-2m}}\right),
		& n>2m;\\
		\displaystyle
		\frac{2}{n-2k}
		e^{-\frac{n+2k}{2}u}
		(-\Delta)^k\left(e^{\frac{n-2k}{2}u}\right),& n=2m,
	\end{cases}
	\]
	where \(1\leq k\leq\min\{3,m-1\}\).
	In particular, the cases \(k=1,2,3\) correspond to the second-, fourth-, and sixth-order \(Q\)-curvatures, respectively. These identities allow us to translate the curvature inequalities established earlier into quantitative differential inequalities for higher-order elliptic equations.
	
	\textbf{Proof of Theorem \ref{them:main-5}}. 
	We first consider the non-critical case \(n>2m\). Consider the conformally Euclidean metric 
	\[
	g=u^{\frac4{n-2m}}|dx|^2.
	\]
	By the conformal transformation formula and the equation satisfied by \(u\),
	\[
	Q_g^{(2m)}=\frac{2}{n-2m}u^{-\frac{n+2m}{n-2m}}(-\Delta)^m u=\frac{2}{n-2m}f(x).
	\]
	Thus, the assumptions on \(f\) ensure that \(Q_g^{(2m)}\geq0\) in \(\mathbb R^n\) and that \(g\) satisfies assumption (A).
	Under the additional dimensional restriction in Theorem \ref{them:main-5} \((i)\), Theorem \(\ref{them:main-3}\) further yields
	\[
	Q_g^{(6)}\geq\Theta_{m,n}(Q_g^{(2)})^3>0.
	\]
	To express these estimates in terms of \(u\), we recall the following conformal transformation formulas:
	\[
	Q_g^{(2)}=\frac{2}{n-2}u^{-\frac{n+2}{n-2m}}(-\Delta)\left(u^{\frac{n-2}{n-2m}}\right)
	\]
	and
	\[
	Q_g^{(6)}=\frac{2}{n-6}u^{-\frac{n+6}{n-2m}}(-\Delta)^3\left(u^{\frac{n-6}{n-2m}}\right).
	\]
	Applying Theorem \(\ref{them:main-3}\) and using the above expressions, we obtain that
	\[
	\begin{aligned}
	(-\Delta)^3\left(u^{\frac{n-6}{n-2m}}\right)
	&=\frac{n-6}{2}u^{\frac{n+6}{n-2m}}Q_g^{(6)}\\
	&\geq\frac{n-6}{2}\Theta_{m,n}u^{\frac{n+6}{n-2m}}(Q_g^{(2)})^3\\
	&=\frac{4(n-6)\Theta_{m,n}}{(n-2)^3}u^{-\frac{2n}{n-2m}}
	\left[(-\Delta)\left(u^{\frac{n-2}{n-2m}}\right)\right]^3>0.
	\end{aligned}
	\]
	
	We now turn to the critical case \(n=2m\). Consider the metric \(g=e^{2u}|dx|^2\).
	In this case, the higher-order \(Q\)-curvature satisfies \(Q_g^{(n)}=e^{-nu}(-\Delta)^m u=f(x)\).
	Hence, \(Q_g^{(n)}\geq0\) in \(\mathbb R^n\) and \(Q_g^{(n)}\not\equiv0\). Moreover, the second-order conformal transformation formula gives
	\[
	Q_g^{(2)}=-e^{-2u}\left(\Delta u+\frac{n-2}{2}|\nabla u|^2\right)=\frac2{n-2}e^{-\frac{n+2}{2}u}(-\Delta)e^{\frac{n-2}{2}u}.
	\]
	The assumed differential inequality therefore implies that \(Q_g^{(2)}\geq0\) outside a compact set. Thus, assumption (A) is satisfied. Since \(n=2m\), the dimensional condition in Theorem \(\ref{them:main-3}\) is satisfied. Hence,
	\[
	Q_g^{(6)}\geq\Theta_{m,n}(Q_g^{(2)})^3>0.
	\]
	Using
	\(Q_g^{(6)}=\frac2{n-6}e^{-\frac{n+6}{2}u}(-\Delta)^3e^{\frac{n-6}{2}u}\),
	we derive that
	\[
	\begin{aligned}
	(-\Delta)^3e^{\frac{n-6}{2}u}&=\frac{n-6}{2}e^{\frac{n+6}{2}u}Q_g^{(6)}\\
	&\geq\frac{n-6}{2}\Theta_{m,n}e^{\frac{n+6}{2}u}(Q_g^{(2)})^3\\
	&=\frac{4(n-6)\Theta_{m,n}}{(n-2)^3}e^{-nu}
	\left[(-\Delta)e^{\frac{n-2}{2}u}\right]^3>0.
	\end{aligned}
	\]
	This completes the proof.
	\\
	\qed
	
	\section{Counterexamples to Sixth-Order Q-Curvature Positivity}
\begin{figure}[H]
		\centering
		\includegraphics[width=0.8\textwidth]{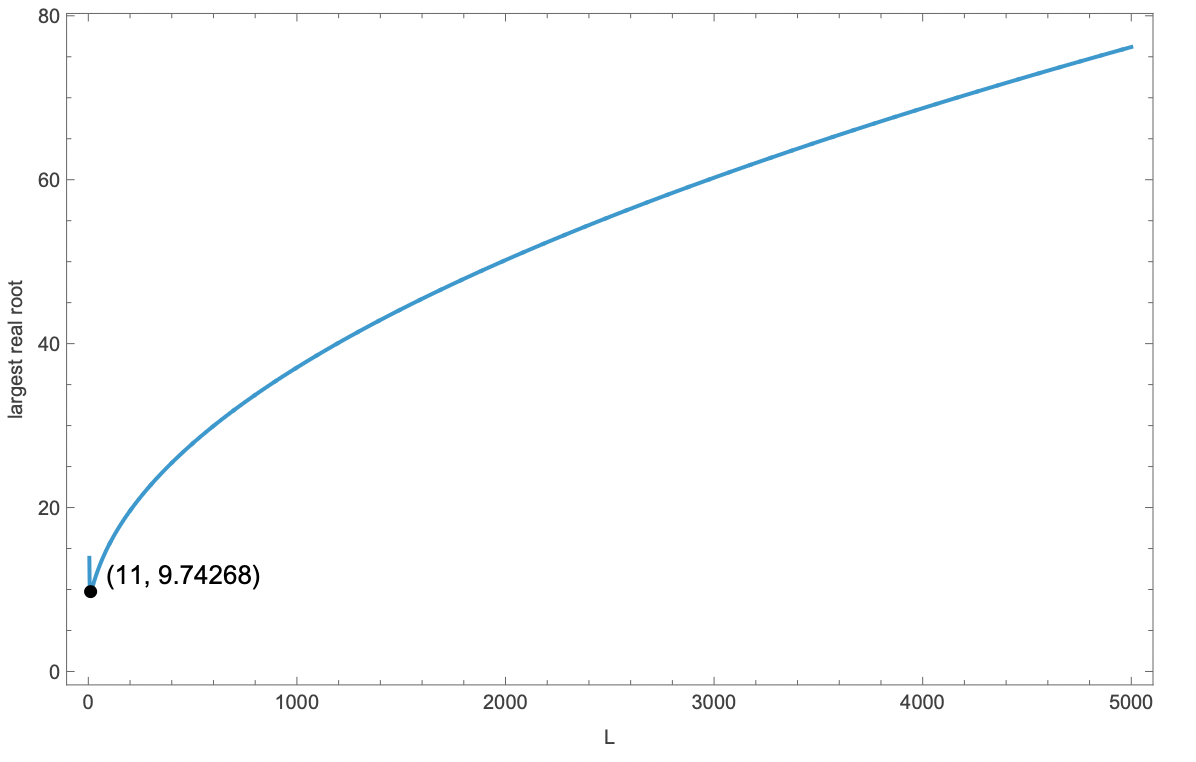}
		\caption{The largest real root \(\Lambda_L\) of \(P_L\) for \(5\le L\le5000\).}
		\label{fig:h_comparison-4}
	\end{figure}

\begin{figure}[H]
		\centering
		\includegraphics[width=0.8\textwidth]{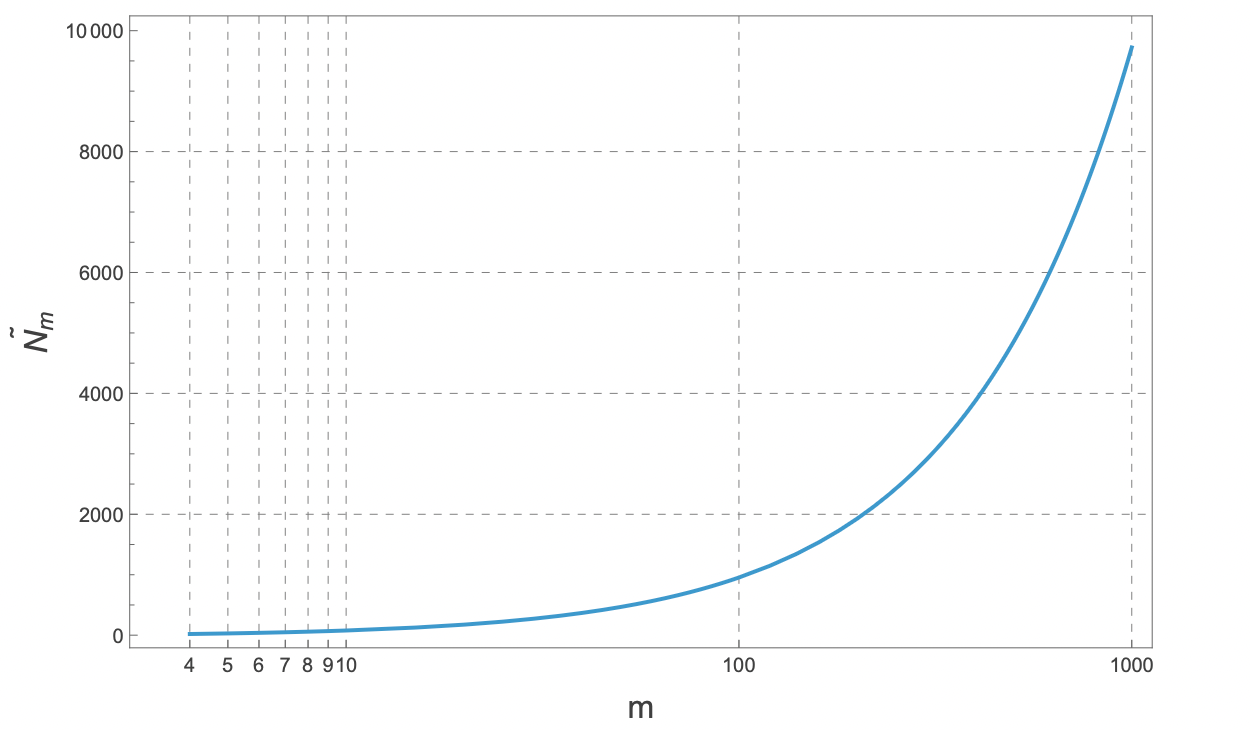}
		\caption{The largest real root \(\widetilde N_m\) of \(H_m\) for \(4\le m\le1000\) on a logarithmic horizontal scale.}
		\label{fig:h_comparison-5}
	\end{figure}

	Throughout this section, \(m,n\in\mathbb{N}\) satisfy \(n>2m\geq8\). The construction of \cite[Theorem 1.2]{Vetois-Zeitler-2026} uses two concentration points with different heights. The relative height of the two peaks plays an essential role in their arguments. Here we keep the same potential-theoretic framework but use a different mechanism to produce non-positivity. Instead of varying the relative heights of two peaks, we place \(L\) identical peaks of equal mass and height at mutually orthogonal unit vectors. This choice gives a particularly simple geometry. Two entries of a sampled six-tuple are either the same center or two distinct centers at distance \(\sqrt2\). Therefore, the limiting integrals can be reduced to counting the index configurations that make the corresponding factors nonzero. Set
	\[
	c_{n,m}:=\frac{\Gamma\left(\frac n2-m+1\right)}{2^{2m-1}(m-1)!
		\Gamma\left(\frac n2\right)\omega_{n-1}}.
	\]
	Fix \(L\in\mathbb N\) with \(2\le L\le n\). For \(j=1,\ldots,L\), let \(e_j\) be the unit vector in the \(j\)-th coordinate direction of \(\mathbb R^n\), namely,
	\[
	e_j=(0,\ldots,0,\underbrace{1}_{j\text{-th component}},0,\ldots,0).
	\]
	Following  \cite[Section 4]{Vetois-Zeitler-2026}, choose \(\rho_\varepsilon\in C_c^\infty(\mathbb R^n)\) for \(0<\varepsilon<\frac{1}{2}\) such that \(\rho_\varepsilon\ge0\), \(\supp\rho_\varepsilon\subset B(0,\varepsilon)\), \(\int_{\mathbb R^n}\rho_\varepsilon(x)\,dx=1,\) and \(\rho_\varepsilon\mathcal L\overset{*}{\rightharpoonup}
	\delta_0\) weakly in the sense of measures as \(\varepsilon\to0\), where \(\mathcal L\) is the Lebesgue measure in \(\mathbb R^n\). Define
	\begin{equation}\label{eqn:h-L-varepsilon-define}
		h_{L,\varepsilon}(x):=\varepsilon(1+|x|^2)^{-\frac{n+2m}{2}}
		+\sum_{j=1}^{L}\rho_\varepsilon(x-e_j).
	\end{equation}
	Since \(\rho_\varepsilon\mathcal L\) converges weakly to \(\delta_0\), translating the measure by \(e_j\) gives
	\[
	\rho_\varepsilon(x-e_j)\mathcal L
	\overset{*}{\rightharpoonup}
	\delta_{e_j}
	\]
	as \(\varepsilon\to 0\) weakly in the sense of measures for each \(j=1,\ldots,L\), where \(\delta_{e_j}\) is the Dirac measure at \(e_j\). At the same time, the first term in \((\ref{eqn:h-L-varepsilon-define})\) converges to zero in \(L^1(\mathbb R^n)\) as \(\varepsilon\to0\). Hence,
	\begin{equation}\label{eqn:hLepsilon-weak-delta}
		h_{L,\varepsilon}\mathcal L\overset{*}{\rightharpoonup}
		\sum_{j=1}^{L}\delta_{e_j}
	\end{equation}
	as \(\varepsilon\to0\) weakly in the sense of measures. Now let
	\[
	u_{L,\varepsilon}(x):=c_{n,m}\int_{\mathbb R^n}|x-y|^{2m-n}
	h_{L,\varepsilon}(y)\,dy.
	\]
	As shown in \cite[Section 4]{Vetois-Zeitler-2026}, \(u_{L,\varepsilon}\) is smooth and satisfies
	\[
	(-\Delta)^m u_{L,\varepsilon}=\frac{n-2m}{2}h_{L,\varepsilon}
	\quad\text{in }\mathbb R^n.
	\]
	Moreover, there exists a constant \(C_{L,\varepsilon}>1\) depending on \(n,m,L\), and \(\varepsilon\), such that
	\[
	C_{L,\varepsilon}^{-1}(1+|x|^2)^{\frac{2m-n}{2}}
	\le u_{L,\varepsilon}(x)
	\le C_{L,\varepsilon}(1+|x|^2)^{\frac{2m-n}{2}}
	\]
	for all \(x\in\mathbb R^n\). Let \(g_{L,\varepsilon}\) be the metric on \(\mathbb R^n\) defined by
	\[
	g_{L,\varepsilon}:=\phi_{L,\varepsilon}^{2}|dx|^2,
	\]
	where
	\[
	\phi_{L,\varepsilon}:=u_{L,\varepsilon}^{\frac{2}{n-2m}}.
	\]
	Then 
	\[
	\begin{aligned}
		Q_{g_{L,\varepsilon}}^{(2m)}&=\frac{2}{n-2m}\phi_{L,\varepsilon}^{-\frac{n+2m}{2}}(-\Delta)^m(\phi_{L,\varepsilon}^\frac{n-2m}{2})\\
		&=\frac{2}{n-2m}u_{L,\varepsilon}^{-\frac{n+2m}{2}\cdot\frac{2}{n-2m}}(-\Delta)^mu_{L,\varepsilon}\\
		&=\frac{2}{n-2m}u_{L,\varepsilon}^{-\frac{n+2m}{n-2m}}\cdot \frac{n-2m}{2}h_{L,\varepsilon}\\
		&=u_{L,\varepsilon}^{-\frac{n+2m}{n-2m}}h_{L,\varepsilon}
	\end{aligned}
	\]
	satisfies
	\begin{equation}\label{eqn:Q-g-varep-estimate}
		\widetilde{C}_{L,\varepsilon}^{-1}\le Q_{g_{L,\varepsilon}}^{(2m)}\le
		\widetilde{C}_{L,\varepsilon}
	\end{equation}
	for  some \(\widetilde{C}_{L,\varepsilon}>0\) depending on \(n,m,L\) and \(\varepsilon\). In particular, \(g_{L,\varepsilon}\) satisfies assumption (A) for every \(0<\varepsilon<\frac{1}{2}\). Let
	\(I_1^{(\varepsilon)},\ldots,I_8^{(\varepsilon)}\) denote the eight integrals in Proposition \ref{pro:vetois-Zeitler-result-pro2.1} corresponding to \(g_{L,\varepsilon}\). Then
	\[
	\frac{\widetilde{Q}_{g_{L,\varepsilon}}^{(6)}(0)}
	{2(m-3)}=\sum_{k=1}^{8}a_k I_k^{(\varepsilon)},
	\]
	where \(a_1,\ldots,a_8\) are the coefficients introduced in (\ref{eqn:sumI-Q6g}). We now determine the limit of each \(I_k^{(\varepsilon)}\) (\(k=1,\dots,8\)) as \(\varepsilon\to0\). By the definition of \(d\mu_g^{(2m)}\) and the expression for \(Q_{g_{L,\varepsilon}}^{(2m)}\), we have
	\[
	\begin{aligned}
		d\mu_{g_{L,\varepsilon}}^{(2m)}(y)
		&=c_{n,m}Q_{g_{L,\varepsilon}}^{(2m)}(y)u_{L,\varepsilon}(y)^{\frac{n+2m}{n-2m}}\,dy\\
		&=c_{n,m}u_{L,\varepsilon}(y)^{-\frac{n+2m}{n-2m}}h_{L,\varepsilon}(y)u_{L,\varepsilon}(y)^{\frac{n+2m}{n-2m}}\,dy\\
		&=c_{n,m}h_{L,\varepsilon}(y)\,dy.
	\end{aligned}
	\]
	Moreover, it follows from (\ref{eqn:hLepsilon-weak-delta}) that
	\[
	d\mu_{g_{L,\varepsilon}}^{(2m)}\overset{*}{\rightharpoonup}c_{n,m}\sum_{j=1}^{L}\delta_{e_j}
	\]
	as \(\varepsilon \to 0\). To justify the passage to the limit in the sixfold integrals at \(x=0\), we first note that the translated peaks are supported away from the origin. For \(0\le q\le3\), we have
	\[
	\varepsilon\int_{\mathbb R^n}|y|^{2m-n-2q}(1+|y|^2)^{-\frac{n+2m}{2}}\,dy
	=O(\varepsilon).
	\]
	The condition \(2m\ge8>2q\) ensures integrability near the origin, while the radial integrand decays like \(r^{-n-2q-1}\) at infinity. The corresponding weighted integrals over the translated peaks are uniformly bounded in \(\varepsilon\), since their supports remain away from the origin and each peak has unit mass. Together with
	\[
	Z_{i,j}\le2\left(|y_i|^{-2}+|y_j|^{-2}\right),
	\]
	these estimates show that all terms containing at least one background factor vanish as \(\varepsilon\to0\). The remaining terms involve only the translated peaks. Since their supports lie in a fixed compact set away from the origin, the corresponding integrands are uniformly continuous. Consequently, the limiting values of \(I_1^{(\varepsilon)},\ldots,I_8^{(\varepsilon)}\) are determined by the points in \(\{e_1,\ldots,e_L\}^{6}\). To describe these points, take
	\[
	\alpha:=(\alpha_1,\ldots,\alpha_6)\in\{1,\ldots,L\}^{6}
	\]
	and write
	\[
	e_\alpha:=(e_{\alpha_1},\ldots,e_{\alpha_6}).
	\]
	For every such \(\alpha\), each component \(e_{\alpha_i}\) is one of the unit vectors \(e_1,\ldots,e_L\). Therefore, for every \(i=1,\ldots,6\), we have
	\[
	|e_{\alpha_i}|=1
	\]
	and
	\[
	|e_{\alpha_i}|^{2m-n}=1.
	\] 
	The only nontrivial factors that remain are the terms \(Z_{i,j}\), and their values are determined by the mutual orthogonality of \(e_1,\ldots,e_L\). More precisely,
	\[
	\begin{aligned}
		Z_{i,j}(0,e_\alpha)&=\left|Y_{i,j}(0,e_\alpha)\right|^2\\
		&=\left|V(e_{\alpha_i})-V(e_{\alpha_j})\right|^2\\
		&=\left|\frac{0-e_{\alpha_i}}{|0-e_{\alpha_i}|^2}-\frac{0-e_{\alpha_j}}{|0-e_{\alpha_j}|^2}\right|^2\\
		&=\left|-\frac{e_{\alpha_i}}{|e_{\alpha_i}|^2}+\frac{e_{\alpha_j}}{|e_{\alpha_j}|^2}\right|^2\\
		&=\left|-e_{\alpha_i}+e_{\alpha_j}\right|^2\\
		&=|e_{\alpha_i}|^2+|e_{\alpha_j}|^2-2e_{\alpha_i}\cdot e_{\alpha_j}\\
		&=
		\begin{cases}
			0,& \alpha_i=\alpha_j;\\[2mm]
			2,& \alpha_i\ne\alpha_j.
		\end{cases}
	\end{aligned}
	\]
	Therefore, a product of three \(Z_{i,j}\)-terms can take only the values \(0\) and \(8\). For the metric \(g_{L,\varepsilon}\), the eight integrals in Proposition \ref{pro:vetois-Zeitler-result-pro2.1} are evaluated at \(x=0\). Thus,
	\[
	I_1^{(\varepsilon)}=\int_{\mathbb R^{6n}}Z_{1,2}^{3}\prod_{i=1}^{6}|y_i|^{2m-n}
	\,d\mu_{g_{L,\varepsilon}}^{(2m)}(y).
	\]
	Passing to the limit gives
	\[
	\begin{aligned}
		\lim_{\varepsilon\to0} I_1^{(\varepsilon)}
		&=\lim_{\varepsilon\to0}\int_{\mathbb R^{6n}}Z_{1,2}^{3}(0,y)
		\prod_{i=1}^{6}|y_i|^{2m-n}\,d\mu_{g_{L,\varepsilon}}^{(2m)}(y)\\
		&=c_{n,m}^{6}\sum_{\alpha_1=1}^{L}\cdots\sum_{\alpha_6=1}^{L}Z_{1,2}^{3}(0,e_\alpha)\prod_{i=1}^{6}|e_{\alpha_i}|^{2m-n}\\
		&=c_{n,m}^{6}\sum_{\alpha_1=1}^{L}\cdots\sum_{\alpha_6=1}^{L}Z_{1,2}^{3}(0,e_\alpha)\\
		&=8c_{n,m}^{6}\sum_{\substack{\alpha_1,\ldots,\alpha_6\in\{1,\ldots,L\}\\\alpha_1\ne\alpha_2}}1\\
		&=8c_{n,m}^{6}L^5(L-1).
	\end{aligned}
	\]
	Applying the same calculation to \(I_2^{(\varepsilon)},\ldots,I_8^{(\varepsilon)}\), we obtain that
	\[
	\begin{aligned}
		\lim_{\varepsilon\to0} I_2^{(\varepsilon)}
		&=\lim_{\varepsilon\to0}\int_{\mathbb R^{6n}}Z_{1,2}^{2}(0, y)Z_{1,3}(0,y)\prod_{i=1}^{6}|y_i|^{2m-n}\,d\mu_{g_{L,\varepsilon}}^{(2m)}( y)	
		=8c_{n,m}^{6}L^4(L-1)^2,
	\end{aligned}
	\]
	\[
	\begin{aligned}
		\lim_{\varepsilon\to0} I_3^{(\varepsilon)}
		&=\lim_{\varepsilon\to0}\int_{\mathbb R^{6n}}Z_{1,2}^{2}(0, y)Z_{3,4}(0,y)\prod_{i=1}^{6}|y_i|^{2m-n}\,d\mu_{g_{L,\varepsilon}}^{(2m)}( y)
	=8c_{n,m}^{6}L^4(L-1)^2,
	\end{aligned}
	\]
	\[
	\begin{aligned}
		\lim_{\varepsilon\to0} I_4^{(\varepsilon)}
		&=\lim_{\varepsilon\to0}\int_{\mathbb R^{6n}}Z_{1,2}(0, y)Z_{1,3}(0, y)Z_{1,4}(0, y)\prod_{i=1}^{6}|y_i|^{2m-n}
		\,d\mu_{g_{L,\varepsilon}}^{(2m)}(y)
		=8c_{n,m}^{6}L^3(L-1)^3,
	\end{aligned}
	\]
	\[
	\begin{aligned}
		\lim_{\varepsilon\to0} I_5^{(\varepsilon)}
		&=\lim_{\varepsilon\to0}\int_{\mathbb R^{6n}}Z_{1,2}(0,y)Z_{1,3}(0,y)Z_{2,3}(0,y)\prod_{i=1}^{6}|y_i|^{2m-n}\,d\mu_{g_{L,\varepsilon}}^{(2m)}(y)
		=8c_{n,m}^{6}L^4(L-1)(L-2),
	\end{aligned}
	\]
	\[
	\begin{aligned}
		\lim_{\varepsilon\to0} I_6^{(\varepsilon)}
		&=\lim_{\varepsilon\to0}\int_{\mathbb R^{6n}}Z_{1,2}(0,y)Z_{1,3}(0,y)Z_{3,4}(0,y)\prod_{i=1}^{6}|y_i|^{2m-n}\,d\mu_{g_{L,\varepsilon}}^{(2m)}( y)
		=8c_{n,m}^{6}L^3(L-1)^3,
	\end{aligned}
	\]
	\[
	\begin{aligned}
		\lim_{\varepsilon\to0} I_7^{(\varepsilon)}
		&=\lim_{\varepsilon\to0}\int_{\mathbb R^{6n}}Z_{1,2}(0,y)Z_{1,3}(0,y)Z_{5,6}(0, y)\prod_{i=1}^{6}|y_i|^{2m-n}\,d\mu_{g_{L,\varepsilon}}^{(2m)}(y)
		=8c_{n,m}^{6}L^3(L-1)^3,
	\end{aligned}
	\]
	and
	\[
	\begin{aligned}
		\lim_{\varepsilon\to0} I_8^{(\varepsilon)}
		&=\lim_{\varepsilon\to0}\int_{\mathbb R^{6n}}
		Z_{1,2}(0,y)Z_{3,4}(0,y)Z_{5,6}(0,y)\prod_{i=1}^{6}|y_i|^{2m-n}
		\,d\mu_{g_{L,\varepsilon}}^{(2m)}(y)
		=8c_{n,m}^{6}L^3(L-1)^3.
	\end{aligned}
	\]
	Combining these eight limits with Proposition \ref{pro:vetois-Zeitler-result-pro2.1} yields
	\begin{equation}\label{eqn:lim-0-Q-6-epsilon}
		\begin{aligned}
			\lim_{\varepsilon\to0}
			\widetilde{Q}_{g_{L,\varepsilon}}^{(6)}(0)
			&=2(m-3)\lim_{\varepsilon\to0}\Big(a_1 I_1^{(\varepsilon)}+a_2 I_2^{(\varepsilon)}+a_3 I_3^{(\varepsilon)}+a_4 I_4^{(\varepsilon)}+a_5 I_5^{(\varepsilon)}+a_6 I_6^{(\varepsilon)}\\
			&\quad+a_7 I_7^{(\varepsilon)}+a_8 I_8^{(\varepsilon)}\Big)\\
			&=16(m-3)c_{n,m}^{6}L^3(L-1)
			\Big[L^2a_1+L(L-1)(a_2+a_3)\\
			&\quad+L(L-2)a_5+(L-1)^2(a_4+a_6+a_7+a_8)
			\Big]\\
			&=16(m-3)c_{n,m}^{6}L^3(L-1)\Phi_{m,L}(n),
		\end{aligned}
	\end{equation}
	
	where
	\[
	\Phi_{m,L}(n):=L^2a_1+L(L-1)(a_2+a_3)+L(L-2)a_5+(L-1)^2(a_4+a_6+a_7+a_8).
	\]
	Since \(16(m-3)c_{n,m}^{6}L^3(L-1)>0\), the sign of the limiting value is determined by \(\Phi_{m,L}(n)\). Therefore, the non-positivity problem reduces to finding a suitable value of \(L\) such that \(\Phi_{m,L}(n)<0\). A direct expansion of \(\Phi_{m,L}(n)\) yields
	\[
	\begin{aligned}
		\Phi_{m,L}(n)&=(4-L)n^4+\Big[(11L-39)m-21L+57\Big]n^3\\
		&\quad+\Big[(L^2-38L+133)m^2+(-6L^2+162L-396)m+9L^2-148L+287\Big]n^2\\
		&\quad+4\Big[(11L-47)m^3+(3L^2-87L+216)m^2+(-15L^2+175L-320)m\\
		&\quad+14L^2-103L+153\Big]n\\
		&\quad+4(m-1)\Big[(24-4L)m^3+(50L-126)m^2+(13L^2-130L+213)m\\
		&\quad-19L^2+96L-117\Big].
	\end{aligned}
	\]
	For each fixed integer \(L\geq5\), the leading coefficient of \(\Phi_{m,L}(n)\) is \(4-L<0\). Hence, \(\Phi_{m,L}(n)<0\) for all sufficiently large \(n\), with \(n\geq L\) so that the construction is
	available.
	In view of
	\[
	\lim_{\varepsilon\to0}
	\widetilde{Q}_{g_{L,\varepsilon}}^{(6)}(0)
	=16(m-3)c_{n,m}^{6}L^3(L-1)\Phi_{m,L}(n),
	\]
	we have
	\[
	\lim_{\varepsilon\to0}\widetilde{Q}_{g_{L,\varepsilon}}^{(6)}(0)<0
	\]
	for sufficiently large \(n\). To compare the asymptotic non-positivity thresholds for different values of \(L\), we fix an integer \(L\ge5\) and set \(n=sm\), where \(s>2\). Substituting this into \(\Phi_{m,L}(n)\) and dividing by \(m^4\) yields
	\[
	\begin{aligned}
		\lim_{m\to\infty}\frac{\Phi_{m,L}(sm)}{m^4}
		&=(4-L)s^4+(11L-39)s^3+(L^2-38L+133)s^2+(44L-188)s+
		96-16L\\
		&=:P_L(s).
	\end{aligned}
	\]
	Since \(L\ge5\), the leading coefficient of \(P_L\) is negative. Thus,
	\[
	P_L(s)\rightarrow-\infty
	\quad
	\text{as }s\to\infty.
	\]Let \(\Lambda_L\) denote the largest real root of \(P_L\). Then
	\[
	P_L(s)<0
	\quad
	\text{for all }s>\Lambda_L.
	\]
	Hence \(\Lambda_L\) determines the asymptotic non-positivity threshold. To obtain the smallest such threshold, it remains to compare \(\Lambda_L\) for integers \(L\ge5\). In the following, we examine \(\Lambda_L\) numerically. Figure \ref{fig:h_comparison-4} plots \(\Lambda_L\) for \(5\le L\le5000\), and the following numerical results show that the minimum occurs at \(L=11\). 
	At the same time, for \(L=11\),
	\[
	P_{11}(s)=-7s^4+82s^3-164s^2+296s-80.
	\]
	Since
	\[
	P_{11}(9)=3151>0,
	\]
	and
	\[
	P_{11}(10)=-1520<0,
	\]
	there exists a real root \(\Lambda_{11}\in(9,10)\). Note that
	\[
	P_{11}'(s)=-28s^3+246s^2-328s+296
	\]and
	\[
	P_{11}''(s)=-4(21s^2-123s+82).
	\]For \(s\ge9\), the quadratic \(21s^2-123s+82\) is positive. Hence, \(P_{11}''(s)<0\)
   for \(s\ge9\). Thus, \(P_{11}'\) is strictly decreasing on \([9,\infty)\). Since
	\[
	P_{11}'(9)=-3142<0,
	\]
	it follows that \(P_{11}'(s)<0\) for \(s\ge9\).
	Therefore, \(P_{11}(s)\) is strictly decreasing on \([9,\infty)\), and the root \(\Lambda_{11}\) is the unique root in this interval and no larger real root can appear. Hence, \(\Lambda_{11}\) is the largest real root of \(P_{11}\). Moreover,
	\({\Lambda_{11}}\) is also a simple root. To obtain a numerical approximation of \(\Lambda_{11}\), we use mathematical software to compute this root, which gives
	\(\Lambda_{11}\simeq 9.74268\).
	The above numerical results suggest that \(L=11\) is the optimal choice when \(5\le L\le5000\) and the largest real root of \(P_{11}(s)\) is \(\Lambda_{11}\simeq 9.74268\). Motivated by this observation, we next prove that \(\Lambda_{11}\le \Lambda_L\)
	for every integer \(L\ge5\). Set
	\[
	\widetilde\Lambda:=\Lambda_{11}.
	\]
	Thus, \(\widetilde\Lambda\simeq 9.74268\) is the largest real root of \(P_{11}(s)\). Since \(P_{11}(\widetilde\Lambda)=0\),
	a direct calculation gives
	\[
	P_L(\widetilde\Lambda)=\widetilde\Lambda^2(L-11)(L-L_0),
	\]
	where
	\[
	L_0=\frac{\widetilde\Lambda^4-11\widetilde\Lambda^3+27\widetilde\Lambda^2-44\widetilde\Lambda+16}{\widetilde\Lambda^2}\simeq 10.40269.
	\]
	Since \(10<L_0<11\), we have
	\[
	P_L(\widetilde\Lambda)=\widetilde\Lambda^2(L-11)(L-L_0)\ge0
	\]
	for every integer \(L\ge5\) and the equality holds only for \(L=11\). On the other hand, since \(L\ge5\),
	\[
	P_L(s)\rightarrow-\infty\quad\text{as }s\to\infty.
	\]
	Hence, for every integer \(L\ge5\), \(P_L\) has a real root in \([\widetilde\Lambda,\infty)\). Therefore, \(\Lambda_L\ge\widetilde\Lambda\).
	If \(L\ne11\), then \(P_L(\widetilde\Lambda)>0\), so the corresponding largest real root satisfies \(\Lambda_L>\widetilde\Lambda\).
	Consequently,
	\[
	\Lambda_{11}=\min_{L\in\mathbb N,\;L\ge5}\Lambda_L=\widetilde\Lambda\simeq 9.74268.
	\]
	Thus, we choose \(L=11\) and define
	\[
	\begin{aligned}
		H_m(n)&:=\Phi_{m,11}(n)\\
		&=-7n^4+(82m-174)n^3+(-164m^2+660m-252)n^2\\
		&\quad+(296m^3-1512m^2-840m+2856)n\\
		&\quad-80m^4+1776m^3-272m^2-6864m+5440.
	\end{aligned}
	\]
	Accordingly, (\ref{eqn:lim-0-Q-6-epsilon}) reduces to
	\begin{equation}\label{eqn:lim-Q11var-Hm}
		\lim_{\varepsilon\to0}\widetilde{Q}_{g_{11,\varepsilon}}^{(6)}(0)
		=16\cdot11^3\cdot10(m-3)c_{n,m}^{6}H_m(n)=212960(m-3)c_{n,m}^{6}H_m(n).
	\end{equation}
	Therefore, the sign of the limiting \(\widetilde{Q}_{g_{11,\varepsilon}}^{(6)}(0)\) is completely determined by the polynomial \(H_m\). A direct calculation at \(n=4m-6\) gives
	\[
	H_m(4m-6)=1936(m-2)^2(m-1)^2>0.
	\]
	In contrast, the leading term of \(H_m(n)\) is \(-7n^4\). Hence \(H_m(n)\rightarrow-\infty\) as \(n\to\infty\).
	By continuity, \(H_m\) must possess at least one real root larger than \(4m-6\). Let \(\widetilde N_m\) denote the largest real root of \(H_m\). Then, for every \(n>\widetilde N_m\), we have
	\begin{equation}\label{eqn:Hm-le-0-n-ge-Nm}
		H_m(n)<0.
	\end{equation}
	For convenience, we show the dependence of \(\widetilde N_m\) on \(m\) in Figure \ref{fig:h_comparison-5} using a logarithmic horizontal scale.
	If \(n\in\mathbb N\) and \(n>\widetilde N_m\), we have \(n>4m-6\ge10\).
	Hence, \(n\ge11\) and the first eleven standard basis vectors \(e_1,\ldots,e_{11}\) of \(\mathbb R^n\) are well defined and form a family of mutually orthogonal unit vectors. For such an integer \(n\), (\ref{eqn:lim-Q11var-Hm}) and (\ref{eqn:Hm-le-0-n-ge-Nm}) yield
	\[
	\lim_{\varepsilon\to0}
	\widetilde{Q}_{g_{11,\varepsilon}}^{(6)}(0)<0.
	\]
	Since \(u_{11,\varepsilon}(0)>0\), this implies that \({Q}_{g_{11,\varepsilon}}^{(6)}(0)<0\)
	for sufficiently small \(\varepsilon>0\). For each such \(\varepsilon\), estimate (\ref{eqn:Q-g-varep-estimate}) with \(L=11\) gives a constant \(C_\varepsilon>0\) such that
	\[
	C_\varepsilon^{-1}\le Q_{g_{11,\varepsilon}}^{(2m)}\le C_\varepsilon\quad\text{in }\mathbb R^n.
	\]
	Thus, for sufficiently small \(\varepsilon>0\), the metric \(g_{11,\varepsilon}\) satisfies assumption (A), while
	\({Q}_{g_{11,\varepsilon}}^{(6)}(0)<0\). 
	
\begin{figure}[H]
		\centering
		\includegraphics[width=0.8\textwidth]{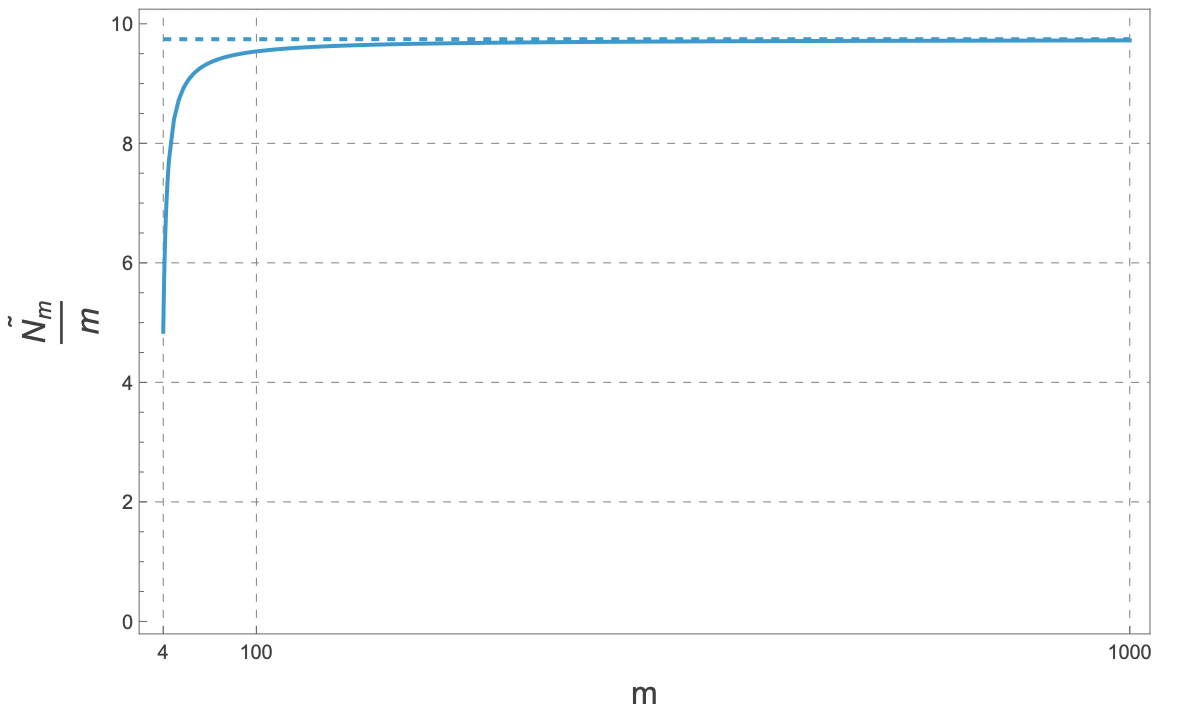}
		\caption{The ratio \(\frac{\widetilde N_m}{m}\) as a function of \(m\).}
		\label{fig:h_comparison-6}
	\end{figure}

	Set
	\[
	s:=\frac{n}{m}.
	\]
	Substituting \(n=ms\) into \(H_m(n)\) and dividing by \(m^4\), we obtain
	\[
	\begin{aligned}
		\frac{H_m(ms)}{m^4}&=-7s^4+82s^3-164s^2+296s-80+
		\frac{-174s^3+660s^2-1512s+1776}{m}\\
		&\quad+\frac{-252s^2-840s-272}{m^2}+\frac{2856s-6864}{m^3}+\frac{5440}{m^4}.
	\end{aligned}
	\]
	Therefore,
	\[
	\frac{H_m(ms)}{m^4}\rightarrow P_{11}(s)
	=-7s^4+82s^3-164s^2+296s-80
	\]
	locally uniformly in \(s\) as \(m\to\infty\). Recall that \(\widetilde\Lambda=\Lambda_{11}\) is the largest real root of \(P_{11}(s)\) with
	\(\widetilde\Lambda\simeq 9.74268\), and \(\widetilde{\Lambda}\) is also a simple root. Let \(s_m\) denote the largest real root of \(\frac{H_m(ms)}{m^4}\). Since the coefficients of \(\frac{H_m(ms)}{m^4}\) converge to the corresponding coefficients of \(P_{11}(s)\) as \(m\to\infty\), the roots vary continuously with the coefficients. Hence, \(s_m\rightarrow \widetilde{\Lambda}\) as \(m\to\infty\).
	On the other hand, if \(\widetilde N_m\) denotes the largest real root of \(H_m(n)\), then the change of variables \(n=ms\) gives \(s_m=\frac{\widetilde N_m}{m}\).
	Consequently,
	\[
	{\widetilde N_m}\sim \widetilde{\Lambda}m
	\quad\text{as }m\to\infty.
	\]
	and \(\widetilde{\Lambda}\simeq 9.74268\).
	To further illustrate this asymptotic behavior, we show the variation of \(\frac{\widetilde N_m}{m}\) with \(m\) in Figure \ref{fig:h_comparison-6} and provide the corresponding numerical values of \(\widetilde N_m\) and \(\frac{\widetilde N_m}{m}\) for selected values of \(m\) in Table \ref{fig:h_comparison-7}.

	\begin{table}[H]
		\centering
		\includegraphics[width=0.4\textwidth]{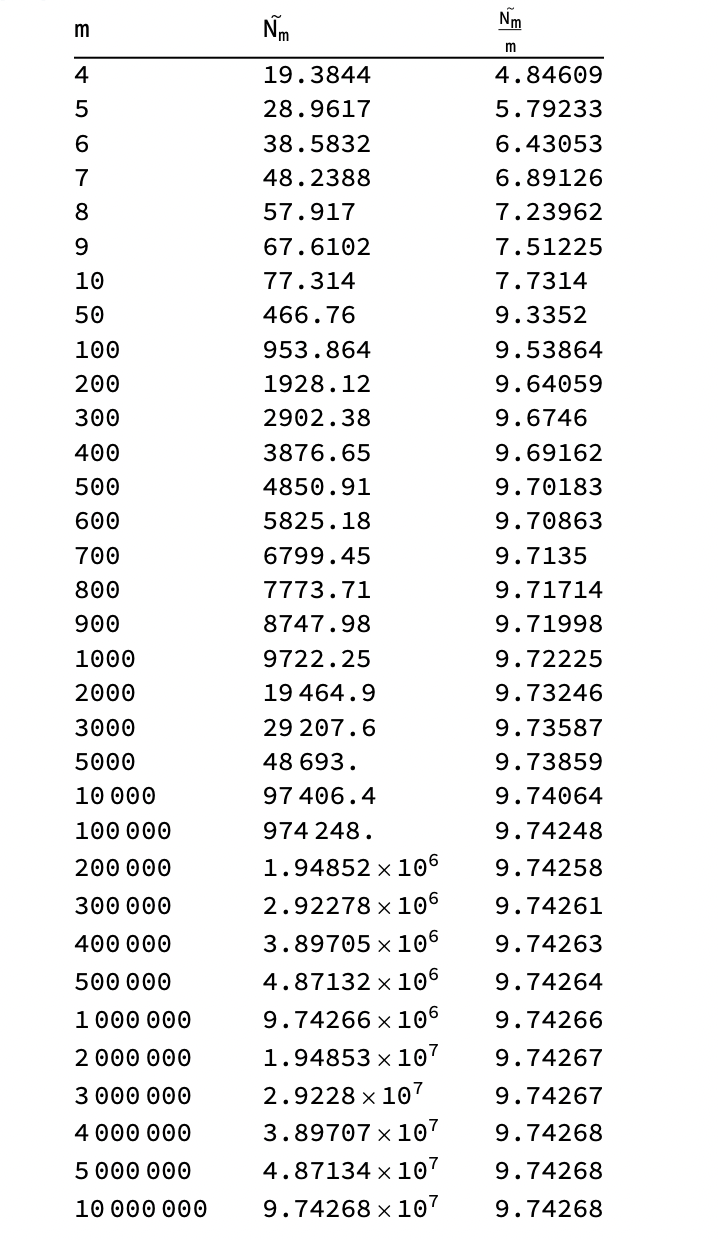}
		\caption{Numerical values of \(\widetilde N_m\) and \(\frac{\widetilde N_m}{m}\) for selected values of \(m\).}
		\label{fig:h_comparison-7}
	\end{table}

	Finally, we apply the stereographic projection argument from \cite[Section 4]{Vetois-Zeitler-2026}. For every \(0<\varepsilon<\frac{1}{2}\), the metric \(g_{11,\varepsilon}\) extends smoothly to a metric on \(\mathbb S^n\) that is conformal to the standard metric. This completes the proof of Theorem \ref{thm:main-2}.
	\qed
	\\

	\noindent\textbf{AI Disclosure.} The authors acknowledge the use of ChatGPT to check the correctness of some proofs and to improve the presentation of the manuscript. All proof strategies and mathematical arguments were developed by the authors. The authors take full responsibility for the entire content of the paper.
	\\[-0.4cm]
	
	\noindent\textbf{Acknowledgment.} The authors would like to thank Sen Cao from Lanzhou University and Yuanhang Chang from the University of Science and Technology of China for their helpful discussions and valuable suggestions on this paper.
	\\[-0.4cm]
	
	\noindent\textbf{Funding.} This work is supported by National Natural Science Foundation of China (Grant No. 11671181,12371402) and Natural Science Foundation of Gansu Province (Grant No. 21JR7RA535, 24JRRA414).
	\\[-0.4cm]
	
	\noindent \textbf{Conflict of interest.}
	On behalf of all authors, the corresponding author states that there is no conflict of interest.
	\\[-0.4cm]
	
	\noindent \textbf{Authors' contributions.}
	All authors have contributed equally.
	\\[-0.4cm]
	
	\noindent \textbf{Data Availability.}
	This study is theoretical, and no empirical data were collected. The numerical values in the tables and figures were computed from the polynomials defined in the article.

	\vspace{0.5cm}
	
	\noindent
	{Jianlun Liu}\\
	School of Mathematics and Statistics, Lanzhou University,
	Lanzhou 730000, China\\
	\textit{E-mail address:} \texttt{jianlunliumath@163.com}
	
	\vspace{0.4cm}
	
	\noindent
	{Hong-Rui Sun}\\
	School of Mathematics and Statistics, Lanzhou University,
	Lanzhou 730000, China\\
	\textit{E-mail address:} \texttt{hrsun@lzu.edu.cn}
	
	\vspace{0.4cm}
	
	\noindent
	{Ziheng Zhang}\\
	School of Mathematical Sciences, Tiangong University,
	Tianjin 300387, China\\
	\textit{E-mail address:} \texttt{zhzh@mail.bnu.edu.cn}


\begin{thebibliography}{99}
		\bibitem{Branson1995}
		T.P. Branson,
		\newblock Sharp inequalities, the functional determinant, and the complementary series,
		\newblock \emph{Trans. Amer. Math. Soc.,} \textbf{347}(1995), 3671-3742.
		
		
		\bibitem{Branson-Orsted1991}
		T. P. Branson and B. {\O}rsted,
		\newblock Explicit functional determinants in four dimensions,
		\newblock \emph{Proc. Amer. Math. Soc.,} \textbf{113}(1991), 669-682.
		
		
		\bibitem{Caffarelli-Gidas-Spruck1989}
		L. A. Caffarelli, B. Gidas and J. Spruck,
		\newblock Asymptotic symmetry and local behavior of semilinear elliptic equations with critical Sobolev growth,
		\newblock \emph{Comm. Pure Appl. Math.,} \textbf{42}(1989), 271-297.
		
		
		\bibitem{Case-2024}
		J.S. Case,
		\newblock The Obata-V\'etois argument and its applications,
		\newblock \emph{J. reine angew. Math.,} \textbf{815}(2024), 23-40.
		
		
		\bibitem{Case-Gover2026}
		J.S. Case and A.R. Gover,
		\newblock The GJMS operators in geometry, analysis and physics,
		\newblock \emph{J. London Math. Soc. (2),} \textbf{113}(2026), e70375.
		
		
		\bibitem{Fazly-Wei-Xu2015}
		M. Fazly, J. Wei and X. Xu,
		\newblock A pointwise inequality for the fourth-order Lane-Emden equation,
		\newblock \emph{Anal. PDE,} \textbf{8}(2015), 1541-1563.
		
		
		\bibitem{Fefferman-Graham1985}
		C. Fefferman and C.R. Graham,
		\newblock Conformal invariants,
		\newblock in \'Elie Cartan et les math\'ematiques d'aujourd'hui
		(Lyon, 25-29 juin 1984),
		\emph{Ast\'erisque,} \textbf{S131}(1985), 95-116.
		
		
		\bibitem{Fefferman-Graham2012}
		C. Fefferman and C.R. Graham,
		\newblock The ambient metric,
		\newblock \emph{Annals of Mathematics Studies,}
		vol. \textbf{178}, Princeton University Press, Princeton, NJ, 2012.
		
		
		\bibitem{Gong-Kim-Wei2025}
		L. Gong, S. Kim and J. Wei,
		\newblock Compactness and non-compactness theorems of the fourth-order and sixth-order constant \(Q\)-curvature problems,
		\newblock arXiv:2502.14237, 2025.
		
		
		\bibitem{GJMS-1992}
		C.R. Graham, R. Jenne, L. Mason and G. Sparling,
		\newblock Conformally invariant powers of the Laplacian. I. Existence,
		\newblock \emph{J. London Math. Soc.,} \textbf{46}(1992), 557-565.
		
		
		\bibitem{Gursky-Malchiodi2015}
		M.J. Gursky and A. Malchiodi,
		\newblock A strong maximum principle for the Paneitz operator and a non-local flow for the \(Q\)-curvature,
		\newblock \emph{J. Eur. Math. Soc.,} \textbf{17}(2015), 2137-2173.
		
		
		\bibitem{Hang-Yang2015}
		F. Hang and P.C. Yang,
		\newblock Sign of Green's function of Paneitz operators and the \(Q\) curvature,
		\newblock \emph{Int. Math. Res. Not.,} \textbf{2015}(2015), 9775-9791.
		
		
		\bibitem{Hang-Yang2016}
		F. Hang and P.C. Yang,
		\newblock \(Q\)-curvature on a class of manifolds with dimension at least 5,
		\newblock \emph{Comm. Pure Appl. Math.,} \textbf{69}(2016), 1452-1491.
		
		
		\bibitem{Li-Wei2025}
		M. Li and J. Wei,
		\newblock A remark on the Case-Gursky-V\'etois identity and its applications,
		\newblock \emph{Proc. Amer. Math. Soc.,} \textbf{153}(2025), 3417-3430.
		
		
		\bibitem{Li-Wei-Xu2025}
		M. Li, J. Wei, and X. Xu,
		\newblock On geometry of \(Q_g^{(2k)}\)-curvature,
		\newblock arXiv:2506.20165, 2025.
		
		
		\bibitem{Li-Xu2025}
		M. Li and X. Xu,
		\newblock On positivity of the \(Q\)-curvatures of conformal metrics,
		\newblock \emph{J. Funct. Anal.,} \textbf{289}(2025), 111011.
		
		
		\bibitem{Lin1998}
		C.S. Lin,
		\newblock A classification of solutions of a conformally invariant fourth order equation in \(\mathbb R^n\),
		\newblock \emph{Comment. Math. Helv.,} \textbf{73}(1998), 206-231.
		
		
		\bibitem{Liu-Lyu-Sun-Zhang2026}
		J. Liu, H. Lyu, H.R. Sun and Z. Zhang,
		\newblock Quantitative positivity of the fourth order Q-curvature and applications to higher order elliptic equations,
		\newblock \emph{Proc. Amer. Math. Soc.,} submitted.
		
		
		\bibitem{Liu-Radulescu-Sun-Zhang2026}
		J. Liu, V. D. R\u{a}dulescu, H.R. Sun and Z. Zhang,
		\newblock Spectral characterization, local rigidity and Morse index of bubble solutions to the CR fractional Yamabe problem,
		\newblock arXiv:2610.39607, 2026.
		
		
		\bibitem{Martinazzi2009}
		L. Martinazzi,
		\newblock Classification of solutions to the higher order Liouville's equation on
		\(\mathbb R^{2m}\),
		\newblock \emph{Math. Z.,} \textbf{263}(2009), 307-329.
		
		
		\bibitem{Mazumdar-Premoselli2025}
		S. Mazumdar and B. Premoselli,
		\newblock Compactness of conformal metrics with constant \(Q\)-curvature of higher order,
		\newblock arXiv:2510.00888, 2025.
		
		
		\bibitem{Mazumdar-Vetois2024}
		S. Mazumdar and J. V\'etois,
		\newblock Existence results for the higher-order \(Q\)-curvature equation,
		\newblock \emph{Calc. Var. Partial Differential Equations,} \textbf{63}(2024), 151.
		
		
		\bibitem{Ngo-Ye2022}
		Q. Ng\^o and D. Ye,
		\newblock Existence and non-existence results for the higher order Hardy-H\'enon equations revisited,
		\newblock \emph{J. Math. Pures Appl.,} \textbf{163}(2022), 265-298.
		
		
		
		\bibitem{Obata-1971}
		M. Obata,
		\newblock The conjectures on conformal transformations of Riemannian manifolds,
		\newblock \emph{J. Differential Geometry,} \textbf{6}(1971), 247-258.
		
		
		
		\bibitem{Steele-2004}
		J.M. Steele,
		\newblock The Cauchy-Schwarz Master Class: An Introduction to the Art of Mathematical Inequalities,
		\newblock Cambridge: Cambridge University Press, 2004.
		
		
		
		
		\bibitem{Vetois2024}
		J. V\'etois,
		\newblock Uniqueness of conformal metrics with constant \(Q\)-curvature on closed Einstein manifolds,
		\newblock \emph{Potential Anal.,} \textbf{61}(2024), 485-500.
		
		
		\bibitem{Vetois-Zeitler-2026}
		J. V\'etois and S. Zeitler,
		\newblock Positivity and non-positivity results for the sixth-order
		\(Q\)-curvature of conformal metrics in \(\mathbb R^n\),
		\newblock arXiv:2607.18205, 2026.
		
		
		\bibitem{Wei-Xu1999}
		J. Wei and X. Xu,
		\newblock Classification of solutions of higher order conformally invariant equations, 
		\newblock \emph{Math. Ann.,} \textbf{313}(1999), 207-228.
		
	\end{thebibliography}
\end{document}